\documentclass[12pt,a4paper]{amsart}
\usepackage[a4paper,textwidth=468pt,hcentering]{geometry}

\usepackage{bm}
\usepackage{amsmath}
\usepackage{amssymb}
\usepackage{amsthm}
\usepackage{graphicx}
\usepackage{parallel}
\usepackage{ stmaryrd }
\usepackage{tikz}
\usepackage{multicol}
\usepackage{ stmaryrd }
\usepackage{ latexsym }
\usepackage[shortlabels]{enumitem}
\setlist[enumerate,1]{leftmargin=30pt,labelsep=6pt}
\setlist[itemize,1]{leftmargin=30pt,labelsep=6pt}
\usepackage{hyperref}

\newcommand{\forkindep}[1][]{
  \mathrel{
    \mathop{
      \vcenter{
        \hbox{\oalign{\noalign{\kern-.3ex}\hfil$\vert$\hfil\cr
              \noalign{\kern-.7ex}
              $\smile$\cr\noalign{\kern-.3ex}}}
      }
    }\displaylimits_{#1}
  }
}

\newcommand{\mipo}{\operatorname{MiPo}}
\newcommand{\Char}{\operatorname{Char}}
\newcommand{\qf}{\operatorname{qf}}
\newcommand{\NN}{\mathbb{N}}

\newcommand{\QQ}{\mathbb{Q}}
\newcommand{\RR}{\mathbb{R}}
\newcommand{\ZZ}{\mathbb{Z}}
\newcommand{\VV}{\mathbb{V}}
\newcommand{\PP}{\mathbb{P}}

\newcommand{\FF}{\mathbb{F}}
\newcommand{\QQp}[1]{\mathbb{Q}[X]_{\operatorname{irr}}^{#1}}
\newcommand{\Kp}[1]{K[X]_{\operatorname{irr}}^{#1}}
\newcommand{\Cc}{{\mathcal{C}}}
\newcommand{\Ccalg}{{\mathcal{C}^{\operatorname{alg}}}}
\newcommand{\Cctrans}{{\mathcal{C}^{\operatorname{tr}}}}
\newcommand{\Ff}{{\mathcal{F}}}
\newcommand{\mm}{\mathcal{M}}
\newcommand{\nn}{\mathcal{N}}

\newcommand{\ii}{\mathcal{I}}
\newcommand{\jj}{\mathcal{J}}
\newcommand{\kk}{\mathcal{K}}
\newcommand{\rr}{\mathcal{R}}

\newcommand{\Hh}{\mathcal{H}}

\newcommand{\stri}{{\scalebox{0.5}{$\triangle$}}}

\newcommand{\Tor}{\operatorname{Tor}}

\newcommand{\set}[1]{{\{#1\}}}
\newcommand{\Set}[1]{{\left\{#1\right\}}}
\newcommand{\spanA}[2]{{{\langle#1\rangle_{#2}}}}

\newcommand{\dotcup}{\mathbin{\dot{\cup}}}
\newcommand{\Fac}{{\operatorname{Fac}}}
\newcommand{\Ker}{{\operatorname{Ker}}}
\newcommand{\ACF}{{\operatorname{ACF}}}
\newcommand{\Image}{{\operatorname{Im}}}

\newcommand{\rk}{{\operatorname{rk}}}

\newcommand{\Diag}{{\operatorname{Diag}}}
\newcommand{\Id}{{\operatorname{Id}}}

\newcommand{\lcm}{{\operatorname{lcm}}}
\newcommand{\Lex}{{\operatorname{Lex}}}

\newcommand{\ud}{{\underline{d}}{}}

\newcommand{\ut}{{\underline{t}}{}}
\newcommand{\uu}{{\underline{u}}{}}
\newcommand{\uv}{{\underline{v}}{}}
\newcommand{\uw}{{\underline{w}}{}}
\newcommand{\ux}{{\underline{x}}{}}
\newcommand{\uy}{{\underline{y}}{}}
\newcommand{\uz}{{\underline{z}}{}}

\newcommand{\uzero}{{\underline{0}}}

\newcommand{\umu}{{\underline{\mu}}}

\newcommand{\utau}{{\underline{\tau}}}

\newcommand{\tiluy}{{\underline{\Tilde{y}}}{}}
\newcommand{\tiluw}{{\underline{\Tilde{w}}}{}}
\newcommand{\li}{{\scalebox{0.5}{{$\operatorname{li}$}}}}
\newcommand{\ld}{{\scalebox{0.5}{{$\operatorname{ld}$}}}}

\newcommand{\lii}{{\scalebox{0.5}{$\operatorname{li}$}}}
\newcommand{\ldd}{{\scalebox{0.5}{$\operatorname{ld}$}}}

\newcommand{\uvvec}{{\underline{\Vec{v}}}}

\newcommand{\uxvec}{{\underline{\Vec{x}}}}

\newcommand{\xvec}{{\Vec{x}}}

\newcommand{\Hfour}{{$(\operatorname{H4})$}}

\newcommand{\TQvs}{{T_{\QQ\operatorname{-vs}}}}

\newcommand{\TKvs}{{T_{K\operatorname{-vs}}}}
\newcommand{\TKvsThe}{{T_{K\operatorname{-vs},\theta}}}
\newcommand{\TKvsTheC}{{T^C_{K\operatorname{-vs},\theta}}}

\newcommand{\RCF}{{\operatorname{RCF}}}

\newcommand{\LK}{{L_K}}
\newcommand{\NIP}{$\operatorname{NIP}$}

\newcommand{\TPtwo}{$\operatorname{TP}_2$}
\newcommand{\NATP}{$\operatorname{NATP}$}

\newcommand{\SOP}{$\operatorname{SOP}$}
\newcommand{\dprk}{\operatorname{dp-rk}}

\newcommand{\LKThe}{{L_{K,\theta}}}
\newcommand{\LRC}{{L_{R_C}}}

\newcommand{\Lr}{L_{\operatorname{r}}}

\newcommand{\standartConstruction}{\hyperref[lemma_standart_construction]{Standard Construction}}
\newcommand{\remainderRule}{\hyperref[lemma_remainder_rule]{Euclidean Division}}
\newcommand{\TrafoLemma}{\hyperref[lemma_trafoo]{Transformation Lemma}}
\newcommand{\remainderRules}
{\hyperref[lemma_remainder_rule]{Euclidean Divisions}}
\newcommand{\projectionRule}{\hyperref[lemma_projection_rule]{Projection Identity}}
\newcommand{\placeholderNotation}{\hyperref[def_placeholder_notation]{Placeholder-Notation}}

\newcommand{\commonBaseTheorem}{\hyperref[theorem_r_c_element]{Common Base Theorem}}

\newcommand{\gcdRule}{\hyperref[lemma_gcd_rule]{Bézout's Identity}}

\newcommand{\sminus}{\text{-}}
\newcommand{\sindex}[1]{{}_{\scalebox{0.45}{$#1$}}}
\newcommand{\drawTextHelper}[5]{
\node[anchor=center, scale = 1] at (#1 * #4 + 0.5 * #4, -#2 * #5 - 0.5 * #5) {$#3$};
}
\newcommand{\drawText}[3]{\drawTextHelper{#1}{#2}{#3}{1.0}{0.6}}

\newcommand{\drawBorderHelper}[6]{\draw [black, line width=0.75]
(#1 * #5 + 0.105, -#2 * #6) --
(#1 * #5, -#2 * #6) --
(#1 * #5, -#4 * #6) --
(#1 * #5 + 0.105, -#4 * #6)
(#3 * #5 - 0.105, -#2 * #6) --
(#3 * #5, -#2 * #6) --
(#3 * #5, -#4 * #6) --
(#3 * #5 - 0.105, -#4 * #6);
\node[anchor=east, scale = 1.0] (Frame) at (0.15, -0.5 * #2 * #6 + -0.5 * #4 * #6 -0.08) {};}
\newcommand{\drawBorder}[4]{\drawBorderHelper{#1}{#2}{#3}{#4}{1.0}{0.6}}

\newcommand{\drawHDotsHelper}[5]{
\draw [line width=1.0, line cap=round, dash pattern=on 0 off 9.5 * #4]
        (#1 * #4, -#2 * #5 - 0.5 * #5) --
        (#1 * #4 + #3 * #4 + 0.02 * #4, -#2 * #5 - 0.5 * #5);    
}
\newcommand{\drawHDots}[3]{\drawHDotsHelper{#1}{#2}{#3}{1.0}{0.6}}

\newcommand{\drawVDotsHelper}[5]{
\draw [line width=1.0, line cap=round, dash pattern=on 0pt off 9.5 * #5]
        (#1 * #4 + 0.5 * #4, -#2 * #5) --
        (#1 * #4 + 0.5 * #4, -#2 * #5 -#3 * #5 - 0.02 * #5);    
}
\newcommand{\drawVDots}[3]{\drawVDotsHelper{#1}{#2}{#3}{1.0}{0.6}}

\newcommand{\drawDDotsHelper}[5]{
\pgfmathsetmacro{\factor}{sqrt(#4^2 + #5^2)}
\draw [line width=1.0, line cap=round, dash pattern=on 0pt off \factor * 9.5]
        (#1 * #4, -#2 * #5) --
        (#1 * #4 + #3 * #4 + 0.02 * #4, -#2 * #5 -#3 * #5 - 0.02 * #5);    
}
\newcommand{\drawDDots}[3]{\drawDDotsHelper{#1}{#2}{#3}{1.0}{0.6}}

\newcommand{\drawHLineHelper}[5]{
\draw [line width=0.25]
        (#1 * #4 + 0.125 * #4, -#2 * #5) --
        (#1 * #4 + #3 * #4 - 0.125 * #4, -#2 * #5);    
}
\newcommand{\drawHLine}[3]{\drawHLineHelper{#1}{#2}{#3}{1.0}{0.6}}

\newcommand{\drawVLineHelper}[5]{
\draw [line width=0.25]
        (#1 * #4, -#2 * #5 -0.125 * #5) --
        (#1 * #4, -#2 * #5 -#3 * #5 + 0.125 * #5);    
}
\newcommand{\drawVLine}[3]{\drawVLineHelper{#1}{#2}{#3}{1.0}{0.6}}

\newtheorem*{notation}{Notation}

\newtheorem{theorem}{Theorem}[section]
\newtheorem*{theorem*}{Theorem}
\newtheorem{theoremi}{Theorem}

\newtheorem{definition}[theorem]{Definition}
\newtheorem{fact}[theorem]{Fact}
\newtheorem{remark}[theorem]{Remark}
\newtheorem{lemma}[theorem]{Lemma}
\newtheorem{corollary}[theorem]{Corollary}
\newtheorem{example}[theorem]{Example}
\newtheorem{observation}[theorem]{Observation}

\newtheorem{subclaim}{Claim}[theorem]
\newtheorem{subdefinition}[subclaim]{Definition}

\newtheorem*{subdefinition*}{Definition}
\newtheorem*{subclaim*}{Claim}

\newenvironment{innerproof}[1][Proof]
{\begin{proof}[#1]}
{\end{proof}}

\title{Model Theory of Generic Vector Space Endomorphisms}

\author{Leon Chini}
\address{Mathematisches Institut, Universität Bonn\\
Endenicher Allee 60\\
D-53115 Bonn, Germany}
\email{lchini@uni-bonn.de}
\urladdr{https://leonchini.github.io}

\begin{document}

\begin{abstract}
This paper deals with the model companion of an endomorphism acting on a vector space, possibly with extra structure.
Given a theory $T$ that \hbox{$\varnothing$-defines} an infinite $K$-vector space ${\mathbb{V}}$ in every model, we define 
$T_\theta := T \cup \{\text{``$\theta$ defines a $K$-endomorphism of $\mathbb{V}$''}\}$.
We then consider extensions of the form
$$
    T_\theta \cup \big\{\sum\nolimits_{k}\bigcap\nolimits_{l}\operatorname{Ker}(\rho_{j, k, l}[\theta]) \!=\! \sum\nolimits_{k}\bigcap\nolimits_{l} \operatorname{Ker}(\eta_{j, k, l}[\theta]) : j \!\in\! \mathcal{J}\big\},
$$
where all sums and intersections are finite, 
and all the $\rho[\theta]$'s and $\eta[\theta]$'s are polynomials over $K$ with $\theta$ plugged in.
Note that properties such as $\theta^2 - 2\operatorname{Id} = 0$ or $\operatorname{Ker}(\theta^n) = \operatorname{Ker}(\theta^{n+1})$ can be expressed in such a form.
We then parametrize the consistent extensions of this form by a family $\{T^C_\theta : C \in \mathcal{C}\}$ and characterize the existentially closed models of each $T^C_\theta$.
We also present a sufficient criterion, which depends only on $T$, for when these characterizations are first-order expressible, i.e., for when a model companion of each $T^C_\theta$ exists.
\end{abstract}

\maketitle
\tableofcontents


\section{Introduction}
\noindent This paper deals with the model companion of an endomorphism acting on a vector space, possibly with extra structure.

A structure $\mm$ is \textit{existentially closed in an extension} $\nn \supseteq \mm$ if every existential sentence with parameters in $\mm$ that is true in $\nn$ is also true in $\mm$.
We say that $\mm$ is an \textit{existentially closed model} of a theory $T$ if $\mm$ is existentially closed in any extension $\nn \supseteq \mm$ that is also a model of $T$.
In practice, this means that an existentially closed model $\mm$ satisfies a form of ``Nullstellensatz'', in the sense that any system of equations with parameters in $\mm$ has a solution in $\mm$, as long as the existence of a solution does not contradict the theory $T$.
In the case where the models of $T$ are closed under chains of substructures ($T$ is called \textit{inductive}) and the existentially closed models form an elementary class, we call the theory of this class the \textit{model companion} of $T$.

When dealing with an endomorphism $\theta$ acting on a vector space over a field $K$, one should notice that $\theta$ can be plugged into any polynomial $P$ over $K$ to obtain another endomorphism $P[\theta]$ that acts on the same vector space.
We can express several interesting properties of $\theta$ in a first-order way, such as $\theta^2 - 2\Id = 0$ or $\Ker(P[\theta]) = \set{0}$ for any polynomial $P$ with positive degree.
The idea here is to include (any reasonable combination of) such properties as a set of constraints imposed on a given theory of vector spaces expanded by an endomorphism, study the existentially closed models, and then give a criterion for when these form an elementary class.
Given a real closed field, the positive elements with multiplication form a divisible torsion-free abelian group (i.e., a $\QQ$-vector space), and our results apply in this case.
We are able to give a description of existentially closed real closed fields expanded by an endomorphism of the group of positive elements with the extra properties mentioned above ($\theta^2 - 2 \Id = 0$, $\Ker(P[\theta]) = \set{0}$, etc.).

Our results constitute a generalization of the work of d'Elbée in \cite{dEl25}, where the model companion of algebraically closed fields with a multiplicative endomorphism, called $\operatorname{ACFH}$, is studied.
Modulo the torsion elements, the multiplicative group can be seen as a $\QQ$-vector space; hence it is very close to our setting (see Remark \ref{rk:ACFcase}).
Note that the case studied by d'Elbée corresponds to an endomorphism with no constraints whatsoever on the properties of the endomorphism.
As mentioned by d'Elbée \cite{dEl25}, in the case of an endomorphism of the additive group of an algebraically closed field, the model companion exists if and only if the characteristic of the field is positive.
This work follows a prominent line of research in model theory: generic expansions of algebraically closed fields.
This includes classical examples such as $\operatorname{ACFA}$, the model companion of difference fields, defined by Macintyre \cite{Mac97}, or the theory of differentially closed fields $\operatorname{DCF}$, defined by Robinson in \cite{Rob59}.
For more recent examples, we already mentioned $\operatorname{ACFH}$, but one could also consider the generic exponential map, described in positive logic by Haykazyan and Kirby in \cite{HK21}, or d'Elbée's ``generic expansion by a reduct'' \cite{dEl21b}, which can be used to extend $\operatorname{ACF}$ with a generic predicate for a multiplicative or (in positive characteristic) additive subgroup.

The study of existentially closed models has led to many interactions between model theory and algebra.
Important applications of model theory to Diophantine geometry are Hrushovski's proofs of the Mordell-Lang conjecture \cite{Hru96} and the Manin-Mumford conjecture \cite{Hru01}.
Those striking results are obtained by applying local stability methods within the framework of a certain model-complete theory, which serves as an ambient universal domain.
The proof of the Mordell-Lang conjecture uses $\operatorname{DCF}$, whereas the proof of the Manin-Mumford conjecture occurs within the framework of $\operatorname{ACFA}$.
The model theory of differentially closed fields $\operatorname{DCF}$ seems to be an endless source of striking applications to differential or diophantine geometry; see, e.g., \cite{FJM22} by Freitag, Jaoui, and Moosa or \cite{FS18} by Freitag and Scanlon.

The model theory of modules has been well studied, with foundational work by Ziegler \cite{Zie84} and further development thereafter; see, e.g., the two books \cite{Pre88} and \cite{Pre09} by Prest for an overview.
As the structure of a $K$-vector space with an endomorphism is exactly that of a $K[X]$-module, this area of research is highly relevant to our setting.
It is a well-known result of Baur \cite{Bau76} that in a complete theory of modules over a fixed ring $R$, every formula is equivalent to a Boolean combination of pp-formulas.
If $R$ is a principal ideal domain, these pp-formulas can furthermore be chosen in a particularly nice form.
This can essentially be proved using the existence of a Smith normal form for matrices over a principal ideal domain; see, e.g., Theorem 2.$\ZZ$1 in \cite{Pre88} (note that the proof presented there is attributed to Hodges and the result itself to Eklof and Fisher).
In the proof of our characterization of existentially closed models, a crucial step is to simplify certain existential formulas containing conjunctions of eqautions in the language of a $K$-vector space with an endomorphism.
This simplification also relies on the existence of the Smith normal form for matrices over principal ideal domains.
Thus, in a very loose sense, we generalize this simplification result for pp-formulas to the structure outside the module.
We can also associate the constraints we consider on our endomorphism $\theta$ with certain closed subsets of the Ziegler spectrum (see Remark \ref{remark_ziegler_spectrum}).

We also view this work as a new point of contact between model theory and classical results about endomorphisms of finite-dimensional vector spaces.
Under a ``generalized surjectivity'' condition (which we call \textit{$C$-image-completeness}; see Definition \ref{def_C_image_comple}), which always holds for existentially closed models, we can definably decompose the underlying vector space into direct summands (that depend on the constraints chosen for $\theta$) on which $\theta$ behaves well.
This can be thought of as a decomposition into algebraic subspaces (i.e., subspaces defined by $P[\theta](x) = 0$ for some polynomial $P$ over the field $K$) and a potential ``transcendental'' extra subspace.
Model theory has seen a lot of interest in the study of vector spaces expanded by bilinear forms, which also established new connections between model theory and linear algebra. It started with Granger's thesis \cite{Gra99}, was then rekindled by Harrison-Shermoen \cite{HS13}, and popularized by Chernikov and Ramsey in \cite{CR16}, with a considerable amount of interest thereafter \cite{Dob23}, \cite{KR24a}, \cite{CH16}, \cite{CH21}, and \cite{CH24}. In the same vein, the model theory of Lie algebras has been studied in \cite{dMRS25a} and \cite{dMRS25b} by d'Elbée, Müller, Ramsey, and Siniora.

Model-theoretic bilinear algebra is quite successful in providing new examples to feed Shelah's classification theory; e.g., \cite{CH24} provides examples of strictly $n$-dependent NSOP$_1$ theories. We expect the same phenomenon to occur here. Depending on the starting theory of vector spaces and the given constraints on $\theta$, the model companion will witness various neostability properties.
For example, if one starts with the pure theory of $K$-vector spaces, then depending on the constraints on $\theta$, the model companion will either be strongly minimal, totally transcendental, or strictly stable. If one starts with the theory of ordered divisible abelian groups and $K = \QQ$, the resulting model companion will be \NIP{}, and the $\dprk{}$ will depend on the constraints imposed on $\theta$. Finally, if one takes $\RCF{}$, the theory of real closed fields, and chooses the positive elements with multiplication as a $\QQ$-vector space, the resulting model companion will be \TPtwo{}, \SOP{}, and \NATP{}.
We expect this to generalize to the following trichotomy in the $o$-minimal setting: the model companion either  
\begin{itemize}
    \item does not exist;  
    \item is distal (and hence \NIP{}) and there is a precise formula for the $\dprk{}$ that depends on the set of constraints imposed on $\theta$ and on the definable germs of endomorphisms;  
    \item is \TPtwo{}, \SOP{}, and \NATP{}.  
\end{itemize}
In the last two cases, the model companion will also have an $o$-minimal open core. If the starting theory is o-minimal, the vector space is a subset of the line, and its operations are continuous, then the class of existentially closed models is elementary under exactly the same condition as in \cite{Blo23}, in which Block Gorman gave a criterion for the existence of a model companion of an $o$-minimal theory with a predicate for a subgroup. That all these theories are \NATP{} will follow from a general preservation of \NATP{} result for our construction.
Tameness of the model companion, the $o$-minimal case, a refined quantifier elimination result, the preservation of \NATP{}, and the study of certain reducts of the model companion will appear in forthcoming papers.

We now give more precise statements of our results. Let $L$ be a language and $T$ a model-complete $L$-theory with an infinite $\varnothing$-definable vector space $\VV$ in every model. Given a consistent set $C$ of constraints on an endomorphism, which encodes conditions of the form
$$
    \sum\nolimits_{k}\bigcap\nolimits_{l}\Ker(\rho_{k, l}[\theta]) = \sum\nolimits_{k}\bigcap\nolimits_{l} \Ker(\eta_{k, l}[\theta])
$$
(where all sums and intersections are finite,  
and all the $\rho$'s and $\eta$'s are polynomials over $K$), we define the following theory in the language $L \cup \set{\theta}$:
\begin{align*}
    T^C_\theta := T &\cup \set{\text{``$\theta$ is an endomorphism of $\VV$''}} \\ &\cup \set{\text{``$\theta$ satisfies the constraints in $C$''}}.
\end{align*} 
We will see in Section \ref{sec_baisc_stuff} that any such set of constraints $C$ is either \textit{algebraic}, i.e., equivalent to $\set{\Ker(P[\theta]) = \VV }$ for some $P \in K[X] \setminus \set{0}$, or \textit{transcendental} in the sense of being non-algebraic.  
If such a set $C$ is transcendental, we show that it is equivalent to a set of conditions of the form  
$$
\set{\Ker(f^{n_f}[\theta]) = \Ker(f^{n_f+1}[\theta]) : f \in Q}
$$  
for some set $Q$ of irreducible monic polynomials and a sequence $(n_f : f \in Q)$ of natural numbers.  
Note that the property $\Ker(f^{0}[\theta]) = \Ker(f^{0+1}[\theta])$ is preceisly the injectivity of $f[\theta]$, i.e., $\Ker(f[\theta]) = \set{0}$.  
We can already formulate our characterizations of existentially closed models in two cases:

\begin{theoremi}[Theorem \ref{theorem_C0_char}] \label{theorem_char_A}
    Let $C := \set{\Ker(P[\theta]) = \set{0} : P \in K[X] \!\setminus\! \set{0}}$. A model $(\mm, \theta)$ of $T^C_\theta$ is existentially closed if and only if:
    \begin{enumerate}[(i)]
        \item The map $P[\theta]$ is an automorphism of $\VV$ for every $P\in K[X] \setminus \set{0}$;
        \item $\exists \ux \in \VV : \psi(\theta^0(\ux), \theta^1(\ux), \dots)$ holds for every $L(M)$-formula $\psi(\ux_0, \ux_1, \dots)$ that implies no finite disjunction of non-trivial linear dependencies over $\VV$ (and satisfies $|\ux_i| = |\ux_0|$ for all $i$).
    \end{enumerate}
\end{theoremi}

\begin{theoremi}[Theorem \ref{theorem_Calgirr_char}] \label{theorem_char_B}
    Let $f \in K[X]$ be an irreducible polynomial and assume that $C = \set{\Ker(f[\theta]) = \VV}$. A model $(\mm, \theta)$ of $ T^C_\theta$ is existentially closed if and only if:
    \begin{itemize}
        \item[] $\exists \ux \in \VV : \psi(\theta^0(\ux), \dots, \theta^{\deg(f)-1}(\ux))$ holds for every $L(M)$-formula $\psi(\ux_0,\dots, \ux_{\deg(f)-1})$ that implies no finite disjunction of non-trivial linear dependencies over $\VV$.
    \end{itemize}
\end{theoremi}

\noindent For general sets of constraints $C$, the characterization of existentially closed models is more involved. It is given in Theorem \ref{theorem_big_characterization} and roughly consists of the following:
\begin{enumerate}[(i)]
    \item One first needs to impose the above-mentioned ``generalized surjectivity'' condition, called \textit{$C$-image-completeness} (see Definition \ref{def_C_image_comple}), which depends on $C$. In Theorem \ref{theorem_char_A}, this condition simplifies to point (i), and in Theorem \ref{theorem_char_B}, it can be omitted since it trivially holds for algebraic sets of constraints.
    \item One needs to ask that the sentences $\exists \ux \in \VV : \psi(\theta^0(\ux), \theta^1(\ux), \dots) \wedge S(\ux)$ hold for all \textit{$C$-sequence-systems} $S(\ux)$ and formulas $\psi(\ux_0, \ux_1, \dots)$ that imply no finite disjunction of non-trivial linear dependencies over $\VV$ and are bounded by $S$.
    A $C$-sequence-system is a special conjunction of $L_{K\text{-vs}} \cup \set{\theta}$ equations depending on $C$ (see Definition \ref{def_c_sequence_system}), and $\psi(\ux_0, \ux_1, \dots)$ being bounded by $S(\ux)$ means that $S$ dictates which variables $x_{i, k}$ are actually allowed to appear in $\psi(\ux_0, \ux_1, \dots)$ (see Definition \ref{def_formual_bounded}).
    With $C$ as in Theorem \ref{theorem_char_A}, these conjunctions are always empty by definition, and in Theorem \ref{theorem_char_B}, they are already implied by $T^C_\theta$.
\end{enumerate}
We present a criterion that implies that these characterizations are first-order axiomatizable, i.e., that model companion of $T^C_\theta$ exists for each (reasonable) $C$:  

\begin{theoremi}[Theorem \ref{theorem_first_oder}]
    If for every $L$-formula $\psi(\ux; \uw)$ there is an $L$-formula $\sigma_\psi(\uw)$ such that the following holds for all $\mm \models T$ and $\ud \in M$:
    $$
    \mm \models \sigma_\psi(\ud) \;\; \Leftrightarrow \;\; \parbox{8.5cm}{``There are $\mm' \succ \mm$ and a tuple $\uv' \in \VV'$  linearly\\ ${}$\hspace{5pt} independent over $\VV$ such that $\mm' \models \psi(\uv'; \ud)$'',}
    $$
    then the model companion of $T^C_\theta$ exists for every (reasonable) set of constraints $C$.
\end{theoremi}
\noindent This criterion corresponds to condition $(\operatorname{H}_4)$ from Definition 1.10 in \cite{dEl21b} (one should notice that the other conditions $(\operatorname{H}_1)$, $(\operatorname{H}_2)$, and $(\operatorname{H}^+_3)$ already hold in our setting).
This condition is used to express in a first-order way the statement ``implies no finite disjunction of non-trivial linear dependencies over $\VV$'' that appears in our characterizations.

As mentioned above, for any polynomial $P \in K[X]$, we obtain an endomorphism $P[\theta]$ of $\VV$, which in turn induces a $K[X]$-module structure on $\VV$.
Now, assume, for instance, that $\theta$ satisfies the condition $P[\theta] = 0$ for a given $P \in K[X]$.
It is then more natural to treat $\VV$ as a module over the ring $K[X]/(P)$ rather than over $K[X]$.
Given a fixed set $C$ of constraints on $\theta$, this raises the following question:
\begin{itemize}
    \item[] What is the ``optimal'' ring $R_C$ for which there is an $L_{K\text{-vs}} \cup \set{\theta}$-definable $R_C$-module on every existentially closed model \hbox{$(\mm, \theta) \models T^C_\theta$}?
\end{itemize}
We define such a ring $R_C$ (see Theorem \ref{theorem_r_c_def}) for every (reasonable) set of constraints $C$ and show that, in an existentially closed model of $T^C_\theta$, the $R_C$-module structure is definable, and distinct elements $r \neq r' \in R_C$ define distinct endomorphisms (see Section \ref{sec_rc_module}).
The ``optimality'' of this ring will follow in an upcoming paper, where we use it to extend our language, giving us quantifier elimination in certain cases.
We also describe these rings in algebraic terms (see Corollary \ref{corollary_R_C_as_ring}); for instance, we obtain
\begin{itemize}
    \item $R_C = K[X]$ if $C$ is empty;
    \item $R_C = K(X)$ if $C = \set{\Ker(P[\theta]) = \set{0} : P \in K[X] \setminus \set{0}}$;
    \item $R_C = K[X]/(P)$ if $C = \set{\Ker(P[\theta]) = \VV}$ for some non-constant $P \in K[X]$.
\end{itemize}
In general, these rings turn out to be much more complicated, the reason being that projections to certain kernels or images become definable, and some maps of the form $P[\theta]$ can become invertible on certain images.

\section{Preliminary Results}

Let $L$ be a first-order language, let $T$ be a model-complete $L$-theory, and let $K$ be a field.
Furthermore, assume that the theory of $K$-vector spaces is definable in $T$.
By this we mean that there are $L$-formulas $\Omega_{\VV}(\ux)$, $\Omega_{0}(\ux)$, $\Omega_{+}(\ux_1, \ux_2, \uy)$, and $\Omega_{\lambda\cdot}(\ux, \uy)$ for each $\lambda \in K$ such that, in every model $\mm \models T$, they define an infinite $K$-vector space $(\VV, 0, +, (\lambda \cdot)_{\lambda \in K})^{\mm}$.
In the case $K = \QQ$, one may also view $(\VV, 0, +)^{\mm}$ as a torsion-free divisible abelian group, since these are precisely the $\QQ$-vector spaces.
If a model $\mm \models T$ is given, then $\VV$ denotes $\VV^\mm$.
If the model is denoted with $\mm'$ instead of $\mm$, then we will write $\VV'$ instead of $\VV$.
If no model of $T$ is clear from the context, then we will still use the letter $\VV$ to denote a $K$-vector space.
We may say ``vector space'' instead of ``$K$-vector space'', ``polynomial'' instead of ``$K$-polynomial'', ``linearly independent'' instead of ``$K$-linearly independent'', and so on.

\begin{example}
    \label{example_main}
    Two of our main examples are as follows:
    \begin{enumerate}[(i)]
        \item Let $\LK = \set{0, +, (\lambda \cdot)_{\lambda\in K}}$ and let $\TKvs$ be the theory of $K$-vector spaces, with the obvious formulas.
        Then, for any $\mm \models \TKvs$, we choose $(\VV, 0, +, (\lambda \cdot)_{\lambda\in K}) = \mm$.
        Similarly, we can also work with ordered divisible abelian groups, since these are precisely ordered $\QQ$-vector spaces.
        \item Let $K = \QQ$, $\Lr = \set{0, 1, +, \cdot, <}$, and $T = \RCF$, the theory of real closed fields.
        Then, for any $\rr \models \RCF{}$, we choose $(\VV, 0, +, (q \cdot)_{q\in \QQ})^{\rr}$ to be $(\rr_{>0}, 1, \cdot, (x \mapsto x^q)_{q\in \QQ})$.
    \end{enumerate}
\end{example}

\noindent Notationally, we will treat $\VV$ as a unary set, or even as a separate sort.
Any $\LK$-term or formula can then be viewed as an $L$-definable function or an $L$-formula, respectively.
Note that $0$, $+$, and $(\lambda\cdot)_{\lambda\in K}$ need not belong to $L$, so they are not necessarily $L$-terms.
For a given theory, there can be multiple definable vector spaces, as in the case of $\RCF{}$.
Thus, whenever a theory $T$ is given, we actually mean the tuple $(T, \Omega_{\VV}, \Omega_{0}, \Omega_{+}, (\Omega_{\lambda\cdot})_{\lambda \in K})$.

We now define $L_\theta := L(\theta) := L \cup \set{\theta}$, where $\theta$ is a function symbol not contained in $L$.
This has the drawback that, for example, if $\rr \models \RCF$ and the positive elements are viewed as a $\QQ$-vector space, then $\theta$ must also be defined outside $\VV = \rr_{>0}$.
In this case, one might try to set $\theta(0) = 0$ and $\theta(-1) \in \set{1, -1}$ so as to extend $\theta$ to an endomorphism of $(\rr, 1, \cdot)$, but then the ambient structure is no longer a vector space.
For the sake of uniformity, we instead define $\theta(x)$ to be the neutral element of $\VV$ for all $x \not\in \VV$ and set
$$
T_\theta := T \cup \set{\text{``$\theta_{\restriction \VV}$ is a $(\VV, 0, +, (\lambda \cdot)_{\lambda\in K})$-endomorphism''}} \cup \set{\forall x \not\in \VV: \theta(x) = 0}.
$$
In particular, $\theta^n(x) = 0$ for all $x \not\in \VV$ and all $n > 0$.
For consistency, we also define $\theta^0(x) = 0$ whenever $x \not\in \VV$, and $x + y = 0$ whenever either $x \not\in \VV$ or $y \not\in \VV$.
Practically, we will ignore the behavior of $\theta$ outside of $\VV$ and simply treat $\theta$ as a function defined only on $\VV$.
For example, we set $\Ker(\theta) := \set{v \in \VV : \theta(v) = 0}$, and similarly define the kernel for any function that is an endomorphism of $\VV$.

An alternative approach would be to take $\theta$ as a relation symbol defining the graph of an endomorphism of $\VV$.
For our purposes, the two approaches are equivalent.
The relation-symbol approach could also be used for interpretable vector spaces.
However, for some quantifier elimination results, we need $\theta$ to be a function symbol anyway.
Thus, it is also more natural for us to treat $\theta$ as a function symbol.

Note that if $\VV$ is an $n$-ary set, then we actually need $n$ different $n$-ary function symbols $\theta_1, \dots, \theta_n$ rather than a single unary function symbol $\theta$.
However, as mentioned above, we will treat elements of $\VV$ as singletons in order to simplify notation and will therefore pretend that $\theta$ is unary.

\subsection{Vector spaces with an endomorphism}
\label{sec_vec_with_endo}
\noindent For many arguments, it suffices to work with the $K$-vector space $\VV$ together with the endomorphism $\theta$, that is, with the pair $(\VV, \theta) := (\VV, 0, +, (\lambda\cdot)_{\lambda \in K}, \theta)$, rather than with $(\mm, \theta)$.
Thus, throughout this section, let $(\VV, \theta)$ consist of a $K$-vector space $\VV$ together with an endomorphism $\theta \colon \VV \to \VV$. We let $\LK$ and $\TKvs$ be, respectively, the language and theory of $K$-vector spaces. Consequently, $\LKThe$ and $\TKvsThe$ are, respectively, the language and theory of $K$-vector spaces with an endomorphism.

\begin{notation}
    Given a polynomial $\rho \in K[X]$, we write $(\rho)_i$ for the coefficient of $X^i$ in $\rho$, so for every $d \geq \deg(\rho)$, we have $\rho = \sum_{i=0}^d (\rho)_i \cdot X^i$.
\end{notation}

\noindent As a general convention, we set $\deg(0) := -\infty$. Also, we define a sum $\sum_{k= k_0}^m \dots$ to be empty, i.e., equal to $0$, whenever $m < k_0$.

\begin{definition}
    Given a polynomial $\rho \in K[X]$ and any $d \geq \deg(\rho)$, we let $\rho[\theta]$ denote the endomorphism of $\VV$ defined by
    $
    \rho[\theta](v) := \sum\nolimits_{i=0}^{d} (\rho)_i \cdot \theta^i(v)
    $
    for all $v \in \VV$.
\end{definition}

\begin{lemma}\label{lemma_add_mult_endo}
    Set $K[\theta] := \set{\rho[\theta] : \rho \in K[X]}$.
    Then $(K[\theta], 0[\theta], 1[\theta], +, \circ)$ is a ring satisfying
    \begin{enumerate}[(i)]
        \item $\rho_1[\theta] + \rho_2[\theta] = (\rho_1 + \rho_2)[\theta]$;
        \item $\rho_1[\theta] \circ \rho_2[\theta] = (\rho_1 \cdot \rho_2)[\theta]$;
    \end{enumerate}
    for all $\rho_1, \rho_2 \in K[X]$.
    The map
    $(K[X], 0, 1, +, \cdot) \to (K[\theta], 0[\theta], 1[\theta], +, \circ); \rho \mapsto \rho[\theta]$
    is a surjective ring homomorphism.
    In particular, $(K[\theta], 0[\theta], 1[\theta], +, \circ)$ is commutative.
\end{lemma}
\begin{notation}
    If $(\VV, \theta)$ is clear from the context, we write $\Ker(\rho)$ instead of $\Ker(\rho[\theta])$ and $\Image(\rho)$ instead of $\Image(\rho[\theta])$.
\end{notation}

\noindent The following fact and the corollary below will play an important role in this paper:

\begin{fact}
\label{fact_euclid_domain}
    $(K[X], 0, 1, +, \cdot, \deg)$ is a Euclidean domain.
    In particular, given $\rho_1, \rho_2 \in K[X]$, we can find $\chi_1, \chi_2 \in K[X]$ such that
    $$
    \chi_1 \cdot \rho_1 + \chi_2 \cdot \rho_2 = \gcd(\rho_1, \rho_2).
    $$
\end{fact}

\noindent Given two polynomials $\rho_1, \rho_2 \in K[X]$, we define the greatest common divisor $\gcd(\rho_1, \rho_2)$ to have leading coefficient $1$, unless $\rho_1 = \rho_2 = 0$, in which case we set $\gcd(0, 0) = 0$.
In particular, if $\rho_1 \neq 0$, then $\gcd(\rho_1, 0)$ is the unique monic scalar multiple of $\rho_1$. We also define $\gcd(\set{\rho_i : i \in \ii})$ for arbitrary sets of polynomials in the obvious way, where $\gcd(\varnothing) = 0$.
Similarly, we define the least common multiple $\lcm(\rho_1, \rho_2)$ to have leading coefficient $1$ whenever $\rho_1$ and $\rho_2$ are both nonzero, and we set $\lcm(\rho_1, \rho_2) = 0$ if $\rho_1 = 0$ or $\rho_2 = 0$. Again, we also define $\lcm(\set{\rho_i : i \in \ii})$ for arbitrary sets of polynomials in the obvious way, where $\lcm(\varnothing) = 1$.
These conventions ensure that (iii) and (iv) of Corollary \ref{corollary_polynomials_kernel_facts} below hold, and generalize to arbitrary sets of polynomials.
\begin{corollary}
    \label{corollary_polynomials_kernel_facts}
    Given polynomials $\rho, \rho_1, \rho_2 \in K[X]$, the following hold:
    \begin{enumerate}[(i)]
        \item Given $v \in \Ker(\rho_2)$, we also have $\rho_1[\theta](v) \in \Ker(\rho_2)$.
        \item If $\Ker(\rho^m) = \Ker(\rho^n)$ holds for some $n > m$, then $\Ker(\rho^m) = \Ker(\rho^n)$ holds for all $n \geq m$.
        \item $\Ker(\rho_1) + \Ker(\rho_2) = \Ker(\lcm(\rho_1, \rho_2))$.
        \item $\Ker(\rho_1) \cap \Ker(\rho_2) = \Ker(\gcd(\rho_1, \rho_2))$.
        \item If $\rho_1$ and $\rho_2$ are relatively prime, then $\rho_1[\theta]_{\restriction \Ker(\rho_2)} \colon \Ker(\rho_2) \to \Ker(\rho_2)$ is an automorphism with $\rho_1[\theta]_{\restriction \Ker(\rho_2)}^{-1} = \chi[\theta]_{\restriction \Ker(\rho_2)}$ for some $\chi \in K[X]$ that depends only on $\rho_1$ and $\rho_2$.
        In particular, $\Ker(\rho_2) \subseteq \Image(\rho_1)$.
    \end{enumerate}
    We have analogous statements to (i)-(iv) for images, in particular:
    \begin{enumerate}[(i)]
        \setcounter{enumi}{5}
        \item Given $v \in \Image(\rho_2)$, we also have $\rho_1[\theta](v) \in \Image(\rho_2)$.
        \item $\Image(\rho_1) \cap \Image(\rho_2) = \Image(\lcm(\rho_1, \rho_2))$.
    \end{enumerate}
\begin{proof}
    Points (i) and (vi) follow from the commutativity in Lemma \ref{lemma_add_mult_endo}.
    Point (ii) is clear.
    All remaining points follow essentially from Fact \ref{fact_euclid_domain}.
    We prove the inclusion ``$\subseteq$'' in (iv) as an example.
    Fix $v \in \Ker(\rho_1) \cap \Ker(\rho_2)$.
    By Fact \ref{fact_euclid_domain}, there are $\chi_1, \chi_2 \in K[X]$ such that $\gcd(\rho_1, \rho_2) = \chi_1 \cdot \rho_1 + \chi_2 \cdot \rho_2$.
    Hence
    $$
        \gcd(\rho_1, \rho_2)[\theta](v) = (\chi_1 \cdot \rho_1 + \chi_2 \cdot \rho_2)[\theta](v) = \chi_1[\theta](\rho_1[\theta](v)) + \chi_2[\theta](\rho_2[\theta](v)) = 0
    $$
    by Lemma \ref{lemma_add_mult_endo}, since $\rho_1[\theta](v) = 0 = \rho_2[\theta](v)$.
\end{proof}
\end{corollary}

\noindent If $\rho_1$ and $\rho_2$ are relatively prime, then (iii) and (iv) of Corollary \ref{corollary_polynomials_kernel_facts} yield the decomposition $\Ker(\rho_1 \cdot \rho_2) = \Ker(\rho_1) \oplus \Ker(\rho_2)$.
Also, whenever we write $\rho_1[\theta]_{\restriction \Ker(\rho_2)}^{-1}$, we mean the map defined as in (v) of Corollary \ref{corollary_polynomials_kernel_facts}.

\begin{lemma}[Euclidean Division]\label{lemma_remainder_rule}
    Given two polynomials $\rho, \xi \in K[X]$ with $\xi \neq 0$, we have
    $$\TKvsThe \models \forall x,y : \xi[\theta](x) = y \rightarrow \rho[\theta](x) = r[\theta](x) + \chi[\theta](y)$$
    for the unique polynomials $r, \chi \in K[X]$ with $\rho = r + \chi \cdot \xi$ and $\deg(r) < \deg(\xi)$.
\begin{proof}
    Let $(\VV, \theta)$ be a vector space with an endomorphism, and assume that $\xi[\theta](v) = u$ for some $v, u \in \VV$.
    Then Lemma \ref{lemma_add_mult_endo} yields
    $$
        \rho[\theta](v) = (r + \chi \cdot \xi)[\theta](v) = r[\theta](v) + \chi[\theta](\xi[\theta](v)) = r[\theta](v) + \chi[\theta](u).
    $$
    This shows that the desired implication holds in $\TKvsThe$.
\end{proof}
\end{lemma}
\noindent
The following is an easy consequence of the preceding lemma:

\begin{lemma}\label{lemma_bound_term_light}
    Let $\ux = (x_1, \dots, x_n)$ be a tuple of variables, and assume that a formula
    $$
        E(\ux; \uy) := \bigwedge\nolimits_{k=1}^n \xi_k[\theta](x_k) = \sum\nolimits_{l=1}^{k-1} Q_{k, l}[\theta](x_l) + \mu_k(\uy)
    $$
    is given, where the $\xi_k$ are non-zero polynomials, the $Q_{k, l}$ are polynomials, and the $\mu_k(\uy)$ are $\LKThe$-terms.
    Then any $\LKThe$-term $\tau(\ux; \uy)$ is equivalent modulo $\TKvsThe \cup \set{E(\ux; \uy)}$ to a term of the form
    $$
        \sum\nolimits_{k=1}^n \rho_k[\theta](x_k) + \mu'(\uy)
    $$
    with $\deg(\rho_k) < \deg(\xi_k)$ for all $k$, where $\mu'(\uy)$ is an $\LKThe$-term.
\begin{proof}
    We choose an $\LKThe$-term
    $
        \sum\nolimits_{k=1}^n \rho_k[\theta](x_k) + \mu'(\uy)
    $
    that is equivalent to $\tau(\ux; \uy)$ modulo $\TKvsThe \cup \set{E(\ux; \uy)}$, and for which $(\deg(\rho_n), \dots, \deg(\rho_1))$ is minimal with respect to the lexicographical ordering.
    If $\deg(\rho_k) \geq \deg(\xi_k)$ held for some $k$, then Lemma \ref{lemma_remainder_rule} together with the equation
    $$
        \xi_k[\theta](x_k) = \sum\nolimits_{l=1}^{k-1} Q_{k, l}[\theta](x_l) + \mu_k(\uy)
    $$
    from $E(\ux; \uy)$ would allow us to replace $\rho_k[\theta](x_k)$ by a term of smaller degree in $x_k$, contradicting the minimality of $(\deg(\rho_n), \dots, \deg(\rho_1))$.
\end{proof}
\end{lemma}

\begin{lemma}[Bézout's Identity]\label{lemma_gcd_rule}
    Given two polynomials $\xi_1, \xi_2 \in K[X]$, there are polynomials $\chi_1, \chi_2 \in K[X]$ such that
    $$
        \TKvsThe \models \forall x,y_1,y_2: (\xi_1[\theta](x) = y_1 \wedge \xi_2[\theta](x) = y_2) \rightarrow \gcd(\xi_1, \xi_2)[\theta](x) = \chi_1[\theta](y_1) + \chi_2[\theta](y_2).
    $$
\begin{proof}
    Choose $\chi_1, \chi_2 \in K[X]$ such that $\gcd(\xi_1, \xi_2) = \chi_1 \cdot \xi_1 + \chi_2 \cdot \xi_2$.
    This is possible by Fact \ref{fact_euclid_domain}.
    Now apply Lemma \ref{lemma_add_mult_endo}.
\end{proof}
\end{lemma}

\subsection{Kernel configurations and extensions of $T_\theta$}
\label{sec_baisc_stuff}

\noindent Recall that the goal is to provide a criterion for when $T_\theta$ has a model companion.
As noted in the introduction, we would like to add certain constraints to the map $\theta$.
We briefly illustrate which ones and why.
Take $K = \QQ$, let $\mm := (\RR, 0, +, (q \cdot)_{q \in \QQ}) \models \TQvs$, and consider the following endomorphisms of the underlying $\QQ$-vector space:
\begin{enumerate}[(i)]
    \item $\theta_1 \colon \RR \to \RR; x \mapsto \sqrt{2} \cdot x$.
    We have $\Ker((X^2 - 2)[\theta_1]) = \RR$ in the notation from Section \ref{sec_vec_with_endo}, or equivalently, $\Ker(\theta_1^2 - 2 \cdot \theta_1^0) = \RR$.
    \item $\theta_2 \colon \RR \to \RR; x \mapsto \pi \cdot x$.
    Given any polynomial $\rho \in \QQ[X] \setminus \set{0}$ and any non-zero $r \in \RR$, we have $\rho[\theta_2](r) = \rho(\pi) \cdot r \neq 0$, since $\pi$ is transcendental.
    Hence $\Ker(\rho[\theta_2]) = \set{0}$ for every $\rho \in \QQ[X] \setminus \set{0}$.
    In other words, $\rho[\theta_2]$ is injective for every $\rho \in \QQ[X] \setminus \set{0}$.
    \item Fix $n > 0$ and
    let $B = \set{a_{\alpha, i} : \alpha < |\RR|, i \in \ZZ} \dotcup \set{b_{\alpha, i} : \alpha < |\RR|, 0 \leq i < n}$ be a $\QQ$-basis for $\RR$.
    Given any $\alpha < |\RR|$, define $\theta_3(a_{\alpha, i}) := a_{\alpha, i+1}$ for all $i \in \ZZ$, $\theta_3(b_{\alpha, i}) := b_{\alpha, i+1}$ for all $i < n - 1$, and $\theta_3(b_{\alpha, n-1}) = 0$.
    We have thus defined $\theta_{3\restriction B}$, which uniquely extends to a $\QQ$-endomorphism $\theta_3 \colon \RR \to \RR$.
    This construction ensures that $\Ker(\theta_3^n)$ is not a nice set such as $\RR$ or $\set{0}$.
    However, the equality $\Ker(\theta_3^n) = \Ker(\theta_3^{n+1})$ holds, or equivalently, $\Ker(X^n[\theta_3]) = \Ker(X^{n+1}[\theta_3])$.
    This implies that the induced map $X[\theta_3] \colon \RR / \Ker(X^n[\theta_3]) \to \RR / \Ker(X^n[\theta_3])$ is injective.
\end{enumerate}
As we will see later, $T_\theta$ is inductive (Lemma \ref{lemma_inductive}).
Hence every model of $T_\theta$ is contained in some existentially closed model of $T_\theta$, and $T_\theta$ has a model companion if and only if its existentially closed models form an elementary class.
In particular, for each $i = 1, 2, 3$, there is some existentially closed model $(\mm'_i, \theta'_i) \supseteq (\mm, \theta_i)$ of $T_\theta$.
Ideally, the endomorphism $\theta'_i$ would still satisfy the respective equation above.
However, this is not the case.
Given any model $(\mm, \theta) \models T_\theta$, there is an extension $(\mm', \theta') \models T_\theta$ with elements $v_1, v_2, v_3$ such that \hbox{$(X^2 - 2)[\theta'](v_1) \neq 0$,} $v_2 \neq 0$, $X[\theta'](v_2) = 0$, and $X^{n+1}[\theta'](v_3) = 0 \neq X^n[\theta'](v_3)$.
Together with the existential closedness of each $(\mm'_i, \theta'_i)$, this implies that the respective equations are not satisfied by $\theta'_i$.

Thus, if we want to have $\Ker((X^2 - 2)[\theta'_1]) = \VV$, then we need $(\mm'_1, \theta'_1)$ to be an existentially closed model of $T_\theta \cup \set{\Ker((X^2 - 2)[\theta]) = \VV}$ instead of just $T_\theta$.
In general, we would like to find a criterion for the existence of a model companion for any theory of the form
$$
    T_\theta \cup \Set{\sum\nolimits_k \bigcap\nolimits_l \Ker(\rho_{j, k, l}) = \sum\nolimits_k \bigcap\nolimits_l \Ker(\eta_{j, k, l}) : j \in \jj},
$$
where $\jj$ is some possibly infinite index set, all sums and intersections are finite, and all the $\rho$ and $\eta$ are in $K[X]$.
Also recall the notation above Fact \ref{fact_euclid_domain}: we write $\Ker(\rho)$ and $\Image(\rho)$ instead of $\Ker(\rho[\theta])$ and $\Image(\rho[\theta])$, respectively, whenever $\theta$ is clear from the context.
For now, let $\Ff$ denote the set of all such theories.

Before looking for model companions, we would like to simplify $\Ff$.
Some theories in $\Ff$ may be inconsistent; for example,
$$
    T_\theta \cup \set{\Ker(X) = \Ker(0), \Ker(X) = \Ker(1)}.
$$
Other theories may be equivalent.
For instance, given any $T_1 \in \Ff$, let $T_2 \in \Ff$ be the theory obtained by multiplying all $\rho$ and $\eta$ that appear in $T_1$ by some fixed $\lambda \in K \setminus \set{0}$.
Then $T_1$ and $T_2$ are obviously equivalent.
There are also less obvious examples of equivalent theories.
\begin{align*}
    T_1 &:= T_\theta \cup \set{\Ker(X^5 + 2X^4 + X^3 + 3X^2 + 1) = \Ker(0), \Ker(X^3 + 4X^2 + X + 4) = \Ker(0)}, \\
    T_2 &:= T_\theta \cup \set{\Ker(X^2 + 1) = \Ker(0)}.
\end{align*}
These theories are also equivalent over $K = \QQ$, since
$$
    \gcd(X^5 + 2X^4 + X^3 + 3X^2 + 1, X^3 + 4X^2 + X + 4) = X^2 + 1
$$
by (iv) of Corollary \ref{corollary_polynomials_kernel_facts}.
Finally, equations of the form $\sum\bigcap\Ker(\dots) = \sum\bigcap\Ker(\dots)$ are too unwieldy to work with properly.

In the rest of this section, we define a family $\set{T^C_\theta : C \in \Cc} \subseteq \Ff$ of consistent theories extending $T_\theta$ such that every consistent theory in $\Ff$ is equivalent to some $T^C_\theta$.
Furthermore, we will have $T^{C_1}_\theta \not\equiv T^{C_2}_\theta$ whenever $C_1 \neq C_2$, and the formulas in each $T^C_\theta \setminus T_\theta$ will be much easier to work with.

\begin{lemma} \label{lemma_all_consis}
    Any theory of the form
    $$
        T_\theta \cup \Set{\sum\nolimits_k \bigcap\nolimits_l \Ker(\rho_{j, k, l}) = \sum\nolimits_k \bigcap\nolimits_l \Ker(\eta_{j, k, l}) : j \in \jj}
    $$
    where all sums and intersections are finite, all the $\rho$ and $\eta$ are polynomials over $K$, and $\jj$ is a possibly infinite index set, is either inconsistent or equivalent to a theory of one of the following forms:
    \begin{enumerate}[(i)]
        \item $T_\theta \cup \set{\Ker(\rho) = \VV}$, where $\rho$ is a monic polynomial with $\rho \neq 1$.
        \item $T_\theta \cup \set{\Ker(f^{n_f}) = \Ker(f^{n_f + 1}) : f \in Q}$, where $Q$ is a possibly infinite set of monic irreducible polynomials, and $n_f$ is a natural number for every $f \in Q$.
    \end{enumerate}
\begin{proof}
    Set $T^*_\theta := T_\theta \cup \set{{}\dots{} : j \in \jj}$, where the dots abbreviate the equations from the statement.
    We show step by step that $T^*_\theta$ can be assumed to be of a simpler form until it is either of form (i) or of form (ii).

    First, we eliminate the sums and intersections.
    By (iii) and (iv) of Corollary \ref{corollary_polynomials_kernel_facts}, we have
    $$
        \sum\nolimits_{k = 1}^{m} \bigcap\nolimits_{l = 1}^{n_k} \Ker(\rho_{k, l}) = \Ker(\lcm(\gcd(\rho_{k, l} : 1 \leq l \leq n_k) : 1 \leq k \leq m)).
    $$
    Therefore, we may assume that $T^*_\theta$ is of the form
    $$T_\theta \cup \set{\Ker(\rho_{j, 1}) = \Ker(\rho_{j, 2}) : j \in \jj}.$$
    Now consider an equation $\Ker(\rho_1) = \Ker(\rho_2)$ with $\rho_1, \rho_2 \neq 0$.
    For each $i = 1, 2$, we write $\rho_i = \lambda_i \cdot \prod_{k=1}^m f_k^{n_{i, k}}$ with $\lambda_i \in K \setminus \set{0}$, where the $f_k$ are monic irreducible polynomials and the $n_{i, k}$ are natural numbers.
    Here we allow $n_{i, k} = 0$ so that the same polynomials $f_k$ can be used for both $\rho_1$ and $\rho_2$.
    Again by (iii) and (iv) of Corollary \ref{corollary_polynomials_kernel_facts}, we have $\Ker(\rho_i) = \bigoplus_{k=1}^m \Ker(f_k^{n_{i, k}})$.
    Since for each $j$ we have either $\Ker(f_j^{n_{1, j}}) \subseteq \Ker(f_j^{n_{2, j}})$ or $\Ker(f_j^{n_{1, j}}) \supseteq \Ker(f_j^{n_{2, j}})$, part (ii) of Corollary \ref{corollary_polynomials_kernel_facts} yields
    $$
        \Ker(\rho_1) = \Ker(\rho_2) \quad\Leftrightarrow\quad \bigwedge\nolimits_{\!\!\!\!\underset{n_{1, k} \neq n_{2, k}}{k=1}}^m\Ker(f_k^{\min(n_{1, k}, n_{2, k})}) =  \Ker(f_k^{\min(n_{1, k}, n_{2, k}) + 1}).
    $$
    Hence, we may assume that $T^*_\theta$ is of the form
    $$T_\theta \cup \set{\Ker(f^{n_f}) = \Ker(f^{n_f + 1}) : f \in Q} \cup \set{\Ker(\rho_j) = \Ker(0) : j \in \jj_0}$$
    for some set $Q$ of monic irreducible polynomials, a family $(n_f : f \in Q)$ of natural numbers, and a subset $\set{\rho_j : j \in \jj_0} \subseteq K[X] \setminus \set{0}$.
    If $\jj_0 = \varnothing$, then we are in case (ii), so we are done.

    Now assume that $\jj_0 \neq \varnothing$.
    Since $\Ker(0) = \VV$, part (iv) of Corollary \ref{corollary_polynomials_kernel_facts} shows that we may assume that $T^*_\theta$ is of the form
    $$
        T_\theta \cup \set{\Ker(f^{n_f}) = \Ker(f^{n_f + 1}) : f \in Q} \cup \set{\Ker(\rho) = \VV}
    $$
    for some monic $\rho \in K[X] \setminus \set{0}$.
    Write $\rho = \prod_{k=1}^m f_k^{r_k}$ as a product of its irreducible factors.
    By (iii) and (iv) of Corollary \ref{corollary_polynomials_kernel_facts}, the equation $\Ker(\rho) = \VV$ is equivalent to $\bigoplus_{k=1}^m \Ker(f_k^{r_k}) = \VV$.
    This already implies $\Ker(f^n) = \Ker(f^{n + 1})$ for every monic irreducible $f \in K[X]$ and every $n \geq \max\set{r : f^r \mid \rho}$.
    Hence we may assume that $Q \subseteq \set{f_1, \dots, f_m}$ and that $f^{n_f + 1} \mid \rho$ for each $f \in Q$.
    Without loss of generality, let $f_1 \in Q$.
    Since $n_{f_1} < r_1$, Corollary \ref{corollary_polynomials_kernel_facts} shows that
    \begin{align*}
        \Ker(\rho) = \VV \wedge \Ker\big(f_1^{n_{f_1}}\big) = \Ker\big(f_1^{n_{f_1}+1}\big)
        &\quad\Leftrightarrow\quad \Ker\!\left(f_1^{n_{f_1}} \cdot \prod\nolimits_{k=2}^m f_k^{r_k}\right) = \VV.
    \end{align*}
    Repeating this argument, we may assume that $Q = \varnothing$.
    Thus we may assume that $T^*_\theta$ is of the form $T_\theta \cup \set{\Ker(\rho) = \VV}$, i.e., of form (i).
    Finally, $\rho \neq 1$, since otherwise we would have $\set{0} = \Ker(1) = \VV$, while $T$ states that $\VV$ is infinite and we assumed that $T^*_\theta$ is consistent.
\end{proof}
\end{lemma}

\noindent Thus, any consistent theory of the form $T_\theta \cup \set{\sum\bigcap \Ker(\dots) = \sum\bigcap \Ker(\dots) : j \in \jj}$ is either
\begin{enumerate}[(i)]
    \item \textbf{algebraic}, in the sense that there is a monic polynomial $\rho \in K[X] \setminus \set{1}$ such that $\rho[\theta] = 0$ for every model $(\mm, \theta)$ of it, or
    \item \textbf{transcendental}, in the sense that it is not algebraic.
\end{enumerate}
Recall that our goal is to find a family $\set{T^C_\theta : C \in \Cc}$ of theories such that every consistent theory of the form $T_\theta \cup \set{\sum\bigcap \Ker(\dots) = \sum\bigcap \Ker(\dots) : j \in \jj}$ is equivalent to some $T^C_\theta$, and such that $T^{C_1}_\theta \not\equiv T^{C_2}_\theta$ whenever $C_1 \neq C_2$.
The next step is to define the index set $\Cc$.
By the above, it seems natural to define two sets $\Ccalg$ and $\Cctrans$ that describe the algebraic and transcendental cases first, and then set $\Cc := \Ccalg \dotcup \Cctrans$.
However, we introduce the full set $\Cc$ before the defining $\Ccalg$ and $\Cctrans$, so that later proofs are as uniform as possible.

\begin{definition}
    Let $\Kp{} := \set{f \in K[X] : \text{$f$ is monic and irreducible}}$.
\end{definition}

\noindent We now define the index set $\Cc$:

\begin{definition} \label{def_kernel_conf}
    We call a pair $(c, d)$ a \textbf{kernel configuration} if
    \begin{enumerate}[(i)]
        \item $c \colon \Kp{} \to \NN \cup \set{\infty}$ is a function; and
        \item $d \in \NN_{> 0} \cup \set{\infty}$ is either $\infty$ or satisfies $d = \sum_{f \in \Kp{}} \deg(f) \cdot c(f)$.
    \end{enumerate}
    We let $\Cc$ denote the set of all kernel configurations.
    Given such a kernel configuration $C = (c, d) \in \Cc$, we set $C(f) := c(f)$ for all $f \in \Kp{}$.
    We define the \textbf{degree} of $C$ by $\deg(C) := d$.
    We say that $C$ is \textbf{algebraic} if $\deg(C) < \infty$.
    In this case, we define the \textbf{minimal polynomial} of $C$ by $\mipo(C) := \prod_{f \in \Kp{}} f^{C(f)}$.
    Since $\deg(C) < \infty$, only finitely many factors are different from $1$, so this product is well defined.
    We say that $C$ is \textbf{transcendental} if $\deg(C) = \infty$.
    We let $\Ccalg$ and $\Cctrans$ denote the sets of all algebraic and all transcendental kernel configurations, respectively.
\end{definition}

\noindent Note that the set of kernel configurations depends on the field $K$.
We are now ready to define the family $\set{T^C_\theta : C \in \Cc}$:

\begin{definition} \label{def_T_C_theta}
    Given $C \in \Cc$ and an endomorphism $\theta \colon \VV \to \VV$, we say that $\theta$ is a \textbf{$\mathbf{C}$-endomorphism} if one of the following holds:
    \begin{enumerate}[(i)]
        \item $C$ is algebraic and $\Ker(\mipo(C)) = \VV$, that is, $\mipo(C)[\theta] = 0$.
        \item $C$ is transcendental and $\Ker(f^{C(f)}) = \Ker(f^{C(f)+1})$ for all $f \in \Kp{}$ with $C(f) < \infty$.
    \end{enumerate}
    We define $T^C_\theta := T_\theta \cup \set{\text{``$\theta$ is a $C$-endomorphism''}}$.
\end{definition}

\begin{corollary}
    Every consistent theory of the form
    $$
        T_\theta \cup \Set{\sum\bigcap \Ker(\dots) = \sum\bigcap\Ker(\dots) : j \in \jj}
    $$
    is equivalent to $T^C_\theta$ for some kernel configuration $C \in \Cc$.
\begin{proof}
    By Lemma \ref{lemma_all_consis}, we may assume that the given theory is of one of the following forms:
    \begin{enumerate}[(i)]
        \item $T_\theta \cup \set{\Ker(\rho) = \VV}$, where $\rho$ is a monic polynomial with $\rho \neq 1$.
        In this case, we define $C := (c, d)$ by setting $d := \deg(\rho)$ and
        $$
            c \colon \Kp{} \to \NN \cup \set{\infty}; \quad f \mapsto \max\set{n \in \NN : f^n \mid \rho}.
        $$
        \item $T_\theta \cup \set{\Ker(f^{n_f}) = \Ker(f^{n_f + 1}) : f \in Q}$, where $Q$ is a possibly infinite set of monic irreducible polynomials and $n_f$ is a natural number for every $f \in Q$.
        In this case, we define $C := (c, d)$ by setting
        $$
            c \colon \Kp{} \to \NN \cup \set{\infty}; \quad f \mapsto \begin{cases}
            n_f & \text{if $f \in Q$} \\
            \infty & \text{if $f \not\in Q$}
            \end{cases} \quad \text{and} \quad d := \infty.
        $$
    \end{enumerate}
    In both cases, it is straightforward to verify that $C$ is a kernel configuration (notice that $\mipo(C) = \rho$ in case (i)), and that $T^C_\theta$ is equivalent to the respective theory.
\end{proof}
\end{corollary}
\noindent It remains to show that each theory $T^C_\theta$ is consistent and that we have $T^{C_1}_\theta \not\equiv T^{C_2}_\theta$ whenever $C_1 \neq C_2$.
We will prove both facts at the end of this section.
For now, we record a few basic observations about kernel configurations, $C$-endomorphisms, and the theories $T^C_\theta$.
We begin with the transcendental case:

\begin{remark}
    The condition $\Ker(f^{C(f)+1}) = \Ker(f^{C(f)})$ in part (ii) of Definition \ref{def_T_C_theta} is equivalent to the map $f[\theta] \colon \VV / \Ker(f^{C(f)}) \to \VV / \Ker(f^{C(f)})$ being injective.
    If $C(f) = 0$, then $\Ker(f^{C(f)}) = \Ker(1) = \set{0}$.
    In this case, the condition is therefore equivalent to $f[\theta] \colon \VV \to \VV$ being injective.
    This yields two extreme cases:
    \begin{enumerate}[(i)]
        \item Let $C_0 \in \Cctrans$ be given by $C_0(f) = 0$ for all $f \in \Kp{}$.
        Then $T^{C_0}_\theta$ implies that $\rho[\theta]$ is injective for every nonzero polynomial $\rho \in K[X]$.
        \item Let $C_\infty \in \Cctrans$ be given by $C_\infty(f) = \infty$ for all $f \in \Kp{}$.
        Then $T^{C_\infty}_\theta$ is equivalent to $T_\theta$ and imposes no injectivity requirements on $\rho[\theta]$ beyond the trivial case $\rho \in K \setminus \set{0}$.
    \end{enumerate}
\end{remark}

\begin{remark} \label{remark_alg_kc}
    It is easy to check that the following hold for any algebraic kernel configuration $C \in \Ccalg$:
    \begin{enumerate}[(i)]
        \item Given a monic irreducible polynomial $f$, we have $C(f) = \max\set{n \in \NN : f^n \mid \mipo(C)}$.
        \item Any $C$-endomorphism $\theta$ also satisfies $\Ker(f^{C(f)}) = \Ker(f^{C(f)+1})$ for all $f \in \Kp{}$.
        This is almost condition (ii) of Definition \ref{def_T_C_theta}, except that $C$ is algebraic rather than transcendental.
        \item We have $\deg(\mipo(C)) = \deg(C)$.
        \item Given another $C' \in \Ccalg$, we have $C = C'$ if and only if the minimal polynomials coincide, that is, if and only if $\mipo(C) = \mipo(C')$.
    \end{enumerate}
    Point (ii) allows us to treat algebraic and transcendental kernel configurations uniformly in many cases.
    Together, points (i), (iii) and (iv) show that one can also regard $\Ccalg$ as the set $\set{\!\text{``monic polynomials''}} \setminus \set{1}$.
\begin{proof}
    Point (ii) follows from
    $$
    \VV = \Ker(\mipo(C)) = \Ker\Big(\prod\nolimits_{f \in \Kp{}} f^{C(f)}\Big) = \bigoplus\nolimits_{f \in \Kp{}} \Ker(f^{C(f)}).
    $$
    The remaining points are clear from the definitions.
\end{proof}
\end{remark}

\noindent When we define a kernel configuration $C$, we will almost always specify whether $C$ lies in $\Ccalg$ or $\Cctrans$ and then define either its minimal polynomial $\mipo(C)$ (in the algebraic case) or the function $f \mapsto C(f)$ (in the transcendental case).
In some situations, especially when we define one kernel configuration from another, we may also specify the function $f \mapsto C(f)$ in the algebraic case.

\begin{remark}
    Looking at Definition \ref{def_T_C_theta}, we see that an endomorphism $\theta \colon \VV \to \VV$ can be a $C$-endomorphism for more than one kernel configuration $C$.
    We therefore define a partial order $\preceq$ on $\Cc$ by setting $C_1 \preceq C_2$ if and only if every $C_1$-endomorphism is also a $C_2$-endomorphism.
    One readily checks that $C_1 \preceq C_2$ holds if and only if $\deg(C_1) \leq \deg(C_2)$ and $C_1(f) \leq C_2(f)$ for all $f \in \Kp{}$.
    It follows that the unique transcendental kernel configuration $C_\infty$ with $C_\infty(f) = \infty$ for all $f \in \Kp{}$ is maximal in $(\Cc, \preceq)$ and that the unique transcendental kernel configuration $C_0$ with $C_0(f) = 0$ for all $f \in \Kp{}$ is the minimal transcendental kernel configuration.
    The minimal algebraic kernel configurations are exactly those $C$ for which $\mipo(C)$ is irreducible.
    Finally, for any endomorphism $\theta$, there is a minimal $C_\theta \in \Cc$ such that $\theta$ is a $C_\theta$-endomorphism.
\end{remark}

\begin{notation}
    We introduce a few more notations for working with a kernel configuration $C \in \Cc$:
    \begin{enumerate}[(i)]
        \item Given $f \in \Kp{}$ with $C(f) < \infty$, we write $f^C$ instead of $f^{C(f)}$, $f^{C+k}$ instead of $f^{C(f) + k}$, and so on.
        \item Given a finite set $F \subseteq \Kp{}$ with $C(f) < \infty$ for all $f \in F$, we set \hbox{$F^C := \prod_{f \in F} f^C$}.
        \item We define the following subsets of $\Kp{}$:
        \begin{itemize}
            \item $\Kp{C<\infty} := \set{f \in \Kp{} : C(f) < \infty}$
            \item $\Kp{0<C<\infty} := \set{f \in \Kp{} : 0 < C(f) < \infty}$
            \item $\Kp{C=0} := \set{f \in \Kp{} : C(f) = 0}$
            \item $\Kp{C=\infty} := \set{f \in \Kp{} : C(f) = \infty}$.
        \end{itemize}
    \end{enumerate}
\end{notation}
\begin{observation} \label{obs_mipo_prod}
    With the notation above, for any algebraic kernel configuration $C \in \Ccalg$, we have
    $$
    \mipo(C) = \prod\nolimits_{f \in \Kp{0<C<\infty}} f^C = (\Kp{0<C<\infty})^C.
    $$
    Hence, by (iii) and (iv) of Corollary \ref{corollary_polynomials_kernel_facts}, for any $C$-endomorphism $\theta \colon \VV \to \VV$, we have
    $$
    \VV = \Ker(\mipo(C)) = \bigoplus\nolimits_{f \in \Kp{0<C<\infty}} \Ker(f^C).
    $$
\end{observation}
\noindent We often want to extend a $C$-endomorphism of $\VV$ to a $C$-endomorphism of a larger vector space $\VV'$ in such a way that a certain property is satisfied.
We do this in two steps:
\begin{enumerate}[(i)]
    \item Choose some $B \subseteq \VV'$ linearly independent over $\VV$ (if $B$ is not already clear from the context) and extend $\theta$ to a $C$-endomorphism of $\spanA{B\VV}{K}$ in such a way that the property is satisfied.
    \item Extend the endomorphism from the first step to all of $\VV'$.
\end{enumerate}
This obviously does not work for all properties, but it does for ``existential'' ones such as ``there is $v$ in $\Ker(f^{n+1}) \setminus \Ker(f^n)$''.
For the second step, we will always use the following construction.

\begin{lemma}[Standard Construction]
    \label{lemma_standart_construction}
    Given a $C$-endomorphism $\theta \colon \VV \to \VV$ and a vector space $\VV' \supset \VV$ with $\dim(\VV'/\VV) \geq \aleph_0$, there exists a $C$-endomorphism $\theta' \colon \VV' \to \VV'$ extending $\theta$.
\begin{proof}
    Let $B$ be a basis of $\VV'$ over $\VV$.
    If $C$ is algebraic, define $\rho := \sum_{i=0}^d (\rho)_i \cdot X^i := \mipo(C)$.
    Write
    $$
    B = \set{u_{\alpha}^i : \alpha < \dim(\VV'/\VV), 0 \leq i < \deg(\rho)}
    $$
    with $u_{\alpha}^i \neq u_{\beta}^j$ for $(i, \alpha) \neq (j, \beta)$, and define $\theta'$ outside of $\VV$ by setting $\theta'(u_\alpha^i) := u_\alpha^{i+1}$ for $i < d-1$ and
    $
    \theta'(u_\alpha^{d-1}) := -\sum_{i=0}^{d-1} (\rho)_i \cdot u_\alpha^i.
    $
    If $C$ is transcendental, write
    $$
    B = \set{u_{\alpha}^i : \alpha < \dim(\VV'/\VV), i \in \omega}
    $$
    and define $\theta'(u_\alpha^i) := u_\alpha^{i+1}$.
    In both cases, it is easy to verify that $\theta'$ is a $C$-endomorphism.
\end{proof}
\end{lemma}

\noindent With our \standartConstruction{} we can now prove the remaining properties of the family $\set{T^C_\theta : C \in \Cc}$:

\begin{lemma} \label{lemma_C_properties}
    The following holds:
    \begin{enumerate}[(i)]
        \item The theory $T^C_\theta$ is consistent for every $C \in \Cc$.
        \item Given $C_1 \neq C_2 \in \Cc$, we have $T^{C_1}_\theta \not\equiv T^{C_2}_\theta$.
    \end{enumerate}
\begin{proof}
    For (i), note that the trivial endomorphism $\set{0} \to \set{0}$ is a $C$-endomorphism for every $C \in \Cc$.
    Given a model $\mm \models T$ with $\dim(\VV) \geq \aleph_0$, which clearly exists, we can extend $\set{0} \to \set{0}$ to a $C$-endomorphism $\theta \colon \VV \to \VV$ using our \standartConstruction{}.

    For (ii), it is enough to find either a $C_1$-endomorphism that is not a $C_2$-endomorphism or vice versa.
    To see this, take some $\mm \models T$ with sufficiently large $\VV$, identify the domain of the chosen endomorphism with a subspace of $\VV$, and use our \standartConstruction{} to extend it to all of $\VV$.
    If $C_1$ is transcendental and $C_2$ is algebraic, choose an $\aleph_0$-dimensional vector space $\VV = \spanA{u^i : i \in \NN}{K}$, where the $u^i$ are linearly independent, and define $\theta(u^i) := u^{i+1}$.
    Then $\theta$ is not a $C_2$-endomorphism, since $\Ker(\mipo(C_2)[\theta]) = \set{0} \neq \VV$.
    Moreover, we have \hbox{$\Ker(f^{n+1}[\theta]) = \set{0} = \Ker(f^n[\theta])$} for all $f \in \Kp{}$ and $n \geq 0$, so $\theta$ is a $C_1$-endomorphism.

    Now let $C_1$ and $C_2$ both be algebraic or both be transcendental.
    Without loss of generality, there is some $f \in \Kp{}$ such that $C_1(f) < C_2(f)$.
    Define $\rho := \sum_{i=0}^d (\rho)_i \cdot X^i := f^{C_1(f)+1}$.
    Choose a $d$-dimensional vector space $\VV := \spanA{u^i : 0 \leq i < d}{K}$ and define an endomorphism $\theta$ of $\VV$ by
    $$
    \theta(u^i) := u^{i+1} \text{ for all $i \in \set{0, \dots, d-2}$}
    $$
    and
    $$
    \theta(u^{d-1}) := -\sum\nolimits_{i=0}^{d-1} (\rho)_i \cdot u^i.
    $$
    Notice that $u^0 \notin \Ker(f^{C_1}[\theta])$, but for all $n > C_1(f)$, we have $\Ker(f^n[\theta]) = \VV$ and hence $\Ker(f^n[\theta]) = \Ker(f^{n+1}[\theta])$.
    This shows that $\Ker(f^{C_1}[\theta]) \neq \Ker(f^{C_1+1}[\theta])$, so $\theta$ is not a $C_1$-endomorphism.
    By (ii) of Remark \ref{remark_alg_kc}, this condition is also required for algebraic kernel configurations.
    For any irreducible polynomial $g \in \Kp{} \setminus \set{f}$ and $n \geq 0$, we have $\Ker(g^n[\theta]) = \set{0} = \Ker(g^{n+1}[\theta])$, so in the transcendental case $\theta$ is a $C_2$-endomorphism.
    In the algebraic case, note that we have $f^{C_1+1} \mid \mipo(C_2)$, so
    $$
    \Ker(\mipo(C_2)[\theta]) \supseteq \Ker(f^{C_1+1}[\theta]) = \VV.
    $$
    Thus, $\theta$ is a $C_2$-endomorphism in this case as well.
\end{proof}
\end{lemma}

\noindent Note that (ii) of Lemma \ref{lemma_C_properties} would also follow from the characterization of existentially closed models of $T^C_\theta$, which we state in Section \ref{sec_char}.

If one agrees with the statement that $L_\theta$-sentences of the form \hbox{$\Ker(\mipo(C)) = \VV$} and \hbox{$\Ker(f^{C+1}) = \Ker(f^C)$} are easier to work with than sentences of the form 
$$
\sum\bigcap \Ker(\dots) = \sum\bigcap \Ker(\dots),
$$
then we have shown all the properties of the family $\set{T^C_\theta : C \in \Cc}$ mentioned at the beginning of this section; see the paragraph above Lemma \ref{lemma_all_consis}.

\begin{remark}[Multiplicative endomorphism of ACF]\label{rk:ACFcase}
    Given a divisible abelian group $(\mathbb{G}, 1, \cdot)$, the quotient $\mathbb{G}/\Tor(\mathbb{G})$ is a $\QQ$-vector space (where $\Tor(\mathbb{G})$ is the torsion of $\mathbb{G}$).
    In this case, we say that an endomorphism $\theta$ of $(\mathbb{G}, 1, \cdot)$ is a \textbf{$\mathbf{C}$-endomorphism} if $\theta_{\restriction \mathbb{G}/\Tor(\mathbb{G})}$ is a $C$-endomorphism of the $\QQ$-vector space $\mathbb{G}/\Tor(\mathbb{G})$ as in Definition \ref{def_T_C_theta}.
    Since we have torsion, we can no longer uniquely divide by any $n \in \NN$, so we must work with polynomials in $\ZZ$.
    If $\mathbb{G}$ is definable in a theory $T$ and the torsion subgroup $\Tor(\mathbb{G})$ is not definable, it is not trivial to find consistent theories that imply
    \begin{align}
        T_\theta \cup \set{\text{``$\theta$ is a $C$-endomorphism of $(\mathbb{G}, 1, \cdot)$''}} \cup \set{\forall x \not\in \mathbb{G} : \theta(x) = 1}. \label{tag_lol_theory_torsion}
    \end{align}
    We quickly highlight the most promising approaches for $T = \ACF$. Given any monic polynomial $\rho \in \QQ[X]$, we let $\rho_\ZZ$ be the multiple of $\rho$ in $\ZZ[X]$ with the smallest possible positive leading coefficient.
    \begin{enumerate}[(i)]
        \item If $C$ is algebraic, then one could consider the following theory:
        $$
        \ACF_{\theta} \cup \set{\forall x \in \mathbb{G} : \mipo(C)_\ZZ[\theta](x) = 1}.
        $$
        Using compactness, one can see that this is the only possible approach. More precisely, any consistent theory that implies the ``theory'' (\ref{tag_lol_theory_torsion}) with $T = \ACF{}$ must already imply the theory above.
        However, since the additional sentence also needs to hold on the torsion group (which is, up to isomorphism, $\bigoplus_{p\in \PP \setminus\set{\Char}}\! \QQ_p/\ZZ_p$), one can show that this theory is consistent with ``$\Char = q$'' if and only if $\mipo(C)$ has a root in $\ZZ_p$ for every $p \in \PP\setminus \set{q}$.
        \item If $C$ is transcendental, fix a function $n \colon (\QQp{C<\infty} \times \NN_{>0}) \to \NN$, and consider the theory $\ACF_\theta$ extended by the sentence
        $$
        \forall x \in \mathbb{G} :f_\ZZ^{C+1}[\theta](x) \in \Tor_{\leq m}(\mathbb{G}) \rightarrow f_\ZZ^{C}[\theta](x) \in \Tor_{\leq n(f, m)}(\mathbb{G})
        $$
        for every $m \in \NN_{>0}$ and $f \in \QQp{C<\infty}$. Here $\Tor_{\leq m}(\mathbb{G})$ denotes the set of all torsion elements of order at most $m$.
        For any transcendental kernel configuration, one can find such a function $n$ for which the theory is consistent.
        Using compactness, one can again show that any consistent theory that implies the ``theory'' (\ref{tag_lol_theory_torsion}) with $T = \ACF{}$ must already imply a theory as described above for some $n\colon (\QQp{C<\infty} \times \NN_{>0}) \to \NN$.
    \end{enumerate}
    We believe that our main results generalize to these settings, especially to the one presented in the transcendental case.
    One should notice that the theory $\operatorname{ACFH}$, the model companion of algebraically closed fields with a multiplicative endomorphism, mentioned in the introduction, is the model companion of the theory in (ii) for $C = C_\infty$.

    An important consequence of Hrushovski's ``twisted'' Lang-Weil estimate \cite{Hru22} is that any ultraproduct $\prod_{\mathcal{U}} (\overline{\mathbb{F}}_p, x \mapsto x^p)$ of Frobenius difference fields is a model of $\operatorname{ACFA}$ for any non-principal ultrafilter $\mathcal U$ on the set $\PP$ of prime numbers.
    However, it is quite challenging to describe the theory of the ultraproduct \hbox{$\prod_{\mathcal{U}} (\overline{\mathbb{F}}_p, x \mapsto x^{s_p})$} for an arbitrary sequence $(s_p)_{p \in \PP}$ of natural numbers and a non-principal ultrafilter $\mathcal{U}$ on $\PP$.
    As mentioned in \cite{dEl25}, the theory $\operatorname{ACFH}$ is \textit{not} a good candidate for the theory of $\prod_{\mathcal{U}} (\overline{\mathbb{F}}_p, x \mapsto x^{s_p})$.
    One reason for this is the lack of control over the kernels of the definable endomorphisms.
    As it turns out, any interesting, meaning non-principal, ultraproduct $\prod_{\mathcal{U}} (\overline{\mathbb{F}}_p, x \mapsto x^{s_p})$ will not satisfy any of our kernel configurations besides $C_\infty$.
    Hence, the additional constraints on the kernels in our setting seem too strong to tackle the problem of describing the theories of these ultraproducts.
\end{remark}

\begin{remark}[Kernel Configurations and the Ziegler spectrum] \label{remark_ziegler_spectrum}
    The reader might ask why we are only working with equalities of certain sums and intersections of kernels but do not allow images to appear in them.
    For example (setting \hbox{$T := \TKvs$} for simplicity), we could consider theories of the form
    $$
    T^*_{K\operatorname{-vs},\theta} := \TKvsThe \cup \Set{\sum\nolimits_{k}\bigcap\nolimits_{l}F_{j, k, l}(\rho_{j, k, l}[\theta]) = \sum\nolimits_{k}\bigcap\nolimits_{l} G_{j, k, l}(\eta_{j, k, l}[\theta]) : j \in \jj},
    $$
    where all sums and intersections are finite, all the $\rho$'s and $\eta$'s are polynomials over $K$, and all $F$'s and $G$'s lie in $\set{\Ker, \Image}$. As any $K$-vector space with an endomorphism is a $K[X]$-module and as the sets $\sum\nolimits_{k}\bigcap\nolimits_{l}F_{j, k, l}(\rho_{j, k, l}[\theta])$ and $\sum\nolimits_{k}\bigcap\nolimits_{l} G_{j, k, l}(\eta_{j, k, l}[\theta])$ are defined using pp-formulas, we see that the models of $T^*_{K\operatorname{-vs},\theta}$ form a definable category of $K[X]$-modules (see Section 3.4.1 in \cite{Pre09}).
    
    The Ziegler spectrum $\operatorname{Zg}_{K[X]}$ is the set $\operatorname{pinj}_{K[X]}$ of isomorphism classes of indecomposable pure-injective $K[X]$-modules (excluding the zero module) whose closed sets are of the form
    $$
    [H] = \set{N \in \operatorname{pinj}_{K[X]} : \phi(N) = \psi(N) \text{ for all } (\phi(x), \psi(x)) \in H},
    $$
    where $H$ is an arbitrary set of pp-pairs $(\phi(x), \psi(x))$ in a single free variable, i.e., both $\phi(x)$ and $\psi(x)$ are pp-formulas in a single free variable (formulas of the form $\exists y : E(x; y)$, where $E$ is a conjunction of equations in the language of $K[X]$-modules), and $\psi(x) \models \phi(x)$ holds in the theory of $K[X]$-modules.
    The important fact that we need is that there is a one-to-one correspondence between definable categories of $K[X]$-modules and closed subsets of the Ziegler spectrum $\operatorname{Zg}_{K[X]}$, and that this correspondence is given via $\mathcal{X} \mapsto \mathcal{X} \cap \operatorname{pinj}_{K[X]}$ (see, e.g., Theorem 5.1.1 of \cite{Pre09}).
    Further details can be found in Section 5.1.1 of \cite{Pre09}.

    Now, for any closed and non-empty subset $[H] \subseteq \operatorname{Zg}_{K[X]}$ of the Ziegler spectrum as above, we can define the following theory
    $$
    T^{[H]}_{K\operatorname{-vs},\theta} := \TKvsThe \cup \Set{\forall x : \phi(x) \leftrightarrow \psi(x) : (\phi(x), \psi(x)) \in H}.
    $$
    As the models of $T^*_{K\operatorname{-vs},\theta}$ form a definable category of $K[X]$-modules, \hbox{$T^*_{K\operatorname{-vs},\theta} \equiv T^{[H]}_{K\operatorname{-vs},\theta}$} holds for some closed subset $[H]$ of the Ziegler spectrum by the above.
    Once one knows enough about the existentially closed models of our theories $\TKvsTheC$, one can show that for any non-empty closed set $[H]$ of the Ziegler spectrum there is some kernel configuration $C \in \Cc$ such that
    $$
    \TKvsTheC \quad \subseteq \quad T^{[H]}_{K\operatorname{-vs},\theta} \quad \subseteq \quad \text{model companion of }\TKvsTheC
    $$
    holds.
    Hence the model companions of $T^{[H]}_{K\operatorname{-vs},\theta}$ and $\TKvsTheC$ must be the same. This also generalizes to the case where $T \neq \TKvs$.
    This means that the introduction of images, as mentioned at the beginning of this remark, would only lead to more complexity while not providing anything new.

    We also obtain another description of our kernel configurations. 
    Let $\Hh^{\qf}$ be the set of all pp-pairs $(\phi(x), \psi(x))$ with both $\phi(x)$ and $\psi(x)$ being quantifier-free pp-formulas in a single free variable.
    A quantifier-free pp-formula in a single free variable is just a conjunction $\bigwedge_{k=1}^n \rho_{k}[\theta](x) = \eta_{k}[\theta](x)$ (written in the language of a $K$-vector space with an endomorphism), and hence equivalent to the formula $x \in \Ker(\rho)$ for $\rho := \gcd(\rho_k - \eta_k : k \in \set{1, \dots, n})$. 
    Therefore, for any pp-pair $(\phi(x), \psi(x)) \in \Hh^{\qf}$ the condition $\forall x : \phi(x) \leftrightarrow \psi(x)$ is a condition of the form $\Ker(\rho) = \Ker(\eta)$ for two polynomials $\rho, \eta$. 
    With this and Lemma \ref{lemma_all_consis}, one can now show that there is a one-to-one correspondence between the non-empty closed subsets of the Ziegler spectrum of the form $[H]$, with $H \subseteq \Hh^{\qf}$, and our set of kernel configurations $\Cc$.
\end{remark}

\section{Characterization of Existentially Closed Models}

\noindent Now that we have our family $\set{T^C_\theta : C \in \Cc}$ in place, we can start looking for model companions of $T^C_\theta$.
We begin with one easy special case.

\begin{remark} \label{def_trivial}
    If $C$ is \textbf{trivial}, that is, if $\deg(C) = 1$, then the models of $T_\theta^C$ and $T$ are interdefinable.
    Hence, the model companion of $T_\theta^C$ is $T_\theta^C$ itself.
\begin{proof}
    Unwrapping all definitions, we obtain
    $$T_\theta^C = T_\theta \cup \set{\Ker((X + \lambda)[\theta]) = \VV}$$
    for some $\lambda \in K$.
    This implies that $T_\theta^C \models \theta_{\restriction \VV} = -\lambda \Id_{\VV}$.
    We conclude since $-\lambda \Id_{\VV}$ is definable in $T$ and $T$ is model-complete.
\end{proof}
\end{remark}

\noindent So, for the rest of this thesis, we will almost always assume that $C$ is non-trivial.

\begin{lemma} \label{lemma_inductive}
    For any $C \in \Cc$, the theory $T^C_\theta$ is inductive, i.e., given a $\subseteq$-chain of models $((\mm_\alpha, \theta_\alpha) : \alpha < \kappa)$ of $T^C_\theta$, the union $(\bigcup_{\alpha < \kappa} \mm_\alpha, \bigcup_{\alpha < \kappa} \theta_\alpha)$ is also a model of $T^C_\theta$.
\begin{proof}
    This follows directly from the model completeness of $T$, the fact that $\VV$ is always defined by the same formulas, and the fact that we are dealing with vector space endomorphisms.
\end{proof}
\end{lemma}

\noindent The inductiveness of $T^C_\theta$ has two consequences.
First, every model of $T_\theta^C$ is contained in an existentially closed one.
Second, a model companion of $T_\theta^C$ exists if and only if the class of existentially closed models of $T^C_\theta$ is $L_\theta$-axiomatizable.
If that is the case, the model companion is precisely this axiomatization.
Hence, we need to find a description of existentially closed models that can be first-order axiomatized.
A first step toward such a potentially first-order axiomatizable description is to simplify all existential $L_\theta$-formulas:

\begin{lemma}
\label{lemma_fml}
    For every quantifier free $L_\theta$-formula $\varphi(\uz; \uw)$, there is an $L$-formula $\psi(\ux^0, \ux^1; \uw)$ such that, modulo $T_\theta$, we have
    $$
    \exists \uz : \varphi(\uz; \uw) \quad \equiv \quad \exists \ux \in \VV : \psi(\ux, \theta(\ux); \uw).
    $$
\begin{proof}
    Suppose that we have a formula of the form
    \begin{align}
       \exists \uz: \bigwedge\nolimits^m_{i=1} x_i = s_i(\uz; \ux^0_{<i} \ux^1_{<i}; \uw) \wedge \bigwedge\nolimits_{j=1}^n t_j(\uz; \ux^0 \ux^1; \uw) \not\in \VV \wedge \psi(\uz; \ux^0 \ux^1; \uw), \label{tag_lemma_fml}
    \end{align}
    where the $s_i$ and $t_j$ are $L$-terms and $\psi(\uz; \ux^0 \ux^1; \uw)$ is a quantifier-free $L_\theta(0)$-formula.
    If this formula, and therefore $\psi(\uz; \ux^0 \ux^1; \uw)$, is not an $L$-formula, then there is an $L$-term $s_{m+1}(\uz; \ux^0 \ux^1; \uw)$ such that $\theta(s_{m+1}(\uz; \ux^0 \ux^1; \uw))$ appears in $\psi(\uz; \ux^0 \ux^1; \uw)$.
    Let the formulas $\psi_0(\uz; \ux^0 \ux^1; \uw)$ and $\psi_1(\uz; \ux^0 x^0_{m+1} \ux^1 x^1_{m+1}; \uw)$ be obtained by taking $\psi(\uz; \ux^0 \ux^1; \uw)$ and replacing every occurrence of $\theta(s_{m+1}(\uz; \ux^0 \ux^1; \uw))$ with $0$ and $x^1_{m+1}$, respectively.
    Define the term $t_{n+1}(\uz; \ux^0\ux^1; \uw)$ to be $s_{m+1}(\uz; \ux^0 \ux^1; \uw)$.
    Since the sentence $\forall x : x \not\in \VV \rightarrow \theta(x) = 0$ holds by definition in $T_\theta$, it is now easy to check that the formula
    \begin{align*}
        \exists \ux \in \VV : \exists \uz: & \bigwedge\nolimits^m_{i=1} x_i = s_i(\uz; \ux_{<i} \theta(\ux_{<i}); \uw) \wedge \bigwedge\nolimits_{j=1}^n t_j(\uz; \ux \theta(\ux); \uw) \not\in \VV \wedge \psi(\uz; \ux \theta(\ux); \uw)
    \end{align*}
    is equivalent, modulo $T_\theta$, to the disjunction
    \begin{align*}
        &\Big(\exists \ux \in \VV \!: \exists \uz: \bigwedge\nolimits^m_{i=1} x_i \!=\! s_i(\uz; \ux_{<i} \theta(\ux_{<i}); \uw) \\
        & \hspace{100pt}\wedge \bigwedge\nolimits_{j=1}^{n+1} t_j(\uz; \ux \theta(\ux); \uw) \!\not\in\! \VV \wedge \psi_0(\uz; \ux \theta(\ux); \uw)\Big) \\
        \vee\; & \Big(\exists \ux x_{m+1} \!\in\! \VV \!: \exists \uz :\!\! \bigwedge\nolimits^{\!m+1}_{i=1}\! x_i \!=\! s_i(\uz; \ux_{<i} \theta(\ux_{<i}); \uw) \\
        & \hspace{100pt}\wedge \bigwedge\nolimits_{j=1}^n\! t_j(\uz; \ux \theta(\ux); \uw) \!\not\in\! \VV \wedge \psi_1(\uz; \ux x_{m+1} \theta(\ux x_{m+1}); \uw)\Big).
    \end{align*}
    The formulas $\psi_0(\uz; \ux^0 \ux^1; \uw)$ and $\psi_1(\uz; \ux^0 x^0_{m+1} \ux^1 x^1_{m+1}; \uw)$ each contain at least one fewer occurrence of the symbol $\theta$ than $\psi(\uz; \ux^0 \ux^1; \uw)$.
    Notice that $\exists \uz : \varphi(\uz; \uw)$ is of the form (\ref{tag_lemma_fml}) with $m = n = 0$, and that we obviously have $\exists \uz : \varphi(\uz; \uw) \equiv \exists \ux \in \VV^0 : \exists \uz : \varphi(\uz; \uw)$.
    Hence, by finite recursion, $\exists \uz: \varphi(\uz; \uw)$ is equivalent to a finite disjunction $\bigvee_k \exists \ux_k \in \VV : \psi_k(\ux_k\theta(\ux_k); \uw)$, where the $\psi_k(\ux_k^0\ux_k^1; \uw)$ are $L$-formulas.
    Putting the disjunction inside the $\exists$-quantifier yields the desired result.
\end{proof}
\end{lemma}

\noindent With this simplification, one can already show the following:

\begin{remark}
\label{remark_reduct_exist_closed}
    If $(\mm, \theta)$ is an existentially closed model of $T^C_\theta$, then the vector space structure $(\VV, 0, +, (\lambda \cdot)_{\lambda\in K}, \theta)$ is an existentially closed model of $\TKvsTheC$.
\end{remark}

\subsection{Placeholder notation} \label{sec_first_char}

\noindent The goal of this section is to give a first indication of what our final characterization of existentially closed models of $T_\theta^C$ will look like. The formulas generated in Lemma \ref{lemma_fml} are generally hard to read.
For example, when taking $\varphi(z,w) := \theta^n(z) = w$, even after a few simplifications, we will end up with
$$
\psi((x^0_0, \dots, x^0_{n-1})(x^1_0, \dots, x^1_{n-1}); w) = \Big(w = 0\Big) \vee \Big(\bigwedge\nolimits_{i=0}^{n-2} x^1_{i} = x^0_{i+1} \wedge x^1_{n-1} = w\Big).
$$
The reason for this is simple.
Think of the tuples of variables $\ux^0$ and $\ux^1$ as placeholders for $\theta^0(\ux)$ and $\theta^1(\ux)$.
Since we do not have any placeholders for $\theta^n(\ux)$, we have to increase the length of $\ux^0$ and $\ux^1$ to talk about higher powers of $\theta$.
To fix this, we introduce placeholders for $\theta^n(\ux)$ for all $n \in \omega$.

\begin{definition}[Placeholder notation] \label{def_placeholder_notation}
    Let $\ux = (x_k : k \in \kk)$ be a tuple of variables.
    We define the \textbf{placeholder sequence} \hbox{$\uxvec := (x^i_k : k \in \kk, i \in \omega)$} to be a new tuple of variables.
    We call each $x^i_k$ a \textbf{placeholder variable} or a \textbf{placeholder} for $\theta^i(x_k)$.
    We furthermore define:
    \begin{enumerate}[(i)]
        \item $\ux^i := (x_k^i : k \in \kk)$ for each $i \in \omega$, and
        \item $\xvec_k := (x^i_k : i \in \omega)$ for each $k \in \kk$.
    \end{enumerate}
    We may sometimes write $(\ux^i : i \in \omega)$ or $(\xvec_k : k \in \kk)$ instead of $\uxvec$.
    For a singleton $x$, we similarly define $\xvec := (x^i : i \in \omega)$.
    If a formula $\psi(\uxvec; \uw)$ is given, we define:
    \begin{enumerate}[(i)]
        \setcounter{enumi}{2}
        \item $\psi_\theta(\ux; \uw) := \psi((\theta^i(x_k) : k \in \kk, i \in \omega); \uw)$.
    \end{enumerate}
\end{definition}

\noindent As one can see, the formula $\psi_\theta(\ux; \uw)$ is just the formula $\psi(\uxvec; \uw)$ with all placeholder variables $x^i_k$ replaced by $\theta^i(x_k)$; hence the name placeholder variable.
Even though $\uxvec$ consists of infinitely many variables, only finitely many of them actually appear in a single formula $\psi(\uxvec; \uw)$.

\begin{corollary}
    \label{corollary_equiv_exist_form}
    Modulo $T_\theta$, any existential $L_\theta$-formula $\phi(\uw)$ is equivalent to a formula of the form $\exists \ux \in \VV : \psi_\theta(\ux; \uw)$ for some $L$-formula $\psi(\uxvec; \uw)$.
\begin{proof}
    Let $\phi(\uw) = \exists \uz : \varphi(\uz; \uw)$, and apply Lemma \ref{lemma_fml} to obtain $\psi'(\ux^0\ux^1;\uw)$.
    Let $\psi(\uxvec; \uw)$ be $\psi'(\ux^0\ux^1;\uw)$, but interpreted as a formula in \hbox{$\uxvec = (\ux^i : i \in \omega)$}.
    We clearly have
    $$
    \exists \uz : \varphi(\uz; \uw) \quad \equiv \quad \exists \ux \in \VV : \psi(\ux\theta(\ux); \uw) \quad \equiv \quad \exists \ux \in \VV : \psi_\theta(\ux; \uw)
    $$
    modulo $T_\theta$.
\end{proof}
\end{corollary}

\noindent Notice that the formulas obtained from Corollary \ref{corollary_equiv_exist_form} above have exactly the same problem as the formulas obtained by Lemma \ref{lemma_fml}, and even more unnecessary variables.
However, while the formulas obtained by Lemma \ref{lemma_fml} will always have this problem, our \placeholderNotation{} at least \textit{allows} for neat formulas. With the help of Corollary \ref{corollary_equiv_exist_form}, we already obtain the first ``characterization'' of existentially closed models of $T^C_\theta$:

\begin{remark}
\label{corollary_exist_closed_0}
    Any $(\mm, \theta) \models T^C_\theta$ is existentially closed if and only if, for any $L(M)$-formula $\psi(\uxvec)$, we have $(\mm, \theta) \models \exists \ux \in \VV : \psi_\theta(\ux)$ whenever there is some model \hbox{$(\mm', \theta') \models T^C_\theta$} with $(\mm, \theta) \subseteq (\mm', \theta')$ and $(\mm', \theta') \models \exists \ux \in \VV : \psi_\theta(\ux)$.
\begin{proof}
    Since $T$ itself is model-complete, $\psi(\uxvec)$ is, modulo $T$, equivalent to an existential $L(M)$-formula.
    Therefore, $\psi_\theta(\ux)$ is itself an existential $L_\theta(M)$-formula.
    The direction from left to right follows easily from this, and the other direction is clear from Corollary \ref{corollary_equiv_exist_form} above.
\end{proof}
\end{remark}

\noindent The main reason why we cannot first-order axiomatize the ``characterization'' in Remark \ref{corollary_exist_closed_0} right now is that we do not know when the formula $\psi_\theta(\ux)$ is consistent in some extension $(\mm', \theta') \models T^C_\theta$.
It is obviously not sufficient to check whether $\psi(\uxvec)$ is consistent, since a fixed realization of $\psi(\uxvec)$, say $\uvvec = (v^i_k : 1 \leq k \leq m, i \in \omega)$, may not behave like a tuple of the form $(\theta^i(v_k) : 1 \leq k \leq m, i \in \omega)$:

\begin{example}
    Consider the formula $\psi(\xvec) := x^0 = 0 \wedge x^1 \neq 0$.
    This formula is obviously consistent in every model of $\TKvs$; however, the formula \hbox{$\psi_\theta(x) := \theta^0(x) = 0 \wedge \theta^1(x) \neq 0$} is inconsistent, as $x = 0$ implies $\theta(x) = 0$.
\end{example}

\noindent Whenever $\psi(\uxvec)$ implies an ``equation for $x_k^i$'' and some $x_k^j$ with $j > i$ appears, we might encounter similar problems as in the example above. In general, similar problems may arise whenever the formula $\psi(\uxvec)$ implies some finite disjunction of non-trivial linear dependencies in $\uxvec$ over $\VV$. However, one can also easily check that restricting ourselves to formulas that imply no finite disjunction of non-trivial linear dependencies over $\VV$ also does not work:

\begin{example}
    Choose $f \in \Kp{}$, set $d := \deg(f)$, and let $C$ be an algebraic kernel configuration with $\mipo(C) = f$.
    The formula 
    $$\psi(\xvec) := \sum\nolimits_{i=0}^d (f)_i \cdot x^i \neq 0$$ is consistent in $\TKvs$ and implies no finite disjunction of non-trivial linear dependencies.
    However, the formula \hbox{$\psi_\theta(x) := f[\theta](x) \neq 0$} is inconsistent in $T^C_\theta = T_\theta \cup \set{\Ker(f) = \VV}$.
\end{example}

\noindent Thus, although linear dependencies can cause problems in general, they cannot be avoided. The idea is to ``split'' all $L(M)$-formulas into a ``linearly independent part'' and some linear dependencies. More precisely, one can show that any sentence of the form $\exists \ux \in \VV : \psi_\theta(\ux)$, as in Remark \ref{corollary_exist_closed_0}, is equivalent to a disjunction of sentences of the form 
$$
\exists \ux \in \VV : \psi'_\theta(\ux) \wedge E(\ux),
$$
where $\psi'(\uxvec)$ is an $L(M)$-formula that implies no finite disjunction of non-trivial linear dependencies in $\uxvec$ over $\VV$, and $E(\ux)$ is a conjunction of $\LKThe(\VV)$-equations. Furthermore, one can ensure that these $\psi'(\uxvec)$ contain no unnecessary placeholders in the sense that if the corresponding $E(\ux)$ contains an equation of the form $\rho[\theta](x_k) = \sum_{l=1}^{k-1} \eta_l[\theta](x_l) + u$, then $x^i_k$ does not appear in $\psi'(\uxvec)$ for $i \geq \deg(\rho)$. This ``splitting'' essentially allows us to focus only on the conjunction of $\LKThe(\VV)$-equations $E(\ux)$. A crucial part of characterizing existentially closed models of $T^C_\theta$ will be figuring out the ``correct'' class of such conjunctions for each kernel configuration $C$.

\subsection{$C$-image-completeness}

\noindent So far, we have only simplified existential $L_\theta$-formulas modulo $T_\theta$, and not modulo $T^C_\theta$, where $T^C_\theta = T_\theta \cup \set{\text{``$\theta$ is a $C$-endomorphism''}}$.
One could now try to simplify existential $L_\theta$-formulas modulo $T^C_\theta$.
However, the theory $T^C_\theta$ itself does not enforce enough structure to allow for any significant simplifications.
Instead, we show that the endomorphisms of existentially closed models of $T^C_\theta$ satisfy a property that we call $C$-image-completeness.
Later on, we will simplify existential $L_\theta$-formulas, or rather existential $L_\theta(M)$-sentences, modulo the theory $T_\theta \cup \set{\text{``$\theta$ is $C$-image-complete''}}$.

\begin{definition} \label{def_C_image_comple}
    We say that an endomorphism $\theta \colon \VV \to \VV$ is \textbf{$\mathbf{C}$-image-complete} if it is a $C$-endomorphism and $\Image(f^{C+1}) = \Image(f^C)$ holds for all $f \in \Kp{C<\infty}$.
    We may call a model $(\mm, \theta) \models T_\theta$ \textbf{$\mathbf{C}$-image-complete} if the endomorphism $\theta$ is $C$-image-complete.
\end{definition}

\noindent In the algebraic case, we obtain $C$-image-completeness for free:

\begin{lemma}
    \label{remark_alg_c_complete}
    If $C \in \Ccalg$, then any $C$-endomorphism $\theta \colon \VV \to \VV$ is $C$-image-complete.
\begin{proof}
    Setting $F := \Kp{0<C<\infty}$, Observation \ref{obs_mipo_prod} yields
    $$
    \VV = \Ker(\mipo(C)) = \Ker\Big(\prod\nolimits_{f\in F} f^C\Big) = \bigoplus\nolimits_{f\in F} \Ker(f^C).
    $$
    For any $g \in \Kp{}\setminus F$, we have $C(g) = 0$ and therefore $g \nmid \mipo(C)$.
    Hence the map $g[\theta]$ is an automorphism on $\Ker(\mipo(C)) = \VV$ by (v) of Corollary \ref{corollary_polynomials_kernel_facts}, which implies the equality $\Image(g^0) = \VV = \Image(g^{0+1})$.
    Now let $f, g \in F$ with $f \neq g$ be given.
    Similarly, the map $f[\theta]_{\restriction \Ker(g^C)}$ is an automorphism of $\Ker(g^C)$, which implies $\Ker(g^C) = f^n[\theta](\Ker(g^C))$ for any $n \geq 0$.
    We now obtain
    \begin{align*}
        f^C[\theta](\VV) = \bigoplus_{g \in F} f^C[\theta](\Ker(g^C))
        = \bigoplus_{\substack{g \in F \\ g \neq f}} \Ker(g^C)
        = \bigoplus_{g \in F} f^{C+1}[\theta](\Ker(g^C))
        = f^{C+1}[\theta](\VV),
    \end{align*}
    since $f^{C+1}[\theta](\Ker(f^C)) = \set{0} = f^{C}[\theta](\Ker(f^C))$.
    Thus $\Image(f^C) = \Image(f^{C+1})$, as desired.
\end{proof}
\end{lemma}

\noindent If $C$ is a transcendental kernel configuration, we need existential closedness in order to obtain $C$-image-completeness:

\begin{lemma} \label{lemma_ext_image}
    Let $C \in \Cctrans$ be a transcendental kernel configuration, let $(\mm, \theta)$ be a model of $T^C_\theta$, and let $f \in \Kp{C<\infty}$ be a monic irreducible polynomial over $K$ with $C(f) < \infty$.
    Given any element \hbox{$v_0 \in \Image(f^{C})\setminus \Image(f^{C+1})$}, there is an extension $(\mm', \theta') \supseteq (\mm, \theta)$ with $(\mm', \theta') \models T^C_\theta$ and $v_0 \in \Image(f^{C+1}[\theta'])$.
\begin{proof}
    Write $f = \sum\nolimits_{i=0}^d (f)_i \cdot X^i$ with $d := \deg(f)$.
    We are now going to construct $(\mm', \theta')$.
    For this, choose $\mm' \succ \mm$ such that $\dim(\VV'/\VV) \geq \aleph_0$.
    Next, choose some $u_0 \in \VV$ with $f^C[\theta](u_0) = v_0$ and a set $B = \set{v^i : 0 \leq i < d} \subset \VV'$ linearly independent over $\VV$.
    We now extend $\theta$ to an endomorphism $\theta'$ of $\spanA{B \VV}{K}$ by setting
    $$
    \theta'(v^i) := \begin{cases}
        v^{i+1} & \text{if $i < d-1$}, \\
        u_0 -\sum\nolimits_{j = 0}^{d-1} (f)_j \cdot v^j & \text{if $i = d-1$}.
    \end{cases}
    $$
    Notice the following:
    \begin{enumerate}[(i)]
        \item Any $v \in \spanA{B \VV}{K}$ can be written as $\rho[\theta'](v^0) + u_1$ with \hbox{$\deg(\rho) < \deg(f) = d$} and $u_1 \in \VV$.
        This follows directly from the definition of $B$ and $\theta'$.
        \item Given any polynomial $\rho \in K[X]$, we have $\rho[\theta'](v^0) = r[\theta'](v^0) + \chi[\theta](u_0)$ for the unique polynomials $\chi, r \in K[X]$ with $\rho = \chi \cdot f + r$ and $\deg(r) < \deg(f)$.
        This follows from the \remainderRule{} (Lemma~\ref{lemma_remainder_rule}), since we have $f[\theta'](v^0) = u_0$ by the definition of $\theta'$.
    \end{enumerate}
    In the following, we will show that $\theta'$ is a $C$-endomorphism.
    For this, we verify the equality $\Ker(g^C[\theta']) = \Ker(g^{C+1}[\theta'])$ for all $g \in \Kp{C<\infty}$.
    Since the inclusion ``$\subseteq$'' is clear, we only need to verify $g^{C}[\theta'](v) = 0$ for each $v \in \spanA{B \VV}{K}$ with $g^{C+1}[\theta'](v) = 0$.
    Let such a polynomial $g \in \Kp{C<\infty}$ and element $v\in \spanA{B \VV}{K}$ with $g^{C+1}[\theta'](v) = 0$ be given.
    By (i), we have $v = \rho[\theta'](v^0) + u_1$ with $\deg(\rho)< \deg(f)$ and $u_1 \in \VV$.
    With this, we obtain
    \begin{equation} \label{eq_1}
        0 = g^{C+1}[\theta'](v) = g^{C+1}[\theta'](\rho[\theta'](v^0) + u_1) = (g^{C+1} \cdot \rho)[\theta'](v^0) + g^{C+1}[\theta](u_1)
    \end{equation}
    and distinguish between the cases $g \neq f$ and $g = f$.

    First, assume that $g \neq f$.
    Recall that $\deg(\rho) < \deg(f)$ and that both $f$ and $g$ are irreducible, so $f \nmid g^{C+1} \cdot \rho$ unless $\rho = 0$.
    Assuming $\rho \neq 0$, we obtain $r \neq 0$ for the unique polynomials $\chi$ and $r$ with $\chi \cdot f + r = g^{C+1} \cdot \rho$ and $\deg(r) < \deg(f)$.
    With (\ref{eq_1}) above, (ii), and $\theta'^i(v^0) = v^i$ for $i < d$, we obtain
    \begin{align*}
        0 = g^{C+1}[\theta'](v)
        = r[\theta'](v^0) + \chi[\theta](u_0) + g^{C+1}[\theta](u_1)
        = \!\sum\nolimits_{i=0}^{d-1} (r)_i \cdot v^i \!+\! \underbrace{\chi[\theta](u_0) + g^{C+1}[\theta](u_1)}_{\in \VV},
    \end{align*}
    where $r = \sum_{i=0}^{d-1} (r)_i \cdot X^i$.
    This contradicts the linear independence of $v^0, \dots, v^{d-1}$ over $\VV$, so we must have $\rho = 0$.
    This means $v = u_1 \in \Ker(g^{C+1}[\theta]) \subseteq \VV$.
    Since $\theta$ is a $C$-endomorphism, we conclude $g^{C}[\theta'](v) = g^{C}[\theta](u_1) = 0$.

    Now, assume that $g = f$.
    Plugging $g = f$ and $f[\theta'](v^0) = u_0$ into (\ref{eq_1}) yields
    $$
    0 = (f^{C+1} \cdot \rho)[\theta'](v^0) + f^{C+1}[\theta](u_1) = (f^{C} \cdot \rho)[\theta](u_0) + f^{C+1}[\theta](u_1).
    $$
    This implies $(f^{C} \cdot \rho)[\theta](u_0) = f^{C+1}[\theta](-u_1) \in \Image(f^{C+1})$.
    Since $f$ is an irreducible polynomial and $\deg(\rho) < \deg(f)$, we must have $\gcd(f^{C} \cdot \rho, f^{C+1}) = f^C$, unless $\rho = 0$.
    Assuming $\rho \neq 0$, we can use \gcdRule{} with $\gcd(f^{C} \cdot \rho, f^{C+1}) = f^C$ to show
    $$
    v_0 = f^C[\theta](u_0) =  \chi_1[\theta](f^{C+1}[\theta](-u_1)) + \chi_2[\theta](f^{C+1}[\theta](u_0)) \in \Image(f^{C+1}[\theta]),
    $$
    for some $\chi_1, \chi_2 \in K[X]$.
    We must have $\rho = 0$, as \hbox{$v_0 \in \Image(f^C[\theta]) \setminus \Image(f^{C+1}[\theta])$}.
    Conclude $f^{C}[\theta'](v) = 0$ as in the case where $g \neq f$.

    Now that we have shown that $\theta'$ is a $C$-endomorphism on $\spanA{B \VV}{K}$, we can use our \standartConstruction{} to extend $\theta'$ to a $C$-endomorphism defined on all of $\VV'$. After this step, it is clear that we still have $f^{C+1}[\theta'](v^0) = f^C[\theta'](f[\theta'](v^0)) = f^C[\theta](u_0) = v_0$, and therefore $v_0 \in \Image(f^{C+1}[\theta'])$.
\end{proof}
\end{lemma}

\begin{corollary}\label{corollary_C_image_compl}
    If $(\mm, \theta) \models T^C_\theta$ is existentially closed, then $\theta$ is $C$-image-complete.
\begin{proof}
    By Lemma \ref{remark_alg_c_complete}, we only need to consider the case where $C \in \Cctrans$ is transcendental.
    Fix some $f \in \Kp{C<\infty}$ and choose some $v \in \Image(f^C[\theta])$.
    By Lemma \ref{lemma_ext_image}, there is an extension $(\mm', \theta') \supseteq (\mm, \theta)$ with $(\mm', \theta') \models T^C_\theta$ and $v \in \Image(f^{C+1}[\theta'])$.
    By existential closedness of $(\mm, \theta)$, we must already have $v \in \Image(f^{C+1}[\theta])$.
    The remaining inclusion $\Image(f^{C+1}) \subseteq \Image(f^{C})$ is trivial.
\end{proof}
\end{corollary}

\begin{notation}
    For any $\rho \in K[X] \setminus \set{0}$, we let $\Fac(\rho) := \set{f \in \Kp{} : f \mid \rho}$ denote the set of all irreducible factors of $\rho$.
    As a convention, we set $\Fac(0) = \varnothing$.
\end{notation}

\noindent Also, recall that we defined $F^C := \prod\nolimits_{f\in F} f^{C(f)}$ for all finite subsets \hbox{$F \subseteq \Kp{C<\infty}$}, where $\Kp{C<\infty} = \set{f \in \Kp{} : C(f) < \infty}$.

\begin{lemma} \label{lemma_decomposition}
    If $\theta \colon \VV \to \VV$ is $C$-image-complete and $F \subseteq \Kp{0<C<\infty}$ is a finite set, then the following hold:
    \begin{enumerate}[(i)]
        \item $\Image(\eta) = \Image(F^C)$ for any $\eta \in K[X] \setminus \set{0}$ with $F^C \mid \eta$ and $\Fac(\eta) \subseteq \Kp{C=0} \cup F$.
        \item $\VV = \Image(F^C) \oplus \Ker(F^C) = \Image(F^C) \oplus \bigoplus\nolimits_{f\in F} \Ker(f^C)$.
        \item $\Image(F^C) = \Image(F'^C) \oplus \bigoplus\nolimits_{f \in F' \setminus F} \Ker(f^C)$ for any finite $F'$ with \hbox{$F \subseteq F' \subseteq \Kp{0<C<\infty}$}.
        \item $\Image(g^C) = \Image(F^C) \oplus \bigoplus\nolimits_{f \in F \setminus \set{g}} \Ker(f^C)$ for any $g \in F$.
    \end{enumerate}
    If $C$ is algebraic, we in particular obtain $\VV = \bigoplus\nolimits_{f \in \Kp{0<C<\infty}} \Ker(f^C)$, as already described in Observation \ref{obs_mipo_prod}.
    \begin{proof}
        We only prove (i) and (ii), since (iii) and (iv) follow directly from point (ii) together with (v) of Corollary \ref{corollary_polynomials_kernel_facts}.
        \begin{enumerate}[(i)]
            \item
            By assumption, we can write $\eta = \prod\nolimits_{f\in \Fac(\eta)} f^{n_f}$ with $n_f \geq C(f)$ for all $f \in \Fac(\eta)$ (without loss $\eta$ is monic).
            By $C$-image-completeness we obtain $\Image(f^{n_f}) = \Image(f^C)$, and especially $\Image(f^{n_f}) = \VV$ in the case where $C(f) = 0$. With (vii) of Corollary \ref{corollary_polynomials_kernel_facts}, we obtain
            $$
            \Image(\eta) = \bigcap\nolimits_{f\in \Fac(\eta)} \Image(f^{n_f}) = \bigcap\nolimits_{f\in F} \Image(f^C) = \Image(F^C).
            $$
            \item As always, $\Ker(F^C) = \bigoplus\nolimits_{f\in F} \Ker(f^C)$ follows from (iii) and (iv) of Corollary \ref{corollary_polynomials_kernel_facts}.
            Let $v \in \VV$ be given.
            By $C$-image-completeness and (i), we know that $F^C[\theta](v) \in \Image(F^{2C})$, where
            $$
            F^{2C} = \prod\nolimits_{f\in F} f^{2C(f)} = (F^C)^2.
            $$
            Hence there is some $u \in \VV$ with $F^{2C}[\theta](u) = F^{C}[\theta](v)$.
            Thus \hbox{$v - F^C[\theta](u) \in \Ker(F^C)$}, and therefore $v \in \Image(F^C) + \Ker(F^C)$.

            Now assume that $v \in \Image(F^C) \cap \Ker(F^C)$.
            Since $v \in \Image(F^C)$, there is some $u \in \VV$ with $F^C[\theta](u) = v$.
            Since $v \in \Ker(F^C)$, we obtain $u \in \Ker(F^{2C})$.
            Since $\theta$ is a $C$-endomorphism, we have $\Ker(F^{2C}) = \Ker(F^C)$.
            Therefore $v = F^C[\theta](u) = 0$. \qedhere
        \end{enumerate}
    \end{proof}
\end{lemma}

\subsection{An optimal ring of definable scalars}  \label{sec_rc_module}
\noindent This section deals with the $\LKThe$-definable endomorphisms of $\VV$ in existentially closed (or rather $C$-image-complete) models of $T^C_\theta$.
Let any model $(\VV, \theta) \models \TKvsTheC$ be given.
Notice the following:
\begin{enumerate}[(i)]
    \item We can define a $K[X]$-module structure on $\VV$ by setting $\rho \cdot v := \rho[\theta](v)$ for all $\rho \in K[X]$ and $v \in \VV$.
    \item Let $C \in \Ccalg$ be algebraic with $\mipo(C) = \rho$. In this case, we can also define a $K[X]/(\rho)$-module structure on $\VV$ as $\eta_1[\theta](v) = \eta_2[\theta](v)$ if $\eta_1 - \eta_2$ is a multiple of $\rho$.
\end{enumerate}
\noindent Given $(\VV,\theta)$ where $\theta$ is a $C$-endomorphism $C \in \Ccalg$ with $\mipo(C) = \rho$, it clearly makes more sense to treat $\VV$ as a $K[X]/(\rho)$-module than to treat $\VV$ as a $K[X]$-module.
This raises the following question:
\begin{itemize}
    \item [] Given an existentially closed $(\VV, \theta) \models \TKvsTheC$, what is the ``optimal'' ring $R_C$ of definable endomorphisms on $\VV$?
\end{itemize}
In order to answer this question, one must decide what ``optimal'' means. We want the following three properties to hold:
\begin{enumerate}[(i)]
    \item The ring $R_C$ should only depend on the kernel configuration $C$.
    \item Every $\LKThe$-definable endomorphism should be in $R_C$.
    \item Given an existentially closed  $(\VV, \theta) \models \TKvsTheC$ and $r \neq r' \in R_C$, there should be some $v \in \VV$ with $r(v) \neq r'(v)$.
\end{enumerate}
So what we are actually doing here is, in the context of model theory of modules, finding the ring of definable scalars (see Section 6.1 of \cite{Pre09}) for the class of existentially closed models of $\TKvsTheC$ treated as $K[X]$-modules.
The main ingredient in finding $R_C$ is the $C$-image-completeness from the previous section and the decomposition $\VV = \Image(F^C) \oplus \Ker(F^C)$ (Lemma \ref{lemma_decomposition}) that it implies.

\begin{lemma} \label{theorem_new_functions}
    In the theory \hbox{$\TKvsThe \cup \set{\text{``$\theta$ is $C$-image-complete''}}$}, the following endomorphisms are definable by pp-formulas (i.e., formulas of the form $\psi(x, y) := \exists \uz : E(xy;\uz)$, where $E(xy; \uz)$ is a conjunction of $K[X]$-equations) in the language of $K[X]$-modules:
    \begin{enumerate}[(i)]
        \item For $F \subseteq \Kp{0<C<\infty}$ finite, we define the \textbf{projection to the image of $\bm{F^C[\theta]}$} by
        $$
        \pi_{\Image(F^C)}(x) := \text{``the unique $u \!\in\! \Image(F^C)$ for which there is $v \in \Ker(F^C)$ with $x \!=\! u \!+\! v$''.}
        $$
        \item For $F \subseteq \Kp{0<C<\infty}$ finite, we define the \textbf{projection to the kernel of $\bm{F^C[\theta]}$} by
        $$
        \pi_{\Ker(F^C)}(x) := \text{``the unique $v \!\in\! \Ker(F^C)$ for which there is $u \in \Image(F^C)$ with $x \!=\! u \!+\! v$''.}
        $$
        We clearly have $\pi_{\Ker(F^C)} = 1[\theta] - \pi_{\Image(F^C)}$.
        \item For every monic polynomial $\eta \in K[X]$ with $\Fac(\eta) \subseteq \Kp{C<\infty}$, we define the \textbf{pseudo-inverse of $\bm{\eta[\theta]}$} by
        $$
        \eta[\theta]^{-1}(x) := \text{``the unique $u \in \Image(\Fac(\eta)^C)$ with $\eta[\theta](u) = \pi_{\Image(\Fac(\eta)^C)}(x)$''.}
        $$
        Notice that $\Fac(\eta)^C = (\Fac(\eta) \cap \Kp{0<C<\infty})^C$.
        In practice, we will also use $\eta[\theta]^{-1}$ for (non-zero) non-monic polynomials by setting $\eta[\theta]^{-1} := \lambda^{-1} \cdot (\eta/\lambda)[\theta]^{-1}$ for the leading coefficient $\lambda$ of $\eta$.
    \end{enumerate}
\begin{proof}
    The well-definedness of the functions in (i) and (ii) follows from (ii) of Lemma \ref{lemma_decomposition}, since we have $\VV = \Image(F^C) \oplus \Ker(F^C)$.
    We now show the well-definedness of the function in (iii).
    Let $v \in \VV$ be given in $(\VV, \theta) \models \TKvsThe \cup \set{\text{``$\theta$ is $C$-image-complete''}}$, and define $F := \Fac(\eta) \cap \Kp{0<C<\infty}$.
    By definition, we have $\pi_{\Image(F^C)}(v) \in \Image(F^C)$.
    Since we also have \hbox{$\Image(F^C) = \Image(F^C \cdot \eta)$} by (i) of Lemma \ref{lemma_decomposition}, there is some element $u \in \Image(F^C)$ with $\eta[\theta](u) = \pi_{\Image(F^C)}(v)$.
    If there were another $u' \in \Image(F^C)$ with $\eta[\theta](u') = \pi_{\Image(F^C)}(v)$, then $u-u' \in \Ker(\eta)$.
    Since $\theta$ is a $C$-endomorphism, we have $\Ker(\eta) \subseteq \Ker(F^C)$, and therefore $u-u' \in \Image(F^C) \cap \Ker(F^C) = \set{0}$.
    We conclude that this $u$ is unique.

    The pp-formula $\exists z_1z_2: x = F^C \cdot z_1 + z_2 \wedge y = F^C \cdot z_1 \wedge F^C \cdot z_2 = 0$ defines the function in (i).
    Here $F^C \cdot z_1$ translates to $F^C[\theta](z_1)$ in our usual language $\LKThe$.
    It is easy to find the corresponding pp-formulas in the other cases.
    All these functions must be $K[X]$-endomorphisms, as $K[X]$ is commutative, and therefore also $K$-endomorphisms, i.e., endomorphisms of the $K$-vector space $\VV$.
\end{proof}
\end{lemma}

\noindent Let $(\VV, \theta)$ be a vector space with a $C$-image-complete endomorphism.
The projections from Lemma \ref{theorem_new_functions} are the canonical projections from $\VV = \Image(F^C) \oplus \Ker(F^C)$ to $\Image(F^C)$ and $\Ker(F^C)$, respectively.
In particular, we have $\pi_{\Image(F^C)}(v) = v$ for any $v \in \Image(F^C)$ and $\pi_{\Ker(F^C)}(v) = v$ for any $v \in \Ker(F^C)$.

\begin{notation}
    Given $f \in \Kp{0<C<\infty}$, we write $\pi_{\Ker(f^C)} := \pi_{\Ker(\set{f}^C)}$.
\end{notation}

\begin{lemma}[Projection Identity] \label{lemma_projection_rule}
For any finite set $F \subseteq \Kp{0<C<\infty}$ of irreducible polynomials $f$ over $K$ with $0 < C(f) < \infty$, we have
$$
1[\theta] = \pi_{\Image(F^C)} + \pi_{\Ker(F^C)} = \pi_{\Image(F^C)} + \sum\nolimits_{f\in F} \pi_{\Ker(f^C)}.
$$
If $C$ is algebraic, we especially obtain $1[\theta] = \sum\nolimits_{f\in \Kp{0<C<\infty}} \pi_{\Ker(f^C)}$.
\end{lemma}

\begin{remark} \label{rem_no_cancel}
Note that the notation $\eta[\theta]^{-1}$ might be misleading, as in general we only have
$$
    \eta[\theta] \circ \eta[\theta]^{-1} = \pi_{\Image(\Fac(\eta)^C)}
$$
and not $\eta[\theta] \circ \eta[\theta]^{-1} = \Id$, as one might expect.
This means that $\eta[\theta]^{-1}$ is not an actual inverse of $\eta[\theta]$, but it is as close as we can get.
If we take $\eta = f^C$ for some $f \in \Kp{0<C<\infty}$ and $(\VV, \theta) \models \TKvsThe \cup \set{\text{``$\theta$ is $C$-image-complete''}}$ with $\Ker(f^C)$ infinite, then the actual preimage of $0$ under $f^C[\theta]$ is $\Ker(f^C)$, whereas for any $v \in \Ker(f^C) \setminus \set{0}$, the preimage of $v$ under $f^C[\theta]$ is empty.
This means that there is no reasonable way to define a function $f^C[\theta]^{-1}$ on $\Ker(f^C)$.
\end{remark}

\begin{theorem} \label{theorem_r_c_def}
    Let $R_C$ be the set of all endomorphisms that are definable in the theory \hbox{$\TKvsThe \cup \set{\text{``$\theta$ is $C$-image-complete''}}$} and $\set{+, \circ}$-generated by
    \begin{align*}
        \set{\rho[\theta] : \rho \in K[X]} &\cup \set{\eta[\theta]^{-1} : \eta \text{ monic with }\Fac(\eta) \subseteq \Kp{C<\infty}}.
    \end{align*}
    The structure $(R_C, 0[\theta], 1[\theta], +, \circ)$ is a unitary commutative ring with the same characteristic as $K$.
    Furthermore, all elements of $R_C$ can be defined by a pp-formula in the language of $K[X]$-modules.
\begin{proof}
    Clearly, $(R_C, 0[\theta], 1[\theta], +, \circ)$ is a subring of the ring of all definable endomorphisms in \hbox{$\TKvsThe \cup \set{\text{``$\theta$ is $C$-image-complete''}}$}.
    Since $1[\theta] = \theta^0$ is the identity, it follows immediately that $R_C$ is unitary.
    Because we work with $\TKvsThe$, it is clear that $R_C$ has the same characteristic as the field $K$.
    By Lemma \ref{theorem_new_functions}, all generating functions are given by pp-formulas in the language of $K[X]$-modules.
    This implies that all elements $r \in R_C$ are given by a pp-formulas, and are therefore $K[X]$-endomorphisms.
    One can easily check that any such pp-formula commutes with every $K[X]$-linear function (see e.g. the proof of 1.1.7 in \cite{Pre09}), and hence $R_C$ is commutative.
\end{proof}
\end{theorem}

\noindent By Remark \ref{rem_no_cancel}, we immediately see that $\pi_{\Image(F^C)} = F^C[\theta] \circ F^C[\theta]^{-1}$ is in $R_C$ for any finite $F \subseteq \Kp{0<C<\infty}$.
With this, we also see that $\pi_{\Ker(F^C)} = 1[\theta] - \pi_{\Image(F^C)}$ is in $R_C$ for any such $F$.

\begin{notation}
    Instead of $(R_C, 0[\theta], 1[\theta], +, \circ)$, we may write $(R_C, 0, \Id, +, \circ)$ or $(R_C, 0, 1, +, \cdot)$, depending on whether we treat $R_C$ as a ring of endomorphisms or as a purely algebraic object.
\end{notation}

\noindent Since the ring $R_C$ depends on the kernel configuration $C$, it is clear that it also depends on the field $K$, because the kernel configuration also depends on $K$.

Also note that the elements of $R_C$ are, by definition, definable functions in the theory $\TKvsThe \cup \set{\text{``$\theta$ is $C$-image-complete''}}$.
By this, we mean that the elements of $R_C$ are equivalence classes of $\LKThe$-formulas modulo the theory \hbox{$\TKvsThe \cup \set{\text{``$\theta$ is $C$-image-complete''}}$} that define an endomorphism in every model of $\TKvsThe \cup \set{\text{``$\theta$ is $C$-image-complete''}}$.
So, in order to prove $r = r'$, we need to show
$$
    r^{(\VV, \theta)} = r'^{(\VV, \theta)} \quad \text{for all $(\VV, \theta) \models \TKvsThe \cup \set{\text{``$\theta$ is $C$-image-complete''}}$},
$$
and, in order to prove $r \neq r'$, we need to find $(\VV, \theta) \models \TKvsThe \cup \set{\text{``$\theta$ is $C$-image-complete''}}$ with
$
    r^{(\VV, \theta)} \neq r'^{(\VV, \theta)}.
$

\begin{lemma} \label{lemma_substructure_R_C}
    Let $(\VV, \theta)$ be $C$-image-complete, and let $U \subseteq \VV$ be non-empty and closed under addition and $R_C$, i.e., for any $r \in R_C$ and $u \in U$, we have $r^{(\VV, \theta)}(u) \in U$.
    Then $(U, \theta_{\restriction U})$ is $C$-image-complete, and for any $r \in R_C$ and $u \in U$, we obtain
    $$
    r^{(\VV, \theta)}(u) = r^{(U, \theta_{\restriction U})}(u).
    $$
\begin{proof}
    The equation $\Ker(\rho[\theta_{\restriction U}]) = \Ker(\rho[\theta]) \cap U$ is clear for every polynomial $\rho \in K[X]$, and $\Image(f^{C+i}[\theta_{\restriction U}]) = \Image(f^{C+i}[\theta]) \cap U$ follows for every $f \in \Kp{C<\infty}$ and $i \in \set{0, 1}$, as $U$ is closed under $f^{C+i}[\theta]^{-1}$.
    Therefore, $(U, \theta_{\restriction U})$ is $C$-image-complete by definition.
    Now $(r^{(\VV, \theta)})_{\restriction U} = r^{(U, \theta_{\restriction U})}$ follows, since each $r \in R_C$ is defined by a pp-formula.
\end{proof}
\end{lemma}

\noindent Given $(\mm, \theta) \models T_\theta \cup \set{\text{``$\theta$ is $C$-image-complete''}}$, the structure $(\VV, \theta)$ is obviously a model of $\TKvsThe \cup \set{\text{``$\theta$ is $C$-image-complete''}}$, and the structure $(\VV, \theta)$ is clearly definable in $(\mm, \theta)$.
Thus, the map $r^{(\VV, \theta)}$ is also definable in $(\mm, \theta)$ for every $r \in R_C$.

\begin{remark} \label{rem_diff_endo}
    If $(\mm, \theta) \models T^C_\theta$ is existentially closed, then $r$ and $r'$ define different endomorphisms on $\VV$ for $r \neq r' \in R_C$.
\begin{proof}
    By Remark \ref{remark_reduct_exist_closed}, we know that $(\VV, \theta)$ is an existentially closed model of $\TKvsTheC$.
    Since $r \neq r'$, there are a model $(\VV', \theta') \models \TKvsThe \cup \set{\text{``$\theta$ is $C$-image-complete''}}$ and some $v \in \VV'$ with $r(v) \neq r'(v)$.
    Now it is easy to check that
    $
    (\VV \oplus \VV', \theta \oplus \theta')
    $
    is $C$-image-complete and extends both $(\VV, \theta)$ and $(\VV', \theta')$.
    Using Lemma \ref{lemma_substructure_R_C}, one can easily verify that we have $r((0, v)) \neq r'((0, v))$ in
    $
    (\VV \oplus \VV', \theta \oplus \theta').
    $
    By existential closedness of $(\VV, \theta)$, we can also find some $u \in \VV$ with $r(u) \neq r'(u)$.
\end{proof}
\end{remark}

\noindent With the above we have seen that $R_C$ has two of the ``optimality'' properties, as discussed at the beginning of this section.
The proof of the remaining ``optimality'' property, namely that given an existentially closed $(\VV, \theta) \models \TKvsTheC$, every $\LKThe$-definable endomorphism is in $R_C$, requires a better understanding of existentially closed models of $\TKvsTheC$ and will be proven in a future paper.

The goal for the rest of this section is to describe $(R_C, 0, 1, +, \cdot)$ algebraically as well as possible and to establish some basic facts for the upcoming chapters.
We will later show the following in Corollary \ref{corollary_R_C_as_ring}, along with some other examples:
\begin{enumerate}[(i)]
    \item $(R_C, 0, 1, +, \cdot) \simeq (K[X]/(\mipo(C)), 0, 1, +, \cdot)$ for all $C \in \Ccalg$.
    This also implies that our ring $(R_C, 0, \Id, +, \circ)$ is a field for all algebraic kernel configurations $C$ with $\mipo(C)$ irreducible.
    \item $(R_{C_\infty}, 0, 1, +, \cdot) \simeq (K[X], 0, 1, +, \cdot)$ for the unique kernel configuration $C_\infty \in \Cctrans$ that satisfies \hbox{$C_\infty(f) = \infty$} for all $f \in \Kp{}$.
    \item $(R_{C_0}, 0, 1, +, \cdot) \simeq (K(X), 0, 1, +, \cdot)$ for the unique kernel configuration $C_0 \in \Cctrans$ with $C_0(f) = 0$ for all $f \in \Kp{}$.
\end{enumerate}
Points (ii) and (iii) can already be guessed from the definition of $(R_C, 0, \Id, +, \circ)$.
The first step is to give some computation rules for the multiplication in $R_C$:

\begin{lemma}
    \label{lemma_communative}
    The following equations hold in $(R_C, 0, \Id, +, \circ)$:
    \begin{enumerate}[(i)]
        \item $\rho_1[\theta] \circ \rho_2[\theta] = (\rho_1 \cdot \rho_2)[\theta]$.
        \item $\rho[\theta] \circ \pi_{\Image(F^C)} = \rho[\theta]$ if $F^C \mid \rho$.
        \item $\rho[\theta] \circ \pi_{\Ker(F^C)} = r[\theta] \circ \pi_{\Ker(F^C)}$ for the remainder $r \in K[X]$ of $\rho$ upon division by $F^C$.
        \item $\rho[\theta] \circ \eta[\theta]^{-1} = \frac{\rho}{\gcd(\rho, \eta)}[\theta] \circ \frac{\eta}{\gcd(\rho, \eta)}[\theta]^{-1} \circ \pi_{\Image(\Fac(\eta)^C)}$.
        \item $\pi_{\Image(F_1^C)} \circ \pi_{\Image(F_2^C)} = \pi_{\Image((F_1 \cup F_2)^C)}$.
        \item $\pi_{\Image(F_1^C)} \circ \pi_{\Ker(F_2^C)} = \pi_{\Ker((F_2 \setminus F_1)^C)}$.
        We also obtain
        $
        \pi_{\Image(F_1^C)} \circ \pi_{\Ker(F_2^C)} = \chi[\theta] \circ \pi_{\Ker(F_2^C)}
        $
        for some polynomial $\chi \in K[X]$.
        This $\chi$ depends only on $F_1^C$ and $F_2^C$, and it is $1$ if $F_1 \cap F_2 = \varnothing$.
        \item $\pi_{\Image(F^C)} \circ \eta[\theta]^{-1} = \eta[\theta]^{-1}$ if $F \subseteq \Fac(\eta)$.
        \item $\pi_{\Ker(F_1^C)} \circ \pi_{\Ker(F_2^C)} = \pi_{\Ker((F_1 \cap F_2)^C)}$.
        We also obtain
        $
        \pi_{\Ker(F_1^C)} \circ \pi_{\Ker(F_2^C)} = \chi[\theta] \circ \pi_{\Ker(F_2^C)}
        $
        for some $\chi \in K[X]$ that depends only on $F_1^C$ and $F_2^C$.
        \item $\pi_{\Ker(F^C)} \circ \eta[\theta]^{-1} = \chi[\theta] \circ \pi_{\Ker(F^C)}$ for some $\chi \in K[X]$ with $\chi = 0$ if $F \subseteq \Fac(\eta)$.
        This $\chi$ depends only on $\eta$ and $F^C$.
        \item $\eta_1[\theta]^{-1} \circ \eta_2[\theta]^{-1} = (\eta_1 \cdot \eta_2)[\theta]^{-1}$.
    \end{enumerate}
    Here all $\rho$'s lie in $K[X]$, all $F$'s are finite subsets of $\Kp{0<C<\infty}$, and all $\eta$'s are monic polynomials with $\Fac(\eta) \subseteq \Kp{C<\infty}$.
\begin{proof}
    Throughout this proof, fix some $C$-image-complete $(\VV, \theta)$.
    We omit (i), as it is just Lemma \ref{lemma_add_mult_endo}, and (ii), as it is trivial.
    We also omit (vii), as it is clear from the definition of $\eta[\theta]^{-1}$, and both (viii) and (ix), since they follow from the \projectionRule{} and the previous points.
    \begin{enumerate}[(i)]
        \setcounter{enumi}{2}
        \item Given $u \in \Image(F^C)$ and $v \in \Ker(F^C)$, we have $F^C[\theta](v) = 0$, so by the \remainderRule{}:
        $$
        \rho[\theta](\pi_{\Ker(F^C)}(u + v)) = \rho[\theta](v) = r[\theta](v) = r[\theta](\pi_{\Ker(F^C)}(u + v))
        $$
        for the unique polynomials $r, \chi \in K[X]$ with $\deg(r) < \deg(F^C)$ and $\rho = r + \chi \cdot F^C$.
        \setcounter{enumi}{4}
        \item Let $i \in \set{1, 2}$.
        By Lemma \ref{theorem_new_functions}, $\pi_{\Image(F_i^C)} \colon \VV \to \Image(F_i^C)$ is the canonical projection to $\Image(F_i^C)$.
        Using (ii) and (iii) of Lemma \ref{lemma_decomposition}, we obtain:
        \begin{align*}
            \VV &= \Image((F_1 \cup F_2)^C) \oplus \bigoplus\nolimits_{f \in F_1 \cup F_2} \Ker(f^C), \\
            \Image(F_i^C) &= \Image((F_1 \cup F_2)^C) \oplus \bigoplus\nolimits_{f \in F_{3-i} \setminus F_i} \Ker(f^C).
        \end{align*}
        We immediately see that the map $\pi_{\Image(F_1^C)} \circ \pi_{\Image(F_2^C)}$ is the canonical projection from $\VV = \Image((F_1 \cup F_2)^C) \oplus \Ker((F_1 \cup F_2)^C)$ to $\Image((F_1 \cup F_2)^C)$, i.e., $\pi_{\Image((F_1 \cup F_2)^C)}$.
        \item One can verify \hbox{$\pi_{\Image(F_1^C)} \circ \pi_{\Ker(F_2^C)} = \pi_{\Ker((F_2 \setminus F_1)^C)}$} using similar arguments as in (v).

        We now prove that the equation $\pi_{\Image(F_1^C)} \circ \pi_{\Ker(F_2^C)} = \chi[\theta] \circ \pi_{\Ker(F_2^C)}$ holds for some polynomial $\chi \in K[X]$.
        Recall that $F_1^C[\theta]_{\restriction \Ker((F_2 \setminus F_1)^C)}$ is an automorphism of the subspace $\Ker((F_2 \setminus F_1)^C)$ by (v) of Corollary \ref{corollary_polynomials_kernel_facts}, with $\chi'[\theta]_{\restriction \Ker((F_2 \setminus F_1)^C)}$ as inverse for some $\chi' \in K[X]$.
        Using the decomposition
        $$
        \VV = \Image(F_2^C) \oplus \underbrace{\Ker((F_1 \cap F_2)^C) \oplus \Ker((F_2 \setminus F_1)^C)}_{=\Ker(F_2^C)},
        $$
        one can verify \hbox{$(\chi' \cdot F_1^C)[\theta] \circ \pi_{\Ker(F_2^C)} = \pi_{\Image(F_1^C)} \circ \pi_{\Ker(F_2^C)}$} together with the first equation.
        Since $\chi'$ depends only on $F_1^C$ and $F_2^C$, so does $\chi' \cdot F_1^C$.
        The case $F_1 \cap F_2 = \varnothing$ is clear from the first equation.
        \setcounter{enumi}{9}
        \item Fix some element $v \in \VV$.
        By definition, $\eta_2[\theta]^{-1}(v)$ is the unique $u \in \Image(\Fac(\eta_2)^C)$ with $\eta_2[\theta](u) = \pi_{\Image(\Fac(\eta_2)^C)}(v)$, and $\eta_1[\theta]^{-1}(u)$ is the unique $u' \in \Image(\Fac(\eta_1)^C)$ with $\eta_1[\theta](u') = \pi_{\Image(\Fac(\eta_1)^C)}(u)$.
        Unwrapping these definitions and applying (ii) and (v), we obtain
        \begin{align*}
            (\eta_1 \cdot \eta_2)[\theta](u') &= \pi_{\Image(\Fac(\eta_1 \cdot \eta_2)^C)}(v), \quad \text{and also} \\
            u' &= \pi_{\Image(\Fac(\eta_1 \cdot \eta_2)^C)}(\eta_1[\theta]^{-1}(\eta_2[\theta]^{-1}(v))) \in \Image(\Fac(\eta_1 \cdot \eta_2)^C)
        \end{align*}
        by (vii) and (v).
        Therefore, $u' = \eta_1[\theta]^{-1} \circ \eta_2[\theta]^{-1}(v)$ must indeed be $(\eta_1 \cdot \eta_2)[\theta]^{-1}(v)$.
        \setcounter{enumi}{3}
        \item Setting $\chi := \gcd(\rho, \eta)$, we prove the equation as follows:
        $$
        \begin{array}{rcl}
             \rho[\theta] \circ \eta[\theta]^{-1}\hspace{-10pt} &\underset{\text{(vii)}}{=}& \rho[\theta] \circ \eta[\theta]^{-1} \circ \pi_{\Image(\Fac(\eta)^C)} \\
            &\hspace{-10pt}\underset{\text{(i), (x)}}{=}\hspace{-10pt}& (\rho/\chi)[\theta] \circ \chi[\theta] \circ \chi[\theta]^{-1} \circ (\eta/\chi)[\theta]^{-1} \circ \pi_{\Image(\Fac(\eta)^C)} \\
            &\underset{\ref{rem_no_cancel}}{=}& (\rho/\chi)[\theta] \circ \pi_{\Image(\Fac(\chi)^C)} \circ (\eta/\chi)[\theta]^{-1} \circ \pi_{\Image(\Fac(\eta)^C)} \\
             &\underset{\ref{theorem_r_c_def}}{=}& (\rho/\chi)[\theta] \circ (\eta/\chi)[\theta]^{-1} \circ \pi_{\Image(\Fac(\chi)^C)} \circ \pi_{\Image(\Fac(\eta)^C)} \\
              &\underset{\text{(v)}}{=}& (\rho/\chi)[\theta] \circ (\eta/\chi)[\theta]^{-1} \circ \pi_{\Image(\Fac(\eta)^C)}.
        \end{array}
        $$
        For the application of (v) in the last equation, use $\Fac(\chi) \subseteq \Fac(\eta)$. \qedhere
    \end{enumerate}
\end{proof}
\end{lemma}

\noindent Lemma \ref{lemma_communative} above has many consequences.
First of all, we can shrink the set in Theorem \ref{theorem_r_c_def} that generates $R_C$:

\begin{corollary} \label{corollary_better_generators}
    The ring $(R_C, 0, \Id, +, \circ)$ is $\set{+, \circ}$-generated by the set
    \begin{align*}
        \set{\lambda[\theta] : \lambda \in K} &\cup \set{g[\theta] : g \in \Kp{}} \cup \set{f[\theta]^{-1} : f \in \Kp{C<\infty}}.
    \end{align*}
\end{corollary}

\begin{corollary} \label{corollary_has_rho}
    Given $r \in R_C$ and finite $F \subseteq \Kp{0<C<\infty}$, there is a unique polynomial $\rho \in K[X]$ with $\deg(\rho) < \deg(F^C)$ such that
    $$
    r \circ \pi_{\Ker(F^C)} = \rho[\theta] \circ \pi_{\Ker(F^C)}.
    $$
    We also have $\rho_1[\theta] \circ \pi_{\Ker(F^C)} \neq \rho_2[\theta] \circ \pi_{\Ker(F^C)}$ for any $\rho_1 \neq \rho_2$ with $\deg(\rho_1), \deg(\rho_2) < \deg(F^C)$.
\begin{proof}
    This is trivial if $F$ is empty. The existence of such a $\rho \in K[X]$ is immediate from (iii), (vi), and (ix) of Lemma \ref{lemma_communative}.
    It remains to show $\rho_1[\theta] \circ \pi_{\Ker(F^C)} \neq \rho_2[\theta] \circ \pi_{\Ker(F^C)}$ for two polynomials $\rho_1, \rho_2 \in K[X]$ with $\deg(\rho_i) < \deg(F^C)$ for $i = 1, 2$.
    For this, we construct a $C$-image-complete $(\VV, \theta)$ that satisfies this inequality.

    Let $B = (b_k^i : k \in \omega, 0 \leq i < \deg(F^C))$ be a basis of a countable vector space $\VV$.
    Similar to previous constructions, we define $\theta$ on $B$ by setting
    $$
    \theta(b_k^i) := \begin{cases}
        b_k^{i+1} & \text{if $i < \deg(F^C) - 1$,} \\
        -\sum\nolimits_{j=0}^{\deg(F^C)-1} (F^C)_j \cdot b_k^j & \text{if $i = \deg(F^C) - 1$.}
    \end{cases}
    $$
    and extend this to an endomorphism of all of $\VV$.
    It is easy to check that this defines a $C$-image-complete endomorphism that satisfies the inequality $\rho_1[\theta] \circ \pi_{\Ker(F^C)} \neq \rho_2[\theta] \circ \pi_{\Ker(F^C)}$ for $\rho_1$ and $\rho_2$ as stated.
\end{proof}
\end{corollary}

\begin{theorem}[Common Base Theorem]\label{theorem_r_c_element}
    Given $r_1, \dots, r_q \in R_C$ and any finite subset $F_0 \subseteq \Kp{0<C<\infty}$, we can write
    $$
    r_i = \rho_{F,i}[\theta] \circ \eta[\theta]^{-1} \circ \pi_{\Image(F^C)} + \sum\nolimits_{f\in F} \rho_{f,i}[\theta] \circ \pi_{\Ker(f^C)}
    $$
    for all $i \in \set{1, \dots, q}$, where
    \begin{enumerate}[(i)]
        \item $F \subseteq \Kp{0<C<\infty}$ is finite with $F_0 \subseteq F$,
        \item $\eta$ is a monic polynomial with $\Fac(\eta) \subseteq F \cup \Kp{C=0}$,
        \item the $\rho_{F,i}$ satisfy both $\Fac(\rho_{F,i}) \subseteq F \cup \Kp{C=0}\cup\Kp{C=\infty}$ and $\gcd(\rho_{F,1}, \dots, \rho_{F,q}, \eta) = 1$,
        \item the $\rho_{f,i}$ are polynomials with $\deg(\rho_{f,i}) < \deg(f^C)$.
    \end{enumerate}
\begin{proof}
    By the definition of $R_C$ and Lemma \ref{lemma_communative}, we can write each $r_i$ as
    $$
    \sum\nolimits_{k=1}^m \rho_{i,k}[\theta] \circ \eta_{i,k}[\theta]^{-1} \circ \pi_{\Image(F_{i,k}^C)},
    $$
    where each $\rho_{i,k}$ is a polynomial, each $F_{i,k}$ is a finite subset of $\Kp{0<C<\infty}$, and each $\eta_{i,k}$ is a monic polynomial with $\Fac(\eta_{i,k}) \subseteq F_{i,k} \cup \Kp{C=0}$.
    We define $F := F_0 \cup \bigcup_{i=1}^q \bigcup_{k=1}^m F_{i,k}$.
    Using the \projectionRule{} (Lemma \ref{lemma_projection_rule}), we see that
    $$
    \pi_{\Image(F_{i,k}^C)} = \pi_{\Image(F^C)} + \sum\nolimits_{f \in F \setminus F_{i,k}} \pi_{\Ker(f^C)}.
    $$
    Together with Corollary \ref{corollary_has_rho}, this shows that we can write $r_i$ as
    $$
    \sum\nolimits_{k=1}^m \rho_{i,k}[\theta] \circ \eta_{i,k}[\theta]^{-1} \circ \pi_{\Image(F^C)} + \sum\nolimits_{f\in F} \rho_{f,i}[\theta] \circ \pi_{\Ker(f^C)}
    $$
    for some polynomials $\rho_{f,i}$ as described.
    Using Lemma \ref{lemma_communative} and the linearity of all maps involved, one can easily show that the sum
    $
    \sum\nolimits_{k=1}^m \rho_{i,k}[\theta] \circ \eta_{i,k}[\theta]^{-1} \circ \pi_{\Image(F^C)}
    $
    is equal to
    $$
    \bigg(\underbrace{\sum\nolimits_{k=1}^m \rho_{i,k} \cdot \prod\nolimits_{j=1}^q \prod\nolimits_{l=1, (j,l) \neq (i,k)}^m \eta_{j,l}}_{=: \rho_{F,i}}\bigg)[\theta] \circ \bigg(\underbrace{\prod\nolimits_{j=1}^q \prod\nolimits_{l=1}^m \eta_{j,l}}_{=: \eta}\bigg)[\theta]^{-1} \circ \pi_{\Image(F^C)}.
    $$
    One can again use Lemma \ref{lemma_communative} to ensure that $\gcd(\rho_{F,1}, \dots, \rho_{F,q}, \eta) = 1$ for the $\rho_{F,i}$ and $\eta$ defined in the equation above.
    To ensure that $\Fac(\rho_{F,i}) \subseteq F \cup \Kp{C=0} \cup \Kp{C=\infty}$ holds, we define the superset $F' := F \cup \bigcup_{i=1}^q (\Fac(\rho_{F,i}) \cap \Kp{0<C<\infty})$ and use the fact that
    $$
    \pi_{\Image(F^C)} = \pi_{\Image(F'^C)} + \sum\nolimits_{f \in F' \setminus F} \pi_{\Ker(f^C)}
    $$
    together with Corollary \ref{corollary_has_rho} as above to replace $F$ with $F'$.
\end{proof}
\end{theorem}

\begin{corollary} \label{corollary_R_C_as_ring}
    In the following, ``$\simeq$'' means ``isomorphic as rings'':
\begin{enumerate}[(i)]
    \item $R_C \simeq K[X]/(\mipo(C)) \simeq \bigoplus\nolimits_{f\in \Kp{0<C<\infty}} K[X]/(f^C)$, if $C$ is algebraic.
    \item If $C$ is a transcendental kernel configuration with $C(f) \in \set{0, \infty}$ for all $f \in \Kp{}$, then we have
    $
    R_C \simeq (\Kp{C=0})^{-1}K[X],
    $
    where $(\Kp{C=0})^{-1}K[X]$ is the localization of $K[X]$ at the multiplicative closure of $\Kp{C=0}$.
    \item If $F \subseteq \Kp{0<C<\infty}$ is finite, then
    $
    R_C \simeq R_{C'} \oplus \bigoplus\nolimits_{f\in F} K[X]/(f^C),
    $
    where $C'$ is the unique kernel configuration that is algebraic or transcendental if and only if $C$ is algebraic or transcendental, respectively, and satisfies $C'(f) = C(f)$ if $f \not\in F$, and $C'(f) = 0$ otherwise.

    One should note that $C'$ is not an actual kernel configuration if $C$ is algebraic and $F \!=\! \Kp{0<C<\infty}$; in this case we can simply replace $R_{C'}$ with the trivial ring $\set{0}$.
\end{enumerate}
    The combination of (ii) and (iii) yields
    $$
    R_C \simeq (\Kp{C<\infty})^{-1}K[X] \oplus \bigoplus\nolimits_{f\in \Kp{0<C<\infty}} K[X]/(f^C)
    $$
    if the set $\Kp{0<C<\infty}$ is finite and $C$ is transcendental.
\begin{proof}
    For point (i), set $F := \Kp{0<C<\infty}$.
    With this, the \projectionRule{} (Lemma \ref{lemma_projection_rule}) yields $\Id = \sum\nolimits_{f\in F}\pi_{\Ker(f^C)}$.
    First, notice that the function $\sum\nolimits_{f\in F} (\rho_{1, f} \cdot \rho_{2, f})[\theta] \circ \pi_{\Ker(f^C)}$ is just
    \begin{align*}
        \Big(\sum\nolimits_{f\in F} \rho_{1, f}[\theta] \circ \pi_{\Ker(f^C)}\Big) \circ \Big(\sum\nolimits_{f\in F} \rho_{2, f}[\theta] \circ \pi_{\Ker(f^C)}\Big)
    \end{align*}
    as $\pi_{\Ker(f^C)} \circ \pi_{\Ker(g^C)} = 0$ for $f \neq g \in \Kp{0<C<\infty}$.
    This shows that the following map defines a ring homomorphism:
    $$
    \begin{array}{rcl}
         R_C & \to & \bigoplus\nolimits_{f\in F} K[X]/(f^C) \\
         \sum\nolimits_{f\in F} \rho_f[\theta] \circ \pi_{\Ker(f^C)} & \mapsto & ([\rho_f] : f \in F).
    \end{array}
    $$
    Using Corollary \ref{corollary_has_rho}, we can easily verify that this map is bijective.
    
    Point (ii) is immediate by the \commonBaseTheorem{} (Theorem \ref{theorem_r_c_element}).
    Point (iii) in the algebraic case follows from (i).
    For (iii) in the transcendental case, we use the \commonBaseTheorem{} to write every element $r \in R_C$ as
    $$
    \rho[\theta] \circ \eta[\theta]^{-1} \circ \pi_{\Image((F' \cup F)^C)} + \sum\nolimits_{f\in F'} \rho_f[\theta] \circ \pi_{\Ker(f^C)} + \sum\nolimits_{f\in F} \rho_f[\theta] \circ \pi_{\Ker(f^C)}
    $$
    for some finite $F' \subseteq \Kp{0<C<\infty}$ with $F' \cap F = \varnothing$, which is equivalent to $F' \subseteq \Kp{0<C'<\infty}$.
    We now define a map $R_C \to R_{C'} \oplus \bigoplus\nolimits_{f\in F} K[X]/(f^C)$ by mapping such an element $r$ to
    $$
    (\rho[\theta] \circ \eta[\theta]^{-1} \circ \pi_{\Image(F'^C)} + \sum\nolimits_{f\in F'} \rho_f[\theta] \circ \pi_{\Ker(f^C)})^\frown ([\rho_f] : f \in F).
    $$
    Using the following claim and Lemma \ref{lemma_communative} (notice the respective comments about the $\chi$'s in (vi) and (ix)) we can easily verify that this map is well defined.
\begin{subclaim} \label{lemma_when_equivalent}
    Let $C$ be transcendental and let $r_1, r_2 \in R_C$ be written as in the \commonBaseTheorem{}, i.e., for $i = 1, 2$ we have
    $$
    r_i = \rho_i[\theta] \circ \eta_i[\theta]^{-1} \circ \pi_{\Image(F_i^C)} + \sum\nolimits_{f \in F_i} \rho_{f, i}[\theta] \circ \pi_{\Ker(f^C)}
    $$
    with $F_i \subseteq \Kp{0<C<\infty}$, $\eta_i$ monic with $\Fac(\eta_i) \subseteq F_i \cup \Kp{C=0}$, and $\deg(\rho_{f, i}) < \deg(f^C)$ for all $f \in F_i$.
    We have $r_1 = r_2$ if and only if all of the following conditions hold:
    \begin{enumerate}[(i)]
        \item $\rho_1/\eta_1 = \rho_2/\eta_2$ in $K(X)$,
        \item $\rho_2[\theta] \circ \eta_2[\theta]^{-1} \circ \pi_{\Ker(f^C)} = \rho_{f, 1}[\theta] \circ \pi_{\Ker(f^C)}$ for all $f \in F_1 \setminus F_2$,
        \item $\rho_1[\theta] \circ \eta_1[\theta]^{-1} \circ \pi_{\Ker(f^C)} = \rho_{f, 2}[\theta] \circ \pi_{\Ker(f^C)}$ for all $f \in F_2 \setminus F_1$,
        \item $\rho_{f, 1} = \rho_{f, 2}$ in $K[X]$ for all $f \in F_1 \cap F_2$.
    \end{enumerate}
\begin{innerproof}
    By the \projectionRule{}, we have $1[\theta] = \pi_{\Image((F_1 \cup F_2)^C)} + \sum\nolimits_{f\in F_1 \cup F_2} \pi_{\Ker(f^C)}$.
    Applying this decomposition to both $r_i$, we see that $r_1 = r_2$ holds if and only if all of the following conditions hold:
    \begin{enumerate}[(i)]
        \item $\rho_1[\theta] \circ \eta_1[\theta]^{-1} \circ \pi_{\Image((F_1\cup F_2)^C)} = \rho_2[\theta] \circ \eta_2[\theta]^{-1} \circ \pi_{\Image((F_1\cup F_2)^C)}$,
        \item $\rho_2[\theta] \circ \eta_2[\theta]^{-1} \circ \pi_{\Ker(f^C)} = \rho_{f, 1}[\theta] \circ \pi_{\Ker(f^C)}$ for all $f \in F_1 \setminus F_2$,
        \item $\rho_1[\theta] \circ \eta_1[\theta]^{-1} \circ \pi_{\Ker(f^C)} = \rho_{f, 2}[\theta] \circ \pi_{\Ker(f^C)}$ for all $f \in F_2 \setminus F_1$,
        \item $\rho_{f, 1}[\theta] \circ \pi_{\Ker(f^C)} = \rho_{f, 2}[\theta] \circ \pi_{\Ker(f^C)}$ for all $f \in F_1 \cap F_2$.
    \end{enumerate}
    That (iv) is equivalent to $\rho_{f, 1} = \rho_{f, 2}$ in $K[X]$ for all $f \in F_1 \cap F_2$ is clear by Corollary \ref{corollary_has_rho}.
    Next, we set $F := F_1 \cup F_2$ and show that
    $$
    \rho_1[\theta] \circ \eta_1[\theta]^{-1} \circ \pi_{\Image(F^C)} = \rho_2[\theta] \circ \eta_2[\theta]^{-1} \circ \pi_{\Image(F^C)} \text{ in } R_C
    $$
    if and only if $\rho_1/\eta_1 = \rho_2/\eta_2$ in $K(X)$.
    One can easily check that the direction from left to right holds for every $C_0$-image-complete endomorphism, and since every $C_0$-image-complete endomorphism is also $C$-image-complete, it must generally hold in $R_C$.
    The other direction is just (iv) of Lemma \ref{lemma_communative}.
\end{innerproof}
\end{subclaim}
    
    \noindent Once the well-definedness is shown, one can prove that this map is bijective.
    With a similar argument as in (i), we see that it is indeed a ring homomorphism.
\end{proof}
\end{corollary}

\begin{remark}\label{fact_when_field}
    $R_C$ is a field if and only if either $C = C_0$ or $C$ is algebraic with $\mipo(C)$ irreducible.
\end{remark}

\begin{remark} \label{rem_dim_r_C}
    As $K \subseteq R_C$, the ring $R_C$ is a $K$-vector space, or even a $K$-algebra.
    Furthermore, the following holds:
    \begin{enumerate}[(i)]
        \item If $C$ is algebraic, then $\dim_K(R_C) = \deg(\mipo(C))$.
        \item If $C$ is transcendental, then $\dim_K(R_C) = \max\set{\omega, |\Kp{C<\infty}|}$.
    \end{enumerate}
\begin{proof}
    In the transcendental case, one can use partial fraction decomposition to verify that
    \begin{align*}
        \Big\{X^i[\theta] : i \in \NN \Big\} &\cup \Big\{X^i[\theta] \!\circ\! f^n[\theta]^{-1} : f \in \Kp{C<\infty}, 0 \leq i < \deg(f), n \geq 1 \Big\} \\
        &\cup\Big\{ X^i[\theta] \circ \pi_{\Ker(f^C)} : f \in \Kp{0<C<\infty}, 0 \leq i < \deg(f^C) \Big\}
    \end{align*}
    is a $K$-basis of $R_C$.
    From there, one can easily calculate the cardinality.
    If $C$ is an algebraic kernel configuration, one has to omit the first line and recall $\mipo(C) = \prod_{f \in \Kp{0<C<\infty}} f^C$.
\end{proof}
\end{remark}

\noindent One might be tempted to think that two kernel configurations $C_1, C_2 \in \Cc$ are equal if and only if $R_{C_1}$ and $R_{C_2}$ are isomorphic as rings.
This is, however, not the case:

\begin{remark} \label{rem_r_C_fail}
    We might have $(R_C, 0, 1, +, \cdot) \simeq (R_{C'}, 0, 1, +, \cdot)$ as rings for two different kernel configurations $C \neq C' \in \Cc$.
    This can, for example, happen if $C'$ is obtained from $C$ by shifting by some $\lambda \in K$, i.e., if we have $\deg(C') = \deg(C)$, see Definition \ref{def_kernel_conf}, and $C'(f) = C(f \circ (X - \lambda))$ for all $f \in \Kp{}$.
\end{remark}
 
\subsection{The characterization}
\label{sec_char}

\noindent Recall our first characterization of existentially closed models of $T^C_\theta$ (Corollary \ref{corollary_exist_closed_0}):
A model $(\mm, \theta) \models T^C_\theta$ is existentially closed if and only if, for any $L(M)$-formula $\varphi(\uxvec)$, we have 
$$
(\mm, \theta) \models \exists \ux \in \VV : \varphi_\theta(\ux)
$$ 
whenever there is $(\mm', \theta') \models T^C_\theta$ with both $(\mm, \theta) \subseteq (\mm', \theta')$ and \hbox{$(\mm', \theta') \models \exists \ux \in \VV : \varphi_\theta(\ux)$}.
(Also recall our \placeholderNotation{}: Given a tuple $\ux = (x_1, \dots, x_n)$, we define the tuple $\uxvec = (x^i_k : 1 \leq k \leq n, i \in \omega)$, where each $x^i_k$ is a placeholder for $\theta^i(x_k)$. Given an $L(M)$-formula $\varphi(\uxvec)$, we define $\varphi_\theta(\ux)$ to be the $L_\theta(M)$-formula with all placeholders replaced accordingly.)

As stated (a bit differently) at the end of Section \ref{sec_first_char}, we will now refine this characterization so that one only needs to consider formulas $\varphi(\uxvec)$ for which $\varphi_\theta(\ux)$ is of the form
$$
\psi_\theta(\ux) \wedge E(\ux)
$$
for some $L(M)$-formula $\psi(\uxvec)$ that implies no finite disjunction of non-trivial linear dependencies in $\uxvec$ over $\VV$, and some nice $\LKThe$-formula $E(\ux)$ that dictates which linear dependencies over $\VV$ the sequence $(\theta^i(\ux) : i \in \omega)$ satisfies.
We start by defining the class of these $\LKThe$-formulas:

\begin{definition}
    \label{def_param_c_sequence_system} A \textbf{parametrized $\mathbf{C}$-sequence-system} is an $\LKThe$-formula of the form
    $$
    S(\ux; \uy) = \bigwedge\nolimits_{k=1}^n f_k^{q_k}[\theta](x_{\ldd, k}) = y_k
    $$
    with $\ux := \ux_\li\ux_\ld := (x_{\lii, k} : 1 \leq k \leq m)(x_{\ldd, k} : 1 \leq k \leq n)$ and $\uy = (y_1, \dots, y_n)$ that satisfies the following conditions:
    \begin{enumerate}[(i)]
        \item If $C$ is algebraic, then $m = 0$.
        \item $f_k \in \Kp{0<C}$ and $q_k \in \set{q \in \NN : 0 < q \leq C(f_k)}$ hold for all $k \in \set{1, \dots, n}$.
    \end{enumerate}
    
\end{definition}
\noindent We will always denote parametrized $C$-sequence-systems by the letter $S$.
Given such a parametrized $C$-sequence-system $S(\ux; \uy)$, we assume that everything is as above, that is, $m$, $n$, and the $f_k$'s and $q_k$'s are defined implicitly, and we set $\ux = \ux_\li\ux_\ld$ as above.
When we partition $\ux = \ux_\li\ux_\ld$ as above, we think of:
    \begin{enumerate}[(i)]
        \item $\ux_\li$ as the linearly independent part of $\ux$, since $S(\ux; \uy)$ does not imply any linear dependencies for the sequence $(\theta^i(x_{\lii, k}) : 1 \leq k \leq m, i \in \omega)$;
        \item $\ux_\ld$ as the linearly dependent part of $\ux$, since the formula $S(\ux; \uy)$ implies that the sequence $(\theta^i(x_{\ldd, k}) : i \in \omega)$ is linearly dependent over $\spanA{y_k}{K}$ for each $k \in \set{1, \dots, n}$.
    \end{enumerate}
The names $\ux_\li$ and $\ux_\ld$ for these tuples are abbreviations for linearly independent and linearly dependent.

\begin{definition} \label{def_formual_bounded}
    Let $S(\ux; \uy)= \bigwedge\nolimits_{k=1}^n f_k^{q_k}[\theta](x_{\ldd, k}) = y_k$ be a parametrized $C$-sequence-system as in Definition \ref{def_param_c_sequence_system}.
    If $\psi(\uxvec; \uw)$ is a formula, then we say $\psi(\uxvec; \uw)$ is \textbf{bounded} by $S$ if one of the following equivalent conditions holds:
        \begin{enumerate}[(i)]
            \item For all $k \in \set{1, \dots, n}$, the variable $x^i_{\ldd, k}$ does not appear in $\psi(\uxvec; \uw)$ for \hbox{$i \geq \deg(f_k^{q_k})$}.
            \item For all $k \in \set{1, \dots, n}$, the term $\theta^i(x_{\ldd, k})$ does not appear in $\psi_\theta(\ux; \uw)$ for \hbox{$i \geq \deg(f_k^{q_k})$}.
        \end{enumerate}
    We also say that $\psi_\theta(\ux; \uw)$ is \textbf{bounded} by $S$, if $\psi(\uxvec; \uw)$ is bounded by $S$.
\end{definition}

\noindent In practice, a formula $\psi(\uxvec)$ being bounded by $S$ means that no subterm of the form $\theta^i(x_{\ldd, k})$ that appears in $\psi_\theta(\ux)$ can be replaced by a \remainderRule{} applied with the equation $f_k^{q_k}[\theta](x_{\ldd, k}) = y_k$ from $S(\ux; \uy)$ (as the remainder of $X^i$ divided by $f_k^{q_k}$ is still $X^i$).

\begin{definition}
    \label{def_c_sequence_system} \label{def_compatible} Let $S(\ux; \uy) = \bigwedge_{k=1}^n f_k^{q_k}[\theta](x_{\ldd, k}) = y_k$ be a parametrized $C$-sequence-system as in Definition \ref{def_param_c_sequence_system}, and let $(\VV, \theta) \models \TKvsThe$ be given.
\begin{enumerate}[(i)]
    \item We say that a tuple $\uu = (u_1, \dots, u_n) \in \VV$ is \textbf{compatible} with $S$ if $u_k \in \Ker(f_k^{C-q_k})$ for every $k \in \set{1, \dots, n}$ with $f_k \in \Kp{0<C<\infty}$.
    \item A \textbf{$\mathbf{C}$-sequence-system} over $(\VV, \theta)$ is an $\LKThe(\VV)$-formula of the form
    \hbox{$
    S(\ux) = S'(\ux; \uu)
    $}
    where $S'$ is a parametrized $C$-sequence-system and $\uu \in \VV$ is compatible with $S'$.
\end{enumerate}
\end{definition}

\noindent We will also denote $C$-sequence-systems over some $(\VV, \theta) \models \TKvsTheC$ by the letter $S$.
Notice that a $C$-sequence-system over $(\VV, \theta)$ is also a $C$-sequence-system over any extension $(\VV', \theta')$ that is also a model of $\TKvsThe$. We also say that a formula is bounded by a $C$-sequence-system if it is bounded by the underlying parametrized $C$-sequence-system.
We can now state our characterization of existentially closed models of $T^C_\theta$:

\begin{theorem} \label{theorem_big_characterization}
    A model $(\mm, \theta) \models T^C_\theta$ is existentially closed if and only if it is $C$-image-complete and
    $$
        (\mm, \theta) \models \exists \ux \in \VV : \psi_\theta(\ux) \wedge S(\ux)
    $$
    holds for every $C$-sequence-system $S(\ux)$ over $(\VV, \theta)$ and every $L(M)$-formula $\psi(\uxvec)$ that is bounded by $S$ and does not imply any finite disjunction of non-trivial linear dependencies in $\uxvec$ over $\VV$.
\end{theorem}

\noindent We would like to highlight the characterizations we obtain for the most interesting kernel configurations.
In two cases, namely those in which $R_C$ is a field, we obtain a characterization that does not even need $C$-sequence-systems.
The first case in which $R_C$ is a field is when $C$ is algebraic and $\mipo(C)$ is irreducible, i.e., when we have $T^C_\theta = T_\theta \cup \set{\forall x \in \VV : f[\theta](x) = 0}$ for some $f \in \Kp{}$:

\begin{theorem} \label{theorem_Calgirr_char}
    Let $C$ be an algebraic kernel configuration such that $f := \mipo(C)$ is irreducible, and set $d := \deg(f)$.
    A model $(\mm, \theta) \models T^{C}_\theta$ is existentially closed if and only if we have
    $$
    (\mm, \theta) \models \exists \ux \in \VV: \psi(\theta^0(\ux) \dots \theta^{d-1}(\ux))
    $$
    for all $L(M)$-formulas $\psi(\ux^0 \dots \ux^{d-1})$ (with $|\ux^i| = |\ux^0|$ for all $i$) that do not imply any finite disjunction of non-trivial linear dependencies in $\ux^0 \dots \ux^{d-1}$ over $\VV$.
\begin{proof}
    In this case, a $C$-sequence-system is just of the form $\bigwedge\nolimits_{k=1}^n f[\theta](x_k) = 0$.
    A formula $\psi(\uxvec)$ is bounded by such a $C$-sequence-system if and only if $x^i_k$ does not appear for any $i \geq d$, so we can work with a formula of the form $\psi(\ux^0, \dots, \ux^{d-1})$.
    Since we have $T^C_\theta \models \forall x \in \VV : f[\theta](x) = 0$, we can omit the sequence-system $S$.
    Also, recall that $C$-image-completeness holds trivially for algebraic kernel configurations.
\end{proof}
\end{theorem}

\noindent The second case in which $R_C$ is a field is when $C = C_0$ for the unique transcendental kernel configuration $C_0$ given by $C_0(f) = 0$ for all monic irreducible polynomials $f$.
As we have seen previously, this means that $T^{C_0}_\theta := T_\theta \cup \set{\text{``$\rho[\theta]$ is injective''} : \rho \in K[X] \setminus \set{0}}$.

\begin{theorem} \label{theorem_C0_char}
    A model $(\mm, \theta) \models T^{C_0}_\theta$ is existentially closed if and only if $\rho[\theta]$ is bijective for every $\rho \in K[X] \setminus \set{0}$, and we have
    $$
    (\mm, \theta) \models \exists \ux \in \VV: \psi_\theta(\ux)
    $$
    for all $L(M)$-formulas $\psi(\uxvec)$ that do not imply any finite disjunction of non-trivial linear dependencies in $\uxvec$ over $\VV$.
\begin{proof}
    By definition, all $C_0$-sequence-systems are just $\top$, and $C_0$-image-completeness is equivalent to the condition that, for every $\rho \neq 0$, the map $\rho[\theta]$ is bijective.
\end{proof}
\end{theorem}

\noindent The last case we highlight is the case $C = C_\infty$, where $C_\infty$ is the unique transcendental kernel configuration given by $C_\infty(f) = \infty$ for all $f \in \Kp{}$.
Notice that $T^{C_\infty}_\theta = T_\theta$.
In this case, we can at least work with rather simple $C$-sequence-systems:

\begin{theorem}
    A model $(\mm, \theta) \models T_\theta$ is existentially closed if and only if
    $$
    (\mm, \theta) \models \exists \ux \in \VV: \psi_\theta(\ux) \wedge \bigwedge\nolimits_{k=1}^n f_k^{q_k}[\theta](x_{m+k}) = u_k,
    $$
    where $\ux = (x_1, \dots, x_{m+n})$, 
    holds for all $m, n \geq 0$, $f_1, \dots, f_n \in \Kp{}$, $q_1, \dots, q_n \in \NN_{>0}$, $u_1, \dots, u_n \in \VV$, and $L(M)$-formulas $\psi(\uxvec)$ that do not imply any finite disjunction of non-trivial linear dependencies in $\uxvec$ over $\VV$ and do not contain $x_{m + k}^i$ for $i \geq \deg(f_k^{q_k})$.
\begin{proof}
    One can easily check that $C_\infty$-image-completeness is a trivial condition, and the rest follows from the form of $C_\infty$-sequence-systems.
\end{proof}
\end{theorem}

\subsection{First-order axiomatization}

\noindent We now have a characterization of existentially closed models of $T^C_\theta$.
In order to obtain a model companion of $T^C_\theta$, we need to axiomatize this characterization in a first-order way. The following is clear by Definition \ref{def_c_sequence_system}:
\begin{lemma}
    Given a parametrized $C$-sequence-system $S(\ux; \uy)$, there is a $\LKThe$-formula $\delta_{S}(\uy)$ such that for any $(\mm, \theta) \models T_\theta$ and $\uu \in \VV$, we have
    $$
    (\mm, \theta) \models \delta_{S}(\uu) \quad \Leftrightarrow \quad \text{``$S(\ux; \uu)$ is a $C$-sequence-system''}.
    $$
\end{lemma}
\noindent Looking at Theorem \ref{theorem_big_characterization}, we take the theory $T_\theta \cup \set{\text{``$\theta$ is $C$-image-complete''}}$ and add the sentence
\begin{align*}
    &\forall \uw:\forall \uy \in \VV : \left(\parbox{8.0cm}{``$\psi(\uxvec; \uw)$ implies no finite disjunction of non-\\ ${}$\hspace{5pt}trivial linear dependencies in $\uxvec$ over $\VV$''} \wedge \delta_S(\uy) \right)\\
    &\hspace{230pt}\rightarrow \; \exists \ux \in \VV : \psi_\theta(\ux; \uw) \wedge S(\ux; \uy)
\end{align*}
for every parametrized $C$-sequence-system $S(\ux; \uy)$ and every $L$-formula $\psi(\uxvec; \uw)$ that is bounded by $S$.
It remains only to verify that the statement in quotes can be expressed as a first-order formula in $\uw$.
However, this is not always the case:

\begin{example} \label{non_example}
    Consider $\RCF$, the theory of real closed fields, with $\VV$ being the $\QQ$-vector space induced by the addition.
    The $L_{\operatorname{ring}}(r)$-formula
    $
    \psi(\xvec; r) := x^1 = r \cdot x^0
    $
    implies a finite disjunction of non-trivial linear dependencies in $\xvec$ over $\VV = \rr$ if and only if $r \in \QQ$.
    Therefore,
    $$
    \text{``$\psi(\xvec; w)$ implies no finite disjunction of non-trivial linear dependencies in $\xvec$ over $\VV$''}
    $$
    cannot be expressed by a first-order formula $\sigma_\psi(w)$, since $\QQ$ is clearly not definable in $\RCF$.
    
    Even worse, there is no way to fix this:
    Suppose we could find an expansion \hbox{$T \supseteq \RCF$} in which such a formula $\sigma_\psi(w)$ exists.
    By the above, we would obtain $\neg\sigma_\psi(\rr) = \QQ$ in every model $\rr \models T$.
    This is impossible, since by compactness there must be a model $\rr \models T$ with $|\neg\sigma_\psi(\rr)| > |\QQ|$.
\end{example}

\noindent Our solution to this problem is quite simple: We add the assumption that this is possible.

\begin{definition}[see Definition 1.11 in \cite{dEl21b}] \label{def_hfour}
    We say that $T$ (with the specific choice of $\VV$) satisfies $(\operatorname{H4})$ if, for every $L$-formula $\psi(\ux; \uw)$, there is an $L$-formula $\sigma_\psi(\uw)$ such that, for all $\mm \models T$ and $\ud \in M$, we have $\mm \models \sigma_\psi(\ud)$ if and only if one of the following two equivalent conditions holds:
    \begin{enumerate}[(i)]
        \item The formula $\psi(\ux; \ud)$ implies no finite disjunction of non-trivial linear dependencies in $\ux$ over $\VV$.
        \item There are an elementary extension $\mm' \succ \mm$ and a tuple $\uv' \in \VV'$ linearly independent over $\VV$ such that $\mm' \models \psi(\uv'; \ud)$.
    \end{enumerate}
\end{definition}

\noindent Thus, $\sigma_\psi(\uw)$ is exactly the formula that describes
$$
\text{``$\psi(\ux; \uw)$ implies no finite disjunction of non-trivial linear dependencies in $\ux$ over $\VV$''}.
$$
This gives us the following first-order axiomatization of existentially closed models of $T^C_\theta$ under the assumption that \Hfour{} holds:

\begin{theorem} \label{theorem_first_oder}
If $T$ satisfies \Hfour{}, then $T_\theta^C$ has a model companion $T\theta^C$, i.e., a first-order axiomatization of the class of existentially closed models.
It is axiomatized by the theory \hbox{$T_\theta \cup \set{\text{``$\theta$ is $C$-image-complete''}}$} together with the sentence
$$
\forall \uw:\forall \uy \in \VV : (\sigma_\psi(\uw) \wedge \delta_S(\uy)) \rightarrow \exists \ux \in \VV : \psi_\theta(\ux; \uw) \wedge S(\ux; \uy)
$$
for every parametrized $C$-sequence-system $S(\ux; \uy)$ and every $L$-formula $\psi(\uxvec; \uw)$ that is boun\-ded by $S$.
\end{theorem}

\noindent We do not know whether \Hfour{} is actually necessary for the existence of a model companion.

\begin{example} \label{examples_hfour}
    \Hfour{} holds in the following settings:
    \begin{enumerate}[(i)]
        \item The theory $\TKvs$ with the vector space $(\VV, +, 0, (\lambda \cdot)_{\lambda \in K})$ being the entire structure satisfies \Hfour{}.
        This follows easily from quantifier elimination.
        \item Any complete and model-complete o-minimal theory $T$ extending the theory of divisible ordered abelian groups, with $(\VV, +, 0, (q \cdot)_{q \in \QQ})$ being a subinterval of the line (but not necessarily a subgroup) and continuous operations, satisfies \Hfour{} if and only if there is no infinite definable family of germs of $(\VV, +, 0, (q \cdot)_{q \in \QQ})$-endomorphisms at $0_\VV$; combine Theorem 2.4 and Lemma 2.9 in \cite{Blo23}.
        \item Any complete and model-complete o-minimal expansion of $\RCF{}$ with $(\VV, +, 0, (q \cdot)_{q \in \QQ})$ given by $(R_{>0}, \cdot, 1, (x \mapsto x^q)_{q\in \QQ})$ satisfies \Hfour{} if and only if no partial exponential function is definable.
        This is a special case of (ii); see the proof of Theorem A in \cite{Blo23}.
        \item Let $\FF_q$ be a finite field with $q = p^r$.
        If $T$ expands the theory of an $\FF_q$-vector space and $\VV$ is that vector space, then $T$ satisfies \Hfour{} if and only if it eliminates $\exists^\infty$; see the proof of Theorem 5.2 in \cite{dEl21b}.

        For expansions of the theories $\ACF{}_p$, $\operatorname{SCF}_{p, e}$ ($e$ either finite or infinite), $\operatorname{Psf}_p$, $\operatorname{ACFA}_p$, and $\operatorname{DCF}_p$, this implies that whenever $\FF_q$ is contained in every model as constants, \Hfour{} holds with the $\FF_q$-vector space given by addition.
        For more details, see Example 5.10 in \cite{dEl21b}.
    \end{enumerate}
    Additionally \Hfour{} (or rather the respective generalization) also holds for $\ACF{}$ with the multiplicative group modulo torsion being a $\QQ$-vector space (see Remark \ref{rk:ACFcase}). To see this use Fact 3.9 with $\ux'$ empty in \cite{dEl25}. Generalizing the torsion case, one can also show \Hfour{} for abelian varieties in the setting presented in \cite{dEl21a} (for \Hfour{} see Proposition 2.7 there).
\end{example}

\noindent In future papers, we will see other examples where \Hfour{} holds. We will also see that $T^C_\theta$ does not have a model companion if $C$ is non-trivial and $T$ does not eliminate $\exists^\infty x \in \VV$.

\section{Proof of the Characterization}

\subsection{Transformability}
\label{sec_trafo}

\noindent The goal of this section is to prove that any conjunction of a parametrized $C$-sequence-system, as in Definition \ref{def_param_c_sequence_system}, and some $\LKThe$-equations can, more or less, be ``reduced'' to another parametrized $C$-sequence-system.
We also show that more equations generally result in fewer ``degrees of freedom'' for the solutions.

The results of this section are essential for the proof of the direction from right-to-left of Theorem \ref{theorem_big_characterization} (Section \ref{sec_right_to_left}).
However, we state all results in broader generality to establish some fundamentals for future papers.

\begin{definition}
    We define $\LRC$ as the language of (left) $R_C$-modules (with $R_C$ as defined in Section \ref{sec_rc_module}), i.e., $\LRC = (0, +, (r)_{r \in R_C})$, where each $r \in R_C$ is treated as a unary function symbol.
\end{definition}

\noindent We will write $r(x)$ instead of $r \cdot x$, since $r$ will usually be a function such as $\pi_{\Image(F^C)}$.
Given a model $(\VV, \theta) \models \TKvsThe \cup \set{\text{``$\theta$ is $C$-image-complete''}}$, we define an $\LRC$-structure on $\VV$ using the definable functions of which $R_C$ consists (see Section \ref{sec_rc_module}).
This structure is obviously interdefinable with $(\VV, \theta)$. We now make our ``reduction'' notion from the beginning of this section precise:

\begin{definition} \label{def_transfo}
    Let $E(\ux; \uy)$ and $E'(\ux'; \uy)$ be two conjunctions of $\LRC$-equations.
    We say that $E(\ux; \uy)$ is \textbf{transformable} into $E'(\ux'; \uy)$ if both
    \begin{enumerate}[(i)]
        \item $\TKvsThe \cup \set{\text{``$\theta$ is $C$-image-complete''}} \models \forall \ux\uy: E(\ux; \uy) \rightarrow E'(\underline{\nu}(\ux); \uy) \wedge \ux = \utau(\underline{\nu}(\ux); \uy)$,
        \item $\TKvsThe \cup \set{\text{``$\theta$ is $C$-image-complete''}} \models \forall \ux'\uy: E'(\ux'; \uy) \rightarrow E(\utau(\ux'; \uy); \uy)$
    \end{enumerate}
    hold for some tuple $\underline{\nu}(\ux)$ of $\LRC$-terms and some tuple $\utau(\ux'; \uy)$ with each entry $\tau_k(\ux'; \uy)$ being the sum of an $\LKThe$-term in $\ux'$ and an $\LRC$-term in $\uy$.
    We say that the tuples $\underline{\nu}(\ux)$ and $\utau(\ux'; \uy)$ \textbf{witness} the transformability.
\end{definition}

\noindent When $\utau(\ux'; \uy)$ does not depend on $\uy$, we may simply write $\utau(\ux')$.

\begin{observation}  \label{obser_trafo_nice}
    Transformability satisfies the following properties:
    \begin{enumerate}[(i)]
        \item It is transitive, i.e., if $E(\ux; \uy)$ is transformable into $E'(\ux'; \uy)$ witnessed by $\underline{\nu}(\ux)$ and $\utau(\ux'; \uy)$, and $E'(\ux'; \uy)$ is transformable into $E''(\ux''; \uy)$ witnessed by $\underline{\nu}{}'(\ux')$ and $\utau'(\ux''; \uy)$, then $E(\ux; \uy)$ is transformable into $E''(\ux''; \uy)$ witnessed by $\underline{\nu}{}'(\underline{\nu}(\ux))$ and $\utau(\utau'(\ux''; \uy); \uy)$.
        \item It is \textbf{modular}, i.e., if $E_1(\ux_1; \uy)$ is transformable into $E'_1(\ux'_1; \uy)$, and $E_2(\ux_2; \uy)$ is transformable into $E'_2(\ux'_2; \uy)$, then $E_1(\ux_1; \uy) \wedge E_2(\ux_2; \uy)$ is transformable into the conjunction $E'_1(\ux'_1; \uy) \wedge E'_2(\ux'_2; \uy)$.
        (Notice that we assume here that $\ux_1$ and $\ux_2$ are disjoint tuples of variables, and that the same holds for $\ux'_1$ and $\ux'_2$.)
        \item If $E(\ux; \uy)$ is equivalent to $E'(\ux; \uy)$, then $E(\ux; \uy)$ is also transformable into $E'(\ux; \uy)$, witnessed by $\underline{\nu}(\ux) = \ux$ and $\utau{}(\ux; \uy) := \ux$.
    \end{enumerate}
\end{observation}

\noindent One needs to be a bit careful with (iii) of Observation \ref{obser_trafo_nice}:

\begin{remark}
    Transformability depends on which variables are treated as parameters.
    For example, if $E(\ux; \uy)$ is transformable into $E'(\ux'; \uy)$, then \hbox{$E_*(x_*\ux; \uy) := E(\ux; \uy)$} (which is logically equivalent to $E(\ux; \uy)$) is not transformable into $E'(\ux'; \uy)$ (but it is transformable into $E'_*(x_*\ux'; \uy) := E'(\ux'; \uy)$).
    Also, transformability is not symmetric in general, as we require each entry of $\utau(\ux';\uy)$ to be the sum of an $\LKThe$-term in $\ux'$ and an $\LRC$-term in $\uy$, whereas $\underline{\nu}(\ux)$ can be any tuple of $\LRC$-terms in $\ux$.
\end{remark}

\begin{example} \label{exampl_trafo}
    Given distinct $h_1, \dots, h_m \in \Kp{}$, $r_1, \dots, r_m > 0$, and some $\LRC$-term $\mu(\uy)$, the formula $\big(\prod\nolimits_{k=1}^m h_k^{r_k} \big)[\theta](x) = \mu(\uy)$ is transformable into $\bigwedge_{k=1}^m h_k^{r_k}[\theta](x_k) = \mu(\uy)$.
\begin{proof}
    Assume $m = 2$ and set $\rho = h_1^{r_1}$ and $\eta = h_2^{r_2}$.
    As $\gcd(\rho, \eta) = 1$, Fact \ref{fact_euclid_domain} yields $\chi_1, \chi_2 \in K[X]$ with $1 = \chi_1 \cdot \eta + \chi_2 \cdot \rho$.
    Now, with the tuples $\underline{\nu}(x) := (\eta[\theta](x), \rho[\theta](x))$ and $\tau(x_1, x_2) := \chi_1[\theta](x_1) + \chi_2[\theta](x_2)$, the two implications from Definition \ref{def_transfo} hold:
    \begin{enumerate}[(i)]
        \item $(\rho \cdot \eta)[\theta](x) = \mu(\uy) \rightarrow (\rho[\theta](\nu_1(x)) = \mu(\uy) \wedge \eta[\theta](\nu_2(x)) = \mu(\uy)) \wedge x = \tau(\underline{\nu}(x))$ (clear by definition of the polynomials $\chi_1, \chi_2$).
        \item $\rho[\theta](x_1) = \mu(\uy) \wedge \eta[\theta](x_2) = \mu(\uy) \rightarrow (\rho \cdot \eta)[\theta](\tau(x_1, x_2)) = \mu(\uy)$.
        Assume $v_1, v_2$ and $\uu$ are given in a model of $\TKvsThe \cup \set{\text{``$\theta$ is $C$-image-complete''}}$ with \hbox{$\rho[\theta](v_1) = \mu(\uu) = \eta[\theta](v_2)$}.
        Then $(\rho \cdot \eta)[\theta](\tau(v_1, v_2)) = (\chi_1 \cdot \eta + \chi_2 \cdot \rho)[\theta](\mu(\uu)) = \mu(\uu)$.
    \end{enumerate}
    The general case follows by similar arguments, since Fact \ref{fact_euclid_domain} works for arbitrarily many pairwise relatively prime polynomials, or by induction using the above and Observation \ref{obser_trafo_nice}.
    The case $m=1$ is trivial, and the case $m=0$ follows with $\nu(x)$ being the empty tuple and $\tau(\uy) = \mu(\uy)$.
\end{proof}
\end{example}

\noindent Transformability can be thought of as a stronger form of ``equivalence up to existence'', as $E(\ux; \uy)$ being transformable into $E'(\ux'; \uy)$ implies 
$$
\exists \ux \in \VV : \phi(\ux; \uy) \wedge E(\ux; \uy) \equiv \exists \ux' \in \VV : \phi(\utau(\ux'; \uy); \uy) \wedge E'(\ux'; \uy)
$$ 
for any formula $\phi(\ux; \uy)$.  
This allows us to apply the results of this section to $L_\theta$-formulas.  
We only really care about the case where $E'(\ux'; \uy)$ is (essentially) a parametrized $C$-sequence-system, and $\phi(\ux; \uz)$ is of the form $\psi_\theta(\ux; \uz)$ for some $L$-formula $\psi(\uxvec; \uz)$ (recall our \placeholderNotation{}).  
In this case, we can furthermore replace $\psi_\theta(\utau(\ux'; \uy); \uz)$ with a formula that is bounded by the parametrized $C$-sequence-system:

\begin{lemma} \label{lemma_main_trafo}
    Suppose a conjunction of $\LRC$-equations $E(\ux; \uy)$ is transformable into a formula of the form
    $
    \varphi(\uy) \wedge S(\ux'; \umu{}(\uy)),
    $
    where $S(\ux'; \tiluy)$ is a parametrized $C$-sequence-system, $\umu{}(\uy)$ is a tuple of $\LRC$-terms, and $\varphi(\uy)$ is a conjunction of $\LRC$-equations.
    Given any $L$-formula $\psi(\uxvec; \uw)$, there is another $L$-formula $\psi'(\uxvec{}'; \uw\tiluw)$ that is bounded by $S$, and a tuple of $\LRC$-terms $\ut(\uy)$ such that
    $$
    \exists \ux \in \VV : \psi_\theta(\ux; \uw) \wedge E(\ux; \uy) \quad \equiv \quad \varphi(\uy) \wedge \exists \ux' \in \VV : \psi'_\theta(\ux'; \uw\ut(\uy)) \wedge S(\ux'; \umu(\uy))
    $$
    holds in $T_\theta \cup \set{\text{``$\theta$ is $C$-image-complete''}}$.
\begin{proof}
    Let $S(\ux'; \tiluy) = \bigwedge_{k=1}^n f_k^{q_k}[\theta](x'_{\ldd, k}) = \Tilde{y}{}_k$ be as in Definition \ref{def_param_c_sequence_system}.
    As in this definition, we partition $\ux' = \ux'_\li \ux'_\ld$, where $\ux'_\li := (x'_{\lii, l} : 1 \leq l \leq m)$ and $\ux'_\ld := (x'_{\ldd, l} : 1 \leq l \leq n)$.
    Also, write $\ux = (x_1, \dots, x_r)$, and let ($\underline{\nu}(\ux)$ and) $\utau(\ux'; \uy) = (\tau_1(\ux'; \uy), \dots, \tau_r(\ux'; \uy))$ witness that $E(\ux; \uy)$ is transformable into $\varphi(\uy) \wedge S(\ux'; \umu(\uy))$ (see Definition \ref{def_transfo}).
    By the definition of transformability, each $\tau_k(\ux'; \uy)$ can be written as the sum of an $\LKThe$-term in $\ux'$ and an $\LRC$-term in $\uy$.
    Hence, the same is true for $\theta^i(\tau_k(\ux'; \uy))$.
    Using Lemma \ref{lemma_bound_term_light} (or just \remainderRules{} with the equations $f_k^{q_k}[\theta](x'_{\ldd, k}) = \mu_k(\uy)$), we obtain
    \begin{align}
        \theta^i(\tau_k(\ux'; \uy)) = \sum\nolimits_{l=1}^m \rho_{k, i, \li, l}[\theta](x'_{\li, l}) + \sum\nolimits_{l=1}^n \rho_{k, i, \ld, l}[\theta](x'_{\ld, l}) + t'_{k, i}(\uy) \label{tag_eqeq}
    \end{align}
    modulo the theory $T_\theta \cup \set{\text{``$\theta$ is $C$-image-complete''}} \cup \set{S(\ux'; \umu(\uy))}$ for polynomials that satisfy \hbox{$\deg(\rho_{k, i, \ldd, l}) < \deg(f_l^{q_l})$} and some $\LRC$-term $t'_{k, i}(\uy)$.
    Let $s_k^i(\ux'; \uy)$ denote the right-hand side of (\ref{tag_eqeq}).
    Now we define the $L$-formula $\psi'(\uxvec{}'; \uw\tiluw)$ and the tuple of $\LRC$-terms $\ut(\uy)$ such that $\psi'_\theta(\ux'; \uw \ut(\uy))$ is $\psi_\theta(\ux;\uw)$ with every occurrence of $\theta^i(x_k)$ replaced by $s^i_k(\ux'; \uy)$ for all $i \in \omega$ and $k \in \set{1, \dots, r}$.
    More precisely: We define the index set \hbox{$\ii \!:=\! \set{(k, i) \!:\! \text{``$x^i_k$ appears in $\psi(\uxvec; \uw)$''}}$}, set $\tiluw := (\Tilde{w}{}_{k, i} : (k, i) \in \ii)$ and $\ut(\uy) := (t'_{k, i}(\uy) : (k, i) \in \ii)$, and define $\psi'(\uxvec'; \uw\tiluw)$ as $\psi(\uxvec; \uw)$ with every instance of $x^i_k$ replaced by
    $$
    \sum\nolimits_{l=1}^m \sum\nolimits_{j=0}^{\deg(\rho_{k, i, \li, l})} (\rho_{k, i, \li, l})_j \cdot x_{\li, l}'^j + \sum\nolimits_{l=1}^n \sum\nolimits_{j=0}^{\deg(\rho_{k, i, \ld, l})} (\rho_{k, i, \ld, l})_j \cdot x_{\ld, l}'^j + \Tilde{w}{}_{k, i}.
    $$
    This is an $L$-formula bounded by $S(\ux'; \Tilde{\uy})$ (see Definition \ref{def_formual_bounded}) since \hbox{$\deg(\rho_{k, i, \ldd, l}) < \deg(f^{q_l}_l)$} holds by definition of the $\rho_{k, i, \ldd, l}$'s.

    Since $E(\ux; \uy)$ is transformable into $\varphi(\uy) \wedge S(\ux'; \umu{}(\uy))$, witnessed by (some $\underline{\nu}(\ux)$ and) $\utau(\ux'; \uy)$, the implications in Definition \ref{def_transfo} immediately yield
    $$
    E(\ux; \uy) \equiv \exists \ux' \in \VV : \varphi(\uy) \wedge S(\ux'; \umu{}(\uy)) \wedge \ux = \utau(\ux'; \uy).
    $$
    Using this equivalence, basic first-order logic, and (\ref{tag_eqeq}), we see that
    \begin{align*}
        \exists \ux \in \VV : \psi_\theta(\ux; \uw) \wedge E(\ux; \uy) \quad &\equiv \quad
    \varphi(\uy) \wedge \exists \ux' \in \VV : \psi_\theta(\utau(\ux'; \uy); \uw) \wedge S(\ux'; \umu{}(\uy)) \\
    & \equiv \quad \varphi(\uy) \wedge \exists \ux' \in \VV : \psi'_\theta(\ux'; \uw\ut(\uy)) \wedge S(\ux'; \umu{}(\uy)).
    \end{align*}
    This holds in $T_\theta \cup \set{\text{``$\theta$ is $C$-image-complete''}}$.
\end{proof}
\end{lemma}

\noindent Next, we introduce two quantities that measure the ``degrees of freedom'' of the solutions to a parametrized $C$-sequence-system.

\begin{definition} \label{def_rk_deg_def}
    Let $S(\ux; \uy) = \bigwedge_{k=1}^n f_k^{q_k}[\theta](x_{\ldd, k}) = y_k$ be a parametrized $C$-sequence-system, as in Definition \ref{def_param_c_sequence_system} (i.e., $\ux = \ux_\li\ux_\ld$ and so on).
    We define:
    \begin{enumerate}[(i)]
        \item The \textbf{rank} of $S(\ux; \uy)$ as $\rk(S) := m = |\ux_\li|$.
        \item The \textbf{degree} of $S(\ux; \uy)$ as $\deg(S) := \sum\nolimits_{k=1}^n \deg(f_k^{q_k})$.
    \end{enumerate}
    We define the rank and degree of a $C$-sequence-system (i.e., a parametrized $C$-sequence-system with compatible parameters plugged in) to be the rank and degree of the underlying parametrized $C$-sequence-system.
\end{definition}

\noindent So, for any parametrized $C$-sequence-system, $(\rk(S), \deg(S))$ lies in $\NN^2$.
One should note that $(\NN^2, <_\Lex)$ is a well-ordered set, where $<_\Lex$ is the lexicographical order that treats the first entry as more significant than the second entry, i.e., $(1, 2) <_{\Lex} (2, 1)$.

\begin{definition}
    \label{def_compatible_pair} Let $S(\ux; \tiluy) = \bigwedge_{k=1}^n f_k^{q_k}[\theta](x_{\ldd, k}) = \Tilde{y}{}_k$ be a parametrized $C$-sequence-system as in Definition \ref{def_param_c_sequence_system}.
    Let $\varphi(\uy)$ be a conjunction of $\LRC$-equations, and let $\umu{}(\uy)$ be a tuple of $\LRC$-terms, where $\uy$ is another tuple of variables.
    If
    $$
    \TKvsThe \cup \set{\text{``$\theta$ is $C$-image-complete''}} \models \forall \uy \in \VV : \varphi(\uy) \rightarrow \text{``$\umu{}(\uy)$ is compatible with $S$''}
    $$
    holds, we say that the pair $(\varphi, \umu{})$ is \textbf{compatible} with the parametrized $C$-sequence-system $S$.
\end{definition}

\begin{observation} \label{observation_merge_pram_c_ss}
    Let $S(\ux_1; \uy{}_1), \dots, S(\ux_q; \tiluy_q)$ be parametrized $C$-sequence-systems with each $\ux_i = \ux_{\li, i}\ux_{\ld, i}$ as in Definition \ref{def_param_c_sequence_system} (and with all tuples of variables disjoint).
    Define the tuples $\ux_\li := \ux_{\li, 1}\dots\ux_{\li, q}$, $\ux_\ld := \ux_{\ld, 1}\dots\ux_{\ld, q}$, $\ux := \ux_\li\ux_\ld$, and $\tiluy := \tiluy{}_1\dots\tiluy{}_q$.
    Then the conjunction
    $$
    S(\ux; \tiluy) := \bigwedge\nolimits_{i=1}^q S(\ux_i; \tiluy_i)
    $$
    is also a parametrized $C$-sequence-system.
    Additionally:
    \begin{enumerate}[(i)]
        \item $(\rk(S), \deg(S)) = \sum_{i=1}^q(\rk(S_i), \deg(S_i))$ holds.
        \item For each $i \in \set{1, \dots, q}$, let $(\varphi_i, \umu{}_i)$ be a pair that is compatible with the parametrized $C$-sequence-system $S_i$. Then, with $\varphi(\uy) := \bigwedge_{i=1}^q \varphi_i(\uy)$ and $\umu{}(\uy) := \umu{}_1(\uy)\dots\umu{}_q(\uy)$, the pair $(\varphi, \umu{})$ is compatible with $S$.
    \end{enumerate}
\begin{proof}
    This follows from the definition of a parametrized $C$-sequence-system (Definition \ref{def_param_c_sequence_system}).
    Point (i) follows from Definition \ref{def_rk_deg_def}, and (ii) follows from Definition \ref{def_compatible_pair}, or the reminder above.
\end{proof}
\end{observation}

\noindent Now we take the simplest non-trivial conjunction of $\LRC$-equations and show that it is essentially transformable into a parametrized $C$-sequence-system:

\begin{lemma} \label{lemma_single_eq_trafo}
    Let $\zeta \in K[X]$ be any polynomial and let $\mu(\uy)$ be an $\LRC$-term.
    If $C$ is algebraic, we additionally require $\zeta \neq 0$.
    Then the formula $\zeta[\theta](x) = \mu(\uy)$ is transformable into a formula of the form
    $$
    \varphi(\uy) \wedge S(\ux; \umu{}(\uy)),
    $$
    where $S(\ux; \tiluy)$ is a parametrized $C$-sequence-system, and $(\varphi, \umu{})$ is a pair compatible with $S$.
    Furthermore, $(\rk(S), \deg(S)) \leq_\Lex (0, \deg(\zeta))$ holds if $\zeta \in K[X]\setminus \set{0}$, and otherwise the equality \hbox{$(\rk(S), \deg(S)) = (1, 0)$} holds.
\begin{proof}
    We first assume $\zeta \neq 0$.
    We write $\zeta = \xi_0 \cdot \prod\nolimits_{f \in F} f^{q_f} \cdot \prod\nolimits_{k=1}^n h_k^{r_k}$ with \hbox{$\Fac(\xi_0) \subseteq \Kp{C=0}$}, with $F \subseteq \Kp{0<C<\infty}$ finite, with distinct $h_k \in \Kp{C=\infty}$, and with all $q_f$'s and $r_k$'s positive.
    Now define
    \begin{align*}
        \underline{\nu}(x) &:= \nu_F(x)(\nu_f(x) : f \in F) := \pi_{\Image(F^C)}(x)(\pi_{\Ker(f^C)}(x) : f \in F)
     \\
    & \hspace{160pt}\text{and}\quad\tau(x_F(x_f : f \in F)) := x_F + \sum\nolimits_{f\in F} x_f.
    \end{align*}
    We now show that $\underline{\nu}(x)$ and $\tau(x_F(x_f : f \in F))$ witness that $\zeta[\theta](x) = \mu(\uy)$ is transformable into a formula of the form
    \begin{align*}
        \Big(\prod\nolimits_{k=1}^n h_k^{r_k}\Big)[\theta](x_F) = \mu_F(\uy) \wedge \bigwedge\nolimits_{f\in F} f^{q'_f}[\theta](x_f) = \mu_f(\uy)
    \end{align*}
    with $0 < q'_f \leq C(f)$ and $\mu_f(\uy) \in \Ker(f^C)$ for every $f \in F$.
    Here $\mu_f(\uy) \in \Ker(f^C)$ means that $\mu_f(\uu) \in \Ker(f^C)$ for any tuple $\uu$ in any model of $\TKvsThe \cup \set{\text{``$\theta$ is $C$-image-complete''}}$.
    We start with the first implication of Definition \ref{def_transfo}.
    First, use the \projectionRule{} (Lemma \ref{lemma_projection_rule}) to show that
    \begin{align*}
        \zeta[\theta](x) = \mu(\uy) &\equiv \zeta[\theta](\pi_{\Image(F^C)}(x)) = \pi_{\Image(F^C)}(\mu(\uy))  \wedge \bigwedge\nolimits_{f\in F}\! \zeta[\theta](\pi_{\Ker(f^C)}(x)) = \pi_{\Ker(f^C)}(\mu(\uy)).
    \end{align*}
    Now one can apply Lemma \ref{lemma_communative}, the commutativity of $R_C$, \hbox{$\Ker(f^C) = \Ker(f^{C+1})$}, and (v) of Corollary \ref{corollary_polynomials_kernel_facts}, to show that the right-hand side above is equivalent to
    \begin{align*}
        &\Big(\prod\nolimits_{k=1}^n h_k^{r_k}\Big)[\theta](\pi_{\Image(F^C)}(x)) = \Big(\xi_0 \cdot \prod\nolimits_{f \in F} f^{q_f} \Big)[\theta]^{-1} \circ \pi_{\Image(F^C)}(\mu(\uy)) \\
        &\quad\quad\wedge\bigwedge\nolimits_{f\in F} f^{\min(q_f, C(f))}[\theta](\pi_{\Ker(f^C)}(x)) = \underbrace{\big(\zeta/ f^{q_f}\big)[\theta]_{\restriction \Ker(f^C)}^{-1}}_{=\chi[\theta]_{\restriction \Ker(f^C)} \text{ for some }\chi \in K[X]\hspace{-50pt}} \circ\; \pi_{\Ker(f^C)}(\mu(\uy)).
    \end{align*}
    Since the $x = \tau(\underline{\nu}(x))$ part of the first implication of Definition \ref{def_transfo} is clear by the definition of $\underline{\nu}(x)$ and $\tau(x_F(x_f : f \in F))$, we have shown that implication.
    The second implication of Definition \ref{def_transfo} follows by the same arguments in reverse.

    By the transitivity and modularity of transformability (see (i) and (ii) of Observation \ref{obser_trafo_nice}), and Observation \ref{observation_merge_pram_c_ss} above, we only need to show the lemma for formulas of the following forms:
    \begin{enumerate}[(i)]
        \item $\big(\prod\nolimits_{k=1}^n h_k^{r_k}\big)[\theta](x) = \mu(\uy)$, with all $h_k \in \Kp{C=\infty}$ and all $r_k > 0$.
        \item $f^{q}[\theta](x) = \mu(\uy)$, with $0 < q \leq C(f)$ and $\mu(\uy) \in \Ker(f^C)$.
    \end{enumerate}
    For case (i), notice that, by Example \ref{exampl_trafo}, $\big(\prod\nolimits_{k=1}^n h_k^{r_k}\big)[\theta](x) = \mu(\uy)$ is transformable into $\bigwedge_{k=1}^n h_k^{r_k}[\theta](x_k) = \mu(\uy)$.
    Now $\bigwedge_{k=1}^n h_k^{r_k}[\theta](x_k) = y_k$ is a parametrized $C$-sequence-system which is compatible with any parameters, as all $h_k \in \Kp{C=\infty}$.

    In case (ii), the formula $f^{q}[\theta](x) = \mu(\uy)$ is already of the form $S(x; \mu(\uy))$ for the parametri\-zed $C$-sequence-system $S(x; y) := f^q[\theta](x) = y$.
    In this case,  $\varphi(\uy) := f^{C-q}[\theta](\mu(\uy)) = 0$ implies that $\mu(\uy)$ is compatible with $S$, and holds if and only if $f^{q}[\theta](x) = \mu(\uy)$ is consistent (as we have $\mu(\uy) \in \Ker(f^C)$ and $\Ker(f^C) = \Ker(f^{C+1})$).
    Thus, $f^{q}[\theta](x) = \mu(\uy)$ is equivalent to (and hence transformable into) $\varphi(\uy) \wedge S(x; \mu(\uy))$.

    For the furthermore part, notice that the final parametrized $C$-sequence-system $S(\ux; \uy)$ obtained with Observation \ref{observation_merge_pram_c_ss} has $\rk(S) = 0$ and that we also obtain
    $$
    \deg(S) = \deg\Big(\prod\nolimits_{k=1}^n h_k^{r_k} \cdot \prod\nolimits_{f\in F} f^{\min(q_f, C(f))}\Big) \leq \deg(\zeta).
    $$
    Finally, if $\zeta = 0$, then $\zeta[\theta](x) = \mu(\uy)$ is equivalent to $\varphi(\uy) \wedge S(x)$ with \hbox{$\varphi(\uy) := \mu(\uy) = 0$} and $S(x) := \top$ being the unique parametrized $C$-sequence-system with $(\rk(S), \deg(S)) = (1, 0)$.
    This is a parametrized $C$-sequence-system, as we assume $C$ to be transcendental in the case $\zeta = 0$.
    As this $S$ is compatible with the empty tuple, we conclude.
\end{proof}
\end{lemma}

\noindent Next, we show that the conjunction of a parametrized $C$-sequence-system and some $\LKThe$-equations is again essentially transformable into a parametrized $C$-sequence-system with lexicographically smaller or equal rank and degree.
We also provide a criterion for when the rank and degree become smaller.
This is actually the setting we mostly care about, as our theory $T$ can only state that certain sets imply a finite disjunction of $K$-linear equations.

\begin{lemma}[Transformation Lemma]\label{lemma_trafoo}
    Let $S(\ux; \uy{}_1)$ be a parametrized $C$-sequence-system, let $E(\ux; \uy{}_2)$ be a conjunction of $\LKThe$-equations, and let $\umu{}_1(\uy)$ and $\umu{}_2(\uy)$ be tuples of $\LRC$-terms.
    Then there is another parametrized $C$-sequence-system $S'(\ux'; \uy')$ and a pair $(\varphi', \umu{}')$ compatible with $S'$ such that
    $$
    \text{$S(\ux; \umu{}_1(\uy)) \wedge E(\ux; \umu{}_2(\uy))\quad$ is transformable into $\quad\varphi'(\uy) \wedge S'(\ux'; \umu{}'(\uy))$}.
    $$
    Furthermore, $(\rk(S'), \deg(S')) \leq_\Lex (\rk(S), \deg(S))$ holds, and if $E(\ux; \uy{}_2)$ contains at least one equation that is non-trivial in $\ux$ and bounded by $S$, then this inequality is strict.
\begin{proof}
    We deviate from our usual notation in Definition \ref{def_param_c_sequence_system} for the parametrized $C$-sequence-system $S(\ux; \uy{}_1)$.
    First, we write the tuple $\ux = \ux_\li\ux_\ld$ with $\ux_\li = (x_1, \dots, x_m)$ and \hbox{$\ux_\ld = (x_{m+1}, \dots, x_{m+n})$} instead of naming the entries of these tuples $x_{\lii,k}$ and $x_{\ldd,k}$.
    Second, we write
    $$
    S(\ux; \umu{}_1(\uy)) = \bigwedge\nolimits_{k=1}^n \xi_k[\theta](x_{m+k}) = \mu{}_{1, k}(\uy),
    $$
    where each $\xi_k$ is of the form $f^{q_k}$ with $C(f) > 0$ and $q_k \!\in\! \set{q \in \NN : 0 \!<\! q\! \leq\! C(f)}$.
    Without loss, we can also assume that
    $$
    E(\ux; \umu{}_2(\uy)) = \bigwedge\nolimits_{k=1}^r \sum\nolimits_{l=1}^{m+n} \rho_{k, l}[\theta](x_l) = \mu_{2, k}(\uy),
    $$
    with all $\rho_{k, l}$'s being polynomials over $K$.
    We can now express $S(\ux; \umu{}_1(\uy)) \wedge E(\ux; \umu{}_2(\uy))$ in the language of left $K[X]$-modules, with \hbox{$\xi \cdot x := \xi[\theta](x)$}, as follows:
    $$
    \begin{tikzpicture}[baseline=(Frame.base)]
    \drawText{0}{0}{0}
    \drawText{0}{3}{0}
    \drawText{3}{0}{0}
    \drawText{3}{3}{0}
    \drawHDots{1}{0}{2}
    \drawHDots{1}{3}{2}
    \drawVDots{0}{1}{2}
    \drawVDots{3}{1}{2}

    \drawText{4}{1}{0}
    \drawText{4}{3}{0}
    \drawText{6}{3}{0}
    \drawHDots{5}{3}{1}
    \drawVDots{4}{2}{1}
    \drawDDots{5}{2}{1}

    \drawText{4}{0}{\xi\sindex{1}}
    \drawText{7}{3}{\xi\sindex{n}}
    \drawDDots{5}{1}{2}

    \drawText{5}{0}{0}
    \drawText{7}{0}{0}
    \drawText{7}{2}{0}
    \drawHDots{6}{0}{1}
    \drawVDots{7}{1}{1}
    \drawDDots{6}{1}{1}

    \drawText{0}{4}{\rho\sindex{1,\!1}}
    \drawText{0}{6}{\rho\sindex{r,\! 1}}
    \drawText{3}{4}{\rho\sindex{1,\!m}}
    \drawText{3}{6}{\rho\sindex{r,\! m}}
    \drawText{4}{4}{\rho\sindex{1,\!m\!+\!1}}
    \drawText{4}{6}{\rho\sindex{r,\!m\!+\!1}}
    \drawText{7}{4}{\rho\sindex{1,\! m\!+\!n}}
    \drawText{7}{6}{\rho\sindex{r,\! m\!+\!n}}
    \drawHDots{1}{4}{2}
    \drawHDots{1}{6}{2}
    \drawHDots{5}{4}{2}
    \drawHDots{5}{6}{2}
    \drawVDots{0}{5}{1}
    \drawVDots{3}{5}{1}
    \drawVDots{4}{5}{1}
    \drawVDots{7}{5}{1}

    \drawHLine{0}{4}{8}
    \drawVLine{4}{0}{7}

    \drawBorder{0}{0}{8}{7}
    \end{tikzpicture}
    \;\cdot\;
    \begin{tikzpicture}[baseline=(Frame.base)]
    \drawText{0}{0}{x\sindex{1}}
    \drawText{0}{3}{x\sindex{m}}
    \drawText{0}{4}{x\sindex{m\!+\!1}}
    \drawText{0}{7}{x\sindex{m\!+\!n}}
    \drawVDots{0}{1}{2}
    \drawVDots{0}{5}{2}
    \drawHLine{0}{4}{1}
    \drawBorder{0}{0}{1}{8}
    \end{tikzpicture}
    \;=\;
    \begin{tikzpicture}[baseline=(Frame.base)]
    \drawText{0}{0}{\mu\sindex{1,\!1}(\uy)}
    \drawText{0}{3}{\mu\sindex{1,\!n}(\uy)}
    \drawText{0}{4}{\mu\sindex{2,\!1}(\uy)}
    \drawText{0}{6}{\mu\sindex{2,\!r}(\uy)}
    \drawVDots{0}{1}{2}
    \drawVDots{0}{5}{1}
    \drawHLine{-0.1}{4}{1.2}
    \drawBorder{-0.1}{0}{1.1}{7}
    \end{tikzpicture}
    $$
    where each row above the horizontal line represents a single equation in $S(\ux; \umu{}_1(\uy))$, and each row below the horizontal line represents one equation in $E(\ux; \umu{}_2(\uy))$.
    Treating finite tuples as column matrices, we write the above equation as \hbox{$A \cdot \ux = \umu(\uy)$}, where $A \in K[X]^{(n+r) \times (m+n)}$ denotes the big matrix, and \hbox{$\umu{}(\uy) := \umu{}_1(\uy)\umu{}_2(\uy)$}.
    Notice that, in general, $n+r$ can be greater than, smaller than, or equal to $m+n$.

    By Theorem 3.8 in \cite{Jac85}, we can find two invertible matrices $L \in K[X]^{(n+r)\times (n+r)}$ and $R \in K[X]^{(m+n) \times (m+n)}$ such that the product $L \cdot A \cdot R$ is a diagonal matrix, or more preceisly the Smith normal form of $A$.
    Defining \hbox{$N := \min\set{n+r, m+n}$}, we obtain
    $$
    D := L \cdot A \cdot R = \begin{tikzpicture}[baseline=(Frame.base)]

    \drawText{0}{1}{0}
    \drawText{0}{3}{0}
    \drawText{2}{3}{0}
    \drawHDots{1}{3}{1}
    \drawVDots{0}{2}{1}
    \drawDDots{1}{2}{1}

    \drawText{0}{0}{\zeta\sindex{1}}
    \drawText{3}{3}{\zeta\sindex{N}}
    \drawDDots{1}{1}{2}

    \drawText{1}{0}{0}
    \drawText{4}{3}{0}
    \drawText{6}{0}{0}
    \drawText{6}{3}{0}
    \drawHDots{2}{0}{4}
    \drawHDots{5}{3}{1}
    \drawVDots{6}{1}{2}
    \drawDDots{2}{1}{2}

    \drawBorder{0}{0}{7}{4}
    \end{tikzpicture}\hspace{3pt}
    $$
    (with the visualization illustrating the case $N = n+r < m+n$).
    Now define the tuple \hbox{$\ux_* := (x_{*,1}, \dots, x_{*,m+n})$} and define \hbox{$\umu{}_*(\uy) := (\mu_{*,1}(\uy), \dots, \mu_{*,n+r}(\uy)) := L \cdot \umu{}(\uy)$}.
    Using the invertibility of $L$ and $R$, we see that the maps $\underline{\nu}(\ux) := R^{-1} \cdot \ux$ and $\utau(\ux_*) := R \cdot \ux_*$ witness that the equation $A \cdot \ux = \umu{}(\uy)$, and therefore $S(\ux; \umu{}_1(\uy)) \wedge E(\ux; \umu{}_2(\uy))$, is transformable into the equation $D \cdot \ux_* = \umu{}_*(\uy)$.
    Switching back to our usual language, $D \cdot \ux_* = \umu{}_*(\uy)$ is just the conjunction
    $$
    \bigwedge\nolimits_{k=1}^N \zeta_k[\theta](x_{*,k}) = \mu_{*,k}(\uy) \wedge \bigwedge\nolimits_{k=N+1}^{n+r} \mu_{*,k}(\uy) = 0.
    $$
    We additionally define $\zeta_k := 0$ and $\mu_{*,k}(\uy) := 0$ for $k > N$ in the case $n+r < m+n$.
    Then, in any case, the above is obviously equivalent to
    $$
    \bigwedge\nolimits_{k=1}^{m+n} \zeta_k[\theta](x_{*,k}) = \mu_{*,k}(\uy) \wedge \underbrace{\bigwedge\nolimits_{k=m+n+1}^{n+r} \mu_{*,k}(\uy) = 0}_{=: \varphi_0(\uy)},
    $$
    where the second conjunction is needed only in the case $n+r > m+n$.
    Thus, by (iii) of Observation \ref{obser_trafo_nice}, $D \cdot \ux_* = \umu{}_*(\uy)$ is transformable into $\varphi_0(\uy) \wedge \bigwedge_{k=1}^{m+n} \zeta_k[\theta](x_{*,k}) = \mu_{*,k}(\uy)$.
    By the transitivity of transformability ((i) of Observation \ref{obser_trafo_nice}), we see that the conjunction $S(\ux; \umu{}_1(\uy)) \wedge E(\ux; \umu{}_2(\uy))$ is transformable into $\varphi_0(\uy) \wedge \bigwedge_{k=1}^{m+n} \zeta_k[\theta](x_{*,k}) = \mu_{*,k}(\uy)$.

    Now, by Lemma \ref{lemma_single_eq_trafo}, each formula $\zeta_k[\theta](x_{*,k}) = \mu_{*,k}(\uy)$ is transformable into a formula of the form
    $$
    \varphi_k(\uy) \wedge S_k(\ux_k; \umu{}_{\stri, k}(\uy))
    $$
    where $S_k(\ux_k; \uy_k)$ is a parametrized $C$-sequence-system and $(\varphi_k, \umu{}_{\stri, k})$ is a pair compatible with $S_k$.
    We will later argue why all $\zeta_k$'s are non-zero if $C$ is algebraic.
    Additionally, we obtain the equation $(\rk(S_k), \deg(S_k)) = (1, 0)$ if $\zeta_k = 0$, and otherwise we obtain $(\rk(S_k), \deg(S_k)) \leq_\Lex (0, \deg(\zeta_k))$.

    With Observation \ref{observation_merge_pram_c_ss} and the modularity of transformability (see (ii) of Observation \ref{obser_trafo_nice}), we see that $\varphi_0(\uy) \wedge \bigwedge_{k=1}^{m+n} \zeta_k[\theta](x_{*,k}) = \mu_{*,k}(\uy)$ is transformable into a formula of the form
    $$
    \varphi'(\uy) \wedge S'(\ux'; \umu{}'(\uy))
    $$
    where $S'(\ux'; \uy')$ is another parametrized $C$-sequence-system, and $(\varphi', \umu{}')$ is a pair compatible with $S'$.
    Note that when applying the modularity of transformability, we treat $\varphi_0(\uy)$ as a single formula, which is obviously transformable into itself.
    Additionally,
    $$
    (\rk(S'), \deg(S')) = \sum\nolimits_{k=1}^{n+m}(\rk(S_k), \deg(S_k))
    $$
    holds.
    By the transitivity of transformability ((i) of Observation \ref{obser_trafo_nice}), we conclude that $S(\ux; \umu{}_1(\uy)) \wedge E(\ux; \umu{}_2(\uy))$ is transformable into the formula $\varphi'(\uy) \wedge S'(\ux'; \umu{}'(\uy))$.
    It remains to show the furthermore part, and to show that all $\zeta_k$'s are indeed non-zero in the algebraic case.

    By the proof of Theorem 3.9 in \cite{Jac85}, we see that for $k \leq N$, the product $\prod_{l=1}^k \zeta_l$ is, up to some factor in $K^\times$, the greatest common divisor of all $k$-rowed minors of $A$, i.e., the determinants of all $k \times k$ submatrices of $A$ obtained by omitting rows and columns.
    Recall that for $k > N$, each $\zeta_k$ is $0$, so it is easy to see that whenever $\zeta_k = 0$ and $l > k$, then $\zeta_l = 0$ also holds.
    Now let $q$ be maximal with $\zeta_q \neq 0$ (or $q = 0$ if all $\zeta_k = 0$).
    We obtain
    $$
    (\rk(S'), \deg(S')) \leq_\Lex (m+n-q, \deg(\zeta_1 \cdots \zeta_q))
    $$
    with our formula for $(\rk(S'),\hspace{-0.3pt} \deg(S'))$ and the inequality for each $(\rk(S_k),\hspace{-0.3pt} \deg(S_k))$.
    By taking a look at the upper-right submatrix of $A$, we see that $\xi_1 \cdots \xi_n \neq 0$ is an $n$-rowed minor of $A$, so we must have $q \geq n$.
    In the algebraic case, we have $m = 0$ by Definition \ref{def_param_c_sequence_system}, so $n = m+n$, and therefore $\zeta_k \neq 0$ follows for all $k$.
    Back in the general case, $\rk(S') \leq \rk(S) = m$ follows from $q \geq n$.
    If $\rk(S') = \rk(S)$, then $q = n$.
    Therefore $\zeta_1 \cdots \zeta_q \mid \xi_1 \cdots \xi_n$, and hence $\deg(S') \leq \deg(\xi_1 \cdots \xi_n) = \deg(S)$.
    We conclude
    $$
    (\rk(S'), \deg(S')) \leq_\Lex (\rk(S), \deg(S)).
    $$
    Now assume that $E(\ux; \uy{}_2)$ contains a non-trivial equation in $\ux$ that is bounded by $S$.
    In this case, we need to show $(\rk(S'), \deg(S')) <_\Lex (\rk(S), \deg(S))$.
    By the definition of boundedness (see Definition \ref{def_formual_bounded}), this means that there must be a $k \in \set{1, \dots, r}$ such that $\deg(\rho_{k, l}) < \deg(\xi_{l-m})$ holds for all $l > m$, and such that at least one $\rho_{k, l}$ is non-zero.
    Fix such a $k$ and let $s$ be maximal with $\rho_{k, s} \neq 0$.
    We distinguish between the cases $s \leq m$ and $s > m$.

    First, assume $s \leq m$.
    In this case, we obtain the following $(n+1)\times(n+1)$ submatrix of $A$ by omitting rows and columns:
    $$
    \begin{tikzpicture}[baseline=(Frame.base)]
    \drawText{0}{0}{0}
    \drawText{0}{3}{0}
    \drawVDots{0}{1}{2}

    \drawText{1}{1}{0}
    \drawText{1}{3}{0}
    \drawText{3}{3}{0}
    \drawHDots{2}{3}{1}
    \drawVDots{1}{2}{1}
    \drawDDots{2}{2}{1}

    \drawText{1}{0}{\xi\sindex{1}}
    \drawText{4}{3}{\xi\sindex{n}}
    \drawDDots{2}{1}{2}

    \drawText{2}{0}{0}
    \drawText{4}{0}{0}
    \drawText{4}{2}{0}
    \drawHDots{3}{0}{1}
    \drawVDots{4}{1}{1}
    \drawDDots{3}{1}{1}

    \drawText{0}{4}{\rho\sindex{k,\!s}}
    \drawText{1}{4}{0}
    \drawText{4}{4}{0}
    \drawHDots{2}{4}{2}

    \drawHLine{0}{4}{5}
    \drawVLine{1}{0}{5}

    \drawBorder{0}{0}{5}{5}
    \end{tikzpicture}\;.
    $$
    As $\rho_{k, s} \neq 0$, we see that $A$ has a non-zero $(n+1)$-rowed minor.
    Hence $q > n$, and therefore $\rk(S') < \rk(S)$, as $\rk(S') = m + n - q$ and $\rk(S) = m$.

    Now assume $s > m$.
    If $q > n$, we conclude as above, so we may assume $q = n$.
    We obtain the following $n \times n$ submatrix of $A$ by omitting rows and columns:
    $$
    \begin{tikzpicture}[baseline=(Frame.base)]
    \drawText{0}{1}{0}
    \drawText{0}{3}{0}
    \drawText{2}{3}{0}
    \drawText{1}{0}{0}
    \drawText{3}{2}{0}
    \drawText{4}{3}{0}
    \drawText{5}{3}{0}
    \drawText{3}{0}{0}
    \drawText{4}{0}{0}
    \drawText{5}{0}{0}
    \drawText{8}{0}{0}
    \drawText{8}{3}{0}
    \drawHDots{2}{0}{1}
    \drawHDots{6}{3}{2}
    \drawHDots{6}{0}{2}
    \drawHDots{1}{3}{1}
    \drawVDots{0}{2}{1}
    \drawVDots{3}{1}{1}
    \drawVDots{4}{1}{2}
    \drawVDots{5}{1}{2}
    \drawVDots{8}{1}{2}
    \drawDDots{1}{2}{1}
    \drawDDots{2}{1}{1}

    \drawText{0}{0}{\xi\sindex{1}}
    \drawText{3}{3}{\xi\sindex{s\sminus m \sminus 1}}
    \drawText{5}{4}{\xi\sindex{s\sminus m \!+\! 1}}
    \drawText{8}{7}{\xi\sindex{n}}
    \drawDDots{1}{1}{2}
    \drawDDots{6}{5}{2}

    \drawText{8}{6}{0}
    \drawText{8}{4}{0}
    \drawText{6}{4}{0}
    \drawText{7}{7}{0}
    \drawText{3}{4}{0}
    \drawText{4}{4}{0}
    \drawText{5}{5}{0}
    \drawText{3}{7}{0}
    \drawText{4}{7}{0}
    \drawText{5}{7}{0}
    \drawText{0}{7}{0}
    \drawText{0}{4}{0}
    \drawHDots{1}{7}{2}
    \drawHDots{6}{7}{1}
    \drawHDots{1}{4}{2}
    \drawHDots{7}{4}{1}
    \drawVDots{8}{5}{1}
    \drawVDots{0}{5}{2}
    \drawVDots{3}{5}{2}
    \drawVDots{4}{5}{2}
    \drawVDots{5}{6}{1}
    \drawDDots{7}{5}{1}
    \drawDDots{6}{6}{1}

    \drawText{0}{8}{\rho\sindex{k,\!m\!+\!1}}
    \drawText{3}{8}{\rho\sindex{k,\!s\sminus 1}}
    \drawText{4}{8}{\rho\sindex{k,\!s}}
    \drawText{5}{8}{0}
    \drawText{8}{8}{0}
    \drawHDots{1}{8}{2}
    \drawHDots{6}{8}{2}

    \drawHLine{0}{8}{9}

    \drawBorder{0}{0}{9}{9}
    \end{tikzpicture}.
    $$
    This implies that $\pm \prod_{l=1, l \neq s-m}^n \xi_l \cdot \rho_{k, s}$ is an $n$-rowed minor of $A$.
    As $\xi_1 \cdots \xi_n$ is another $n$-rowed minor of our matrix $A$, and $q = n$, we obtain $\zeta_1 \cdots \zeta_q \mid \prod_{l=1, l \neq s-m}^n \xi_l \cdot \gcd(\xi_{s-m}, \rho_{k, s})$.
    Since the polynomial $\rho_{k, s}$ is non-zero and has smaller degree than $\xi_{s-m}$, we obtain the inequality $\deg(\gcd(\xi_{s-m}, \rho_{k, s})) < \deg(\xi_{s-m})$, and therefore
    $$
    \deg(S') \leq \deg(\zeta_1 \cdots \zeta_q) < \deg(\xi_1 \cdots \xi_n) = \deg(S).
    $$
    This implies $(\rk(S'), \deg(S')) <_\Lex (\rk(S), \deg(S))$, as $\rk(S') = \rk(S)$.
\end{proof}
\end{lemma}

\begin{remark}
    One can actually show that, in the \TrafoLemma{}, the equality
    $$
    (\rk(S'), \deg(S')) = (\rk(S), \deg(S))
    $$
    holds if and only if every row of the matrix $A$ is $K[X]$-generated by the first $n$ rows.
    This is the case if and only if $S(\ux; \uzero)$ implies $E(\ux; \uzero)$.
\end{remark}

\noindent We now state another easy fact about transformability:

\begin{lemma} \label{lemma_trafo_top}
    The formula $\top$ is, as a formula in $\ux$ and $\uy$, transformable into a formula of the form
    $$
    S(\ux'; \uzero),
    $$
    where $S(\ux'; \uy')$ is a parametrized $C$-sequence-system, and $\uzero$ is treated as a tuple of $\LKThe$-terms in $\uy$. Note that the pair $(\top, \uzero)$ is compatible with $S$.
\begin{proof}
    In the transcendental case, there is nothing to show, since $\top$, as a formula in \hbox{$\ux_\li := \ux$}, is already a parametrized $C$-sequence-system.
    In the case where $C$ is algebraic, we define $F = \Kp{0<C<\infty}$.
    The tuples $\underline{\nu}(\ux) := (\pi_{\Ker(f^C)}(\ux) : f \in F)$ and
    $\utau((\ux_f : f \in F)) := \sum\nolimits_{f \in F} \ux_f$ witness that $\top$ is transformable into
    $$
    \bigwedge\nolimits_{f \in F} \bigwedge\nolimits_{k=1}^n f^C[\theta](x_{f, k}) = 0,
    $$
    which is a parametrized $C$-sequence-system with $\ux_\li$ empty, $\ux_\ld = (\ux_f : f \in F)$, and $\uzero$ plugged in.
\end{proof}
\end{lemma}

\noindent We would like to give a quick recap of what is shown in this section, offer a broader perspective, and clarify why we are working with transformability here.
Using the \commonBaseTheorem{} (Theorem \ref{theorem_r_c_element}), one can transform any conjunction of $\LRC$-equations $E(\ux; \uy)$ into a conjunction of $\LKThe$-equations with some $\LRC$-terms in $\uy$ plugged in (we will do so in a future paper).
Combining this with (Lemma \ref{lemma_trafo_top} and) our \TrafoLemma{} (Lemma \ref{lemma_trafoo}) and Lemma \ref{lemma_main_trafo}, one can see that, for any $L$-formula $\psi(\uxvec; \uw)$ and any conjunction of $\LRC$-equations $E(\ux; \uy)$,
\begin{align}
    \exists \ux \in \VV : \psi_\theta(\ux; \uw) \wedge E(\ux; \uy) \quad\!\equiv\quad\! \varphi(\uy) \wedge \exists \ux' \in \VV : \psi'_\theta(\ux'; \uw \ut(\uy)) \wedge S(\ux'; \umu{}(\uy)) \label{tag_equiv_pp_1}
\end{align}
holds for some parametrized $C$-sequence-system $S(\ux'; \tiluy)$, another $L$-formula $\psi'(\uxvec{}'; \uw\tiluw)$ that is bounded by $S$, a tuple $\ut(\uy)$ of $\LRC$-terms, and a pair $(\varphi, \umu{})$ compatible with $S$.
If one ignores the $L$-part, then one obtains
\begin{align}
    \exists \ux : E(\ux; \uy) \quad\equiv\quad \varphi(\uy) \wedge \exists \ux' : S(\ux'; \umu{}(\uy)) \label{tag_equiv_pp_2}
\end{align}
in $\TKvsThe \cup \set{\text{``$\theta$ is $C$-image-complete''}}$.
Note that the formulas on both sides are pp-formulas if stated in the language of $R_C$-modules.
Given any ring $R$, a complete theory of $R$-modules has pp-elimination of quantifiers, which means that every formula is equivalent to a Boolean combination of pp-formulas \cite{Bau76}.
If $R$ is a principal ideal domain, these pp-formulas can furthermore be chosen to be of a nice form; see, e.g., Theorem 2.$\ZZ$1 in \cite{Pre88} (but note that the presented proof is attributed to Hodges and the result itself to Eklof and Fisher).
Note that the proof of Theorem 2.$\ZZ$1 in \cite{Pre88} also uses the Smith normal form, as does the proof of the \TrafoLemma{}, and that those two proofs are quite similar apart from all the technical details.
So, just like for modules over a principal ideal domain, we have proven that for our $R_C$-modules (i.e., models of $\TKvsThe \cup \set{\text{``$\theta$ is $C$-image-complete''}}$ treated as $R_C$-modules), any pp-formula is equivalent to a pp-formula of a neat form.
The key difference is that, in the context of modules, one does not need to keep track of the tuples $\underline{\nu}(\ux)$ and $\utau(\ux'; \uy)$ witnessing transformability, since there is no structure outside of the module.
By contrast, we need the tuple $\utau(\ux'; \uy)$ in order to make Lemma \ref{lemma_main_trafo} work.
The mere equivalence in (\ref{tag_equiv_pp_2}) alone does not guarantee that the equivalence in (\ref{tag_equiv_pp_1}) holds.
The careful reader probably noticed that the tuple $\underline{\nu}(\ux)$ is not used anywhere above.
Indeed, for the most part, one could probably replace the first implication in Definition \ref{def_transfo} with 
$$
\forall \ux\uy: E(\ux; \uy) \rightarrow \exists \ux' : E'(\ux'; \uy) \wedge \ux = \utau(\ux'; \uy)
$$
to get rid of the tuple $\underline{\nu}(\ux)$. However, any ``elementary'' proof of transformability (i.e., such as Example \ref{exampl_trafo} or Lemma \ref{lemma_single_eq_trafo}) would still define this tuple $\underline{\nu}(\ux)$ implicitly to show the implication above. Moreover, the tuple $\underline{\nu}(\ux)$ will be used in our preservation of \NATP{} result for our model companions in a future paper, so we keep it in the definition of transformability.

\subsection{Proof of ``$\Leftarrow$''} \label{sec_right_to_left}

\noindent The goal of this section is to prove the following theorem and use it to show the implication from right to left of Theorem \ref{theorem_big_characterization}:

\begin{theorem} \label{theorem_exist_sentence_reduction}
    Let $(\mm, \theta) \models T^C_\theta$ be $C$-image-complete.
    Any existential $L_\theta(M)$-sentence $\phi$ is, modulo $T_\theta \cup \set{\text{``$\theta$ is $C$-image-complete''}} \cup \Diag(\mm, \theta)$,
    equivalent to a finite disjunction of sentences of the form
    $$
    \exists \ux \in \VV : \psi_\theta(\ux) \wedge S(\ux)
    $$
    where $S(\ux)$ is a $C$-sequence-system, and $\psi(\uxvec{})$ is an $L(M)$-formula that is bounded by $S$ and does not imply any finite disjunction of non-trivial linear dependencies in $\uxvec{}$ over $\VV$.
\end{theorem}

\noindent Here $\Diag(\mm, \theta)$ is the \textbf{diagram} of $(\mm, \theta)$, which is the set of all quantifier-free $L_\theta(M)$-sentences that hold in $(\mm, \theta)$.
Hence, the theory
$$
T_\theta \cup \set{\text{``$\theta$ is $C$-image-complete''}} \cup \Diag(\mm, \theta)
$$
is the theory of all $C$-image-complete extensions of $(\mm, \theta)$.

\begin{proof}[Proof of Theorem \ref{theorem_exist_sentence_reduction}]
    We start with two claims:
\begin{subclaim} \label{lemma_exist_sentence_0}
    Modulo $T_\theta \cup \set{\text{``$\theta$ is $C$-image-complete''}} \cup \Diag(\mm, \theta)$, any existential $L_\theta(M)$-sentence $\phi$ is equivalent to a sentence of the form
    $$
    \exists \ux \in \VV : \psi_\theta(\ux) \wedge S(\ux)
    $$
    where $S(\ux)$ is a $C$-sequence-system, and $\psi(\uxvec{})$ is an $L(M)$-formula that is bounded by $S$.
\begin{innerproof}
    Assume $\phi = \phi_*(\ud)$ for some existential $L_\theta$-formula $\phi_*(\uw)$ and some $\ud \in M$.
    By Corollary \ref{corollary_equiv_exist_form}, the formula $\phi_*(\uw)$ is, modulo $T_\theta$, equivalent to a formula of the form
    $$
    \exists \ux \in \VV : \psi_\theta(\ux; \uw) \wedge \top
    $$
    where $\psi(\uxvec; \uw)$ is some $L$-formula.
    By Lemma \ref{lemma_trafo_top}, $\top$ is, as a formula in $\ux$, transformable into a formula of the form $S(\ux'; \uzero)$ for some parametrized $C$-sequence-system $S(\ux'; \tiluy)$.
    Now use this transformability together with Lemma \ref{lemma_main_trafo} to obtain
    $$
    \exists \ux \in \VV : \psi_\theta(\ux; \uw) \wedge \top \quad \equiv \quad \exists \ux' \in \VV : \psi'_\theta(\ux'; \uw\uzero) \wedge S(\ux'; \uzero)
    $$
    modulo $T_\theta \cup \set{\text{``$\theta$ is $C$-image-complete''}}$, for some $L$-formula $\psi'(\uxvec{}'; \uw\tiluw)$ that is bounded by $S$ (since $\uy$ is empty in our application of Lemma \ref{lemma_main_trafo}, we obtain $\ut(\uy) = \uzero$).
    As $\uzero$ is compatible with any parametrized $C$-sequence-system, $S(\ux'; \uzero)$ is a $C$-sequence-system.
    Plugging $\ud$ back in for $\uw$ yields the desired result.
\end{innerproof}
\end{subclaim}

\begin{subclaim} \label{lemma_exist_sentence_1}
    Let $\phi = \exists \ux \in \VV : \psi_\theta(\ux) \wedge S(\ux) \wedge E(\ux)$ be an $L_\theta(M)$-sentence, where $S(\ux)$ is a $C$-sequence-system, $\psi(\uxvec{})$ is an $L(M)$-formula, and $E(\ux)$ is a non-trivial $\LKThe(\VV)$-equation that is bounded by $S$.
    Then $\phi$ is, modulo \hbox{$T_\theta \cup \set{\text{``$\theta$ is $C$-image-complete''}} \cup \Diag(\mm, \theta)$}, either equivalent to $\bot$ or equivalent to a formula of the form
    $$
    \exists \ux' \in \VV :  \psi'_\theta(\ux') \wedge S'(\ux')
    $$
    where $S'(\ux')$ is a $C$-sequence-system satisfying $(\rk(S'), \deg(S')) <_\Lex (\rk(S), \deg(S))$, and $\psi'(\uxvec{}')$ is an $L(M)$-formula that is bounded by $S'$.
\begin{innerproof}
    Since $\psi(\uxvec{})$ is an $L(M)$-formula, we have $\psi(\uxvec{}) = \psi_*(\uxvec{}; \ud)$ for some $L$-formula $\psi_*(\uxvec{}; \uw)$ and tuple $\ud \in M$.
    Similarly, $S(\ux) = S_*(\ux; \uu_1)$ and $E(\ux) = E_*(\ux; \uu_2)$ for some parametrized $C$-sequence-system $S_*(\ux; \uy_1)$, $\LKThe$-equation $E_*(\ux; \uy_2)$, tuple $\uu_1 \in \VV$ that is compatible with $S_*$, and tuple $\uu_2 \in \VV$.
    Set $\uu = \uu_1\uu_2$ and $\uy := \uy_1\uy_2$, and define $\umu{}_1(\uy) := \uy_1$ and $\umu{}_2(\uy) := \uy_2$.
    With this, we obtain the following equivalence in any extension of $(\mm, \theta)$:
    \begin{align}
        \exists \ux \in \VV \!:\! \psi_\theta(\ux) \wedge S(\ux) \wedge E(\ux)  \equiv  \exists \ux \in \VV \!:\! \psi_{*,\theta}(\ux; \ud) \wedge S_*(\ux; \umu{}_1(\uu)) \wedge E_*(\ux; \umu{}_2(\uu)). \label{tag_first_eqeq}
    \end{align}
    Apply the \TrafoLemma{} (Lemma \ref{lemma_trafoo}) to show that $S_*(\ux; \umu{}_1(\uy)) \wedge E_*(\ux; \umu{}_2(\uy))$ is transformable into a formula of the form $\varphi(\uy) \wedge S'_*(\ux'; \umu{}'(\uy))$, where $(\varphi, \umu{}')$ is a pair compatible with the parametrized $C$-sequence-system $S'_*(\ux';\tiluy')$.
    Recall that $(\varphi, \umu{}')$ being compatible with $S'_*$ means that $\varphi(\uy)$ is a conjunction of $\LRC$-equations that implies that $\umu{}'(\uy)$ is compatible with $S'_*$.
    Then apply Lemma \ref{lemma_main_trafo} with this transformability to obtain
    $$
    \exists \ux \in \VV \!:\! \psi_{*,\theta}(\ux; \uw) \wedge S_*(\ux; \umu{}_1(\uy)) \wedge E_*(\ux; \umu{}_2(\uy)) \equiv \varphi(\uy) \wedge \exists \ux' \in \VV\!:\! \psi'_{*,\theta}(\ux'; \uw\ut(\uy)) \wedge S'_*(\ux'; \umu{}'(\uy)),
    $$
    for some $L$-formula $\psi'_*(\uxvec{}'; \uw\tiluw)$ bounded by $S'_*$ and some tuple $\ut(\uy)$ of $\LRC$-terms.
    With (\ref{tag_first_eqeq}), we can now easily see that our original $L_\theta(M)$-sentence $\exists \ux \in \VV : \psi_\theta(\ux) \wedge S(\ux) \wedge E(\ux)$ is equivalent to
    $$
    \varphi(\uu) \wedge \exists \ux' \in \VV: \psi'_{*,\theta}(\ux'; \ud\ut(\uu)) \wedge S'_*(\ux'; \umu{}'(\uu)).
    $$
    By Lemma \ref{lemma_substructure_R_C}, $\varphi(\uu)$ holds if and only if it holds in any $C$-image-complete extension.
    If $\varphi(\uu)$ does not hold, then the sentence is equivalent to $\bot$.
    Otherwise, we can conclude by setting $\psi'(\uxvec{}') := \psi'_*(\uxvec{}'; \ud\ut(\uu))$ and $S'(\ux') := S'_*(\ux';\umu{}'(\uu))$.

    Notice that, by definition, $E_*(\ux; \uy_2)$ is bounded by $S_*$ and non-trivial in $\ux$, since $E(\ux)$ is bounded by $S$ and non-trivial.
    We also have \hbox{$(\rk(S_*), \deg(S_*)) = (\rk(S), \deg(S))$}.
    Hence, the \TrafoLemma{} ensures that the inequality $(\rk(S'), \deg(S')) <_\Lex (\rk(S), \deg(S))$ holds.
\end{innerproof}
\end{subclaim}

\noindent We can now start with the actual proof of Theorem \ref{theorem_exist_sentence_reduction}.
By Claim \ref{lemma_exist_sentence_0}, we may assume that the given existential $L_\theta(M)$-sentence $\phi$ is of the form
    $$
    \exists \ux \in \VV : \psi_\theta(\ux) \wedge S(\ux)
    $$
    where $S(\ux)$ is a $C$-sequence-system, and $\psi(\uxvec{})$ is an $L(M)$-formula that is bounded by $S$.
    We prove the theorem by induction on $(\rk(S), \deg(S))$.

    The base case $(\rk(S), \deg(S)) = (0, 0)$ is clear, as this implies that the tuple $\ux$ is empty.
    Thus, either $\psi(\uxvec{})$ implies no finite disjunction of non-trivial linear dependencies in $\uxvec{}$ over $\VV$, or it implies the empty disjunction.
    In the first case, there is nothing to show, and in the second case, $\phi$ is equivalent to the empty disjunction.

    For the induction step, assume that we have already shown the statement for any such sentence of the form $\exists \ux' \in \VV : \psi'_\theta(\ux') \wedge S'(\ux')$ with $(\rk(S'), \deg(S')) <_\Lex (\rk(S), \deg(S))$.
    If $\psi(\uxvec{})$ implies no finite disjunction of non-trivial linear dependencies in $\uxvec{}$ over $\VV$, there is nothing to show.
    So assume that $\psi(\uxvec{})$ implies a finite disjunction $\bigvee_{i=1}^q E_i(\uxvec{})$ of non-trivial linear dependencies in $\uxvec{}$ over $\VV$.
    It is clear that we can choose this disjunction such that every entry of $\uxvec{}$ that appears in some $E_i(\uxvec{})$ also appears in $\psi(\uxvec{})$.
    Since $\psi(\uxvec{})$ is bounded by $S$, this means that we can choose the disjunction such that every $E_i(\uxvec{})$ is bounded by $S$ as well.
    By Definition \ref{def_formual_bounded}, this also means that each $E_{i, \theta}(\ux)$ is bounded by $S$.
    Since $T$ is model-complete, $\psi(\uxvec{})$ implies $\bigvee_{i=1}^q E_i(\uxvec{})$ in any extension $\mm' \supseteq \mm$ with $\mm' \models T$.
    Hence, in any $C$-image-complete extension, we obtain
    $$
    \exists \ux \in \VV: \psi_\theta(\ux) \wedge S(\ux) \quad \equiv \quad \bigvee\nolimits_{i=1}^q\exists \ux \in \VV : \psi_\theta(\ux) \wedge S(\ux) \wedge E_{i, \theta}(\ux).
    $$
    Now we can apply Claim \ref{lemma_exist_sentence_1} to each disjunct $\exists \ux \in \VV : \psi_\theta(\ux) \wedge S(\ux) \wedge E_{i, \theta}(\ux)$ and conclude with the induction hypothesis.
\end{proof}

\noindent In a future paper, we will similarly use the tools from Section \ref{sec_trafo} to simplify any $L_\theta$-formula $\phi(\uz; \uw)$ modulo the model companion of $T^C_\theta$.
We now show how one direction of our characterization of existentially closed models follows from Theorem \ref{theorem_exist_sentence_reduction}:

\begin{proof}[Proof of ``$\Leftarrow$'' of Theorem \ref{theorem_big_characterization}]
    Assume $(\mm, \theta) \models T^C_\theta$ is $C$-image-complete, and
    $$
    (\mm, \theta) \models \exists \ux \in \VV : \psi_\theta(\ux) \wedge S(\ux)
    $$
    holds for every $C$-sequence-system $S(\ux)$ over $(\VV, \theta)$ and every $L(M)$-formula $\psi(\uxvec)$ that is bounded by $S$ and does not imply any finite disjunction of non-trivial linear dependencies in $\uxvec$ over $\VV$.
    By Theorem \ref{theorem_exist_sentence_reduction}, any existential $L_\theta(M)$-sentence $\phi$ either holds in the structure $(\mm, \theta)$ or is equivalent to the empty disjunction, i.e., to $\bot$, modulo the theory $T_\theta \cup \set{\text{``$\theta$ is $C$-image-complete''}} \cup \Diag(\mm, \theta)$.
    In the latter case, $\phi$ cannot hold in any extension $(\mm', \theta') \supseteq (\mm, \theta)$ with $(\mm', \theta') \models T^C_\theta$.
    Indeed, if it did, then $\phi$ would also hold in an existentially closed model of $T^C_\theta$ extending $(\mm', \theta')$.
    By Corollary \ref{corollary_C_image_compl}, this existentially closed model is $C$-image-complete, i.e., a model of $T_\theta \cup \set{\text{``$\theta$ is $C$-image-complete''}} \cup \Diag(\mm, \theta)$.
    This contradicts the fact that $\phi$ is equivalent to $\bot$ modulo this theory.
    Therefore $(\mm, \theta)$ must be an existentially closed model of $T^C_\theta$.
\end{proof}

\subsection{Proof of ``$\Rightarrow$''}

\noindent The goal of this section is to prove the direction from left to right of Theorem \ref{theorem_big_characterization}, but for generalized $C$-sequence-systems, which we define now.
A special form of these generalized $C$-sequence-systems will appear in a future paper, where we prove that if $T$ satisfies \Hfour{} and is \NATP{}, then the model companion of $T^C_\theta$ is also \NATP{} for any $C \in \Cc$.

\begin{definition} \label{def_param_gen_c_ss}
    Let $\ux := \ux_\li \ux_\ld := (x_{\lii, k} : 1 \leq k \leq m)(x_{\ldd, k} : 1 \leq k \leq n)$ and $(\VV, \theta) \models \TKvsThe$ be given.
    A \textbf{generalized $\mathbf{C}$-sequence-system} over $(\VV, \theta)$ is an $\LKThe(\VV)$-formula of the form
    \begin{align*}
        S(\ux) = & \bigwedge\nolimits_{k=1}^n \xi_k[\theta](x_{\ldd, k}) = \sum\nolimits_{l=1}^m P_{k, l}[\theta](x_{\lii, l}) + \sum\nolimits_{l=1}^{k-1} Q_{k, l}[\theta](x_{\ldd, l}) + u_{k}
    \end{align*}
    that satisfies the following conditions:
    \begin{enumerate}[(i)]
        \item If $C$ is algebraic, we require $m = 0$.
        \item Each $\xi_k$ is a monic polynomial that either satisfies $\Fac(\xi_k) \subseteq \Kp{C=\infty}$ or is of the form $f^q$ for some $f \in \Kp{0<C<\infty}$ and some $q$ with $0 < q \leq C(f)$.
        \item Given any $k$ with $\xi_k = f^q$ for some $f \in \Kp{0<C<\infty}$ and $q \in \NN$ with $0 < q \leq C(f)$, the following holds:
        $$
        \TKvsThe \cup \Diag(\VV, \theta) \models \forall \ux : S(\ux) \rightarrow f^C[\theta](x_{\ldd, k}) = 0.
        $$
    \end{enumerate}
\end{definition}
\noindent We will always denote generalized $C$-sequence-systems by the letter $S$.
Whenever a generalized $C$-sequence-system $S(\ux)$ over some $(\VV, \theta) \models \TKvsThe$ is given, we will assume that everything is as above, i.e., all $\xi_k$'s, $P_{k, l}$'s, and so on are defined implicitly, and we set $\ux = \ux_\li \ux_\ld$.
As with ``regular'' $C$-sequence-systems, i.e., the ones from Definition \ref{def_compatible}, we see that any generalized $C$-sequence-system over $(\VV, \theta)$ is also a generalized $C$-sequence-system over any $(\VV', \theta') \supseteq (\VV, \theta)$ with $(\VV', \theta') \models \TKvsThe$.
We now generalize boundedness from Definition \ref{def_formual_bounded} to generalized $C$-sequence-systems:

\begin{definition}[Boundedness] \label{bounded_gen_c_ss}
    Let $(\mm, \theta) \models T_\theta$ be given, and fix a generalized $C$-sequence-system $S(\ux) = \bigwedge\nolimits_{k=1}^n \xi_k[\theta](x_{\ldd, k}) = \sum\nolimits_{l=1}^m P_{k, l}[\theta](x_{\lii, l}) + \sum\nolimits_{l=1}^{k-1} Q_{k, l}[\theta](x_{\ldd, l}) + u_{k}$ over $(\VV, \theta)$, as in Definition \ref{def_param_gen_c_ss}.
    We say that an $L(M)$-formula $\psi(\uxvec)$ is \textbf{bounded} by $S$ if one of the following equivalent conditions holds:
    \begin{enumerate}[(i)]
        \item For all $k \in \set{1, \dots, n}$, the variable $x^i_{\ldd, k}$ does not appear in $\psi(\uxvec)$ for \hbox{$i \geq \deg(\xi_k)$}.
        \item For all $k \in \set{1, \dots, n}$, the term $\theta^i(x_{\ldd, k})$ does not appear in $\psi_\theta(\ux)$ for \hbox{$i \geq \deg(\xi_k)$}.
    \end{enumerate}
\end{definition}

\begin{observation} \label{rem_c_ss_is_also_gen}
    Let $(\mm, \theta) \models T_\theta$ be given, and let $S(\ux)$ be a ``regular'' $C$-sequence-system over $(\VV, \theta)$, i.e., as in (ii) of Definition \ref{def_compatible}.
    The following holds:
    \begin{enumerate}[(i)]
        \item $S(\ux)$ is also a generalized $C$-sequence-system, as in Definition \ref{def_param_gen_c_ss}.
        \item Any $L(M)$-formula $\psi(\uxvec{})$ is bounded by $S$ as a ``regular'' $C$-sequence-system, i.e., as in (i) of Definition \ref{def_formual_bounded}, if and only if $\psi(\uxvec{})$ is bounded by $S$ as a generalized $C$-sequence-system, i.e., as in Definition \ref{bounded_gen_c_ss} above.
    \end{enumerate}
\end{observation}

\noindent Now, with Observation \ref{rem_c_ss_is_also_gen} above, the following theorem implies, together with Corollary \ref{corollary_C_image_compl} (i.e., the fact that existentially closed models of $T^C_\theta$ are $C$-image-complete), the implication from left to right of our characterization of existentially closed models of $T^C_\theta$.

\begin{theorem} \label{theorem_stronger_left_to_right}
    Let $(\mm, \theta) \models T^C_\theta$ be existentially closed.
    Then
    $$
    (\mm, \theta) \models \exists \ux \in \VV : \psi_\theta(\ux) \wedge S(\ux)
    $$
    holds for every generalized $C$-sequence-system $S(\ux)$ over $(\VV, \theta)$ and every $L(M)$-formula $\psi(\uxvec)$ that is bounded by $S$ and implies no finite disjunction of non-trivial linear dependencies in $\uxvec$ over $\VV$.
\begin{proof}
    Assume that \hbox{$(\mm, \theta) \models T^C_\theta$} is existentially closed, and let $S(\ux)$ and $\psi(\uxvec)$ be given as in the statement above.
    Since $\psi_\theta(\ux) \wedge S(\ux)$ can be written as an existential formula (by the model completeness of $T$), and since $(\mm, \theta)$ is existentially closed, it suffices to find an extension $(\mm', \theta') \models T^C_\theta$ with $(\mm', \theta') \supseteq (\mm, \theta)$ and
    $$
    (\mm', \theta') \models \exists \ux \in \VV : \psi_\theta(\ux) \wedge S(\ux)
    $$
    in order to prove Theorem \ref{theorem_stronger_left_to_right}.
    We assume that $S(\ux)$ is of the following form
    \begin{align*}
        S(\ux) = & \bigwedge\nolimits_{k=1}^n \xi_k[\theta](x_{\ldd, k}) = \sum\nolimits_{l=1}^m P_{k, l}[\theta](x_{\lii, l}) + \sum\nolimits_{l=1}^{k-1} Q_{k, l}[\theta](x_{\ldd, l}) + u_{k}
    \end{align*}
    with everything as stated in Definition \ref{def_param_gen_c_ss}.
    Thus $\ux = \ux_\li\ux_\ld$, with $\ux_\li = (x_{\lii, 1}, \dots, x_{\lii, m})$ and $\ux_\ld = (x_{\ldd, 1}, \dots, x_{\ldd, n})$.
    We first simplify $S(\ux)$ a bit:

\begin{subclaim} \label{lemma_better_assumptions}
    We can, without loss, assume that $\deg(\xi_k) > 0$ and that $\deg(Q_{k, l}) < \deg(\xi_l)$ for all $k$ and $l < k$.
\begin{innerproof}
    Apply Lemma \ref{lemma_bound_term_light} with the conjunction
    $$
    \bigwedge\nolimits_{j = 1}^{k-1} \xi_{j}[\theta](x_{\ldd, j}) = \sum\nolimits_{l=1}^m P_{j, l}[\theta](x_{\lii, l}) + \sum\nolimits_{l=1}^{j-1} Q_{j, l}[\theta](x_{\ldd, l}) + u_{j}
    $$
    to the term
    $
    \sum\nolimits_{l=1}^m P_{k, l}[\theta](x_{\lii, l}) + \sum\nolimits_{l=1}^{k-1} Q_{k, l}[\theta](x_{\ldd, l}) + u_{k}
    $
    to ensure that $\deg(Q_{k, l}) < \deg(\xi_l)$ holds for all $k$ and $l < k$.

    Now fix $k$ and assume that $\xi_k = 1$.
    Then no variable $x^i_{\ldd, k}$ appears in $\psi(\uxvec)$, since this formula is bounded by $S$.
    Therefore, $x_{\ldd, k}$ does not appear in $\psi_\theta(\ux)$.
    Since we already have $\deg(Q_{l, k}) < \deg(\xi_k) = 0$, the variable $x_{\ldd, k}$ appears only in the equation $1[\theta](x_{\ldd, k}) = \dots$ in $S(\ux)$, and nowhere else in $\psi_\theta(\ux) \wedge S(\ux)$.
    By simple first-order logic, removing $1[\theta](x_{\ldd, k}) = \dots$ from $S(\ux)$ and deleting the entry $x_{\ldd, k}$ from $\ux$ does not change the truth value of the sentence $\exists \ux \in \VV : \psi_\theta(\ux) \wedge S(\ux)$.
    Hence, we can assume that $\deg(\xi_k) > 0$ for all $k$.
\end{innerproof}
\end{subclaim}

\noindent Since $\psi(\uxvec)$ implies no finite disjunction of non-trivial linear dependencies in $\uxvec$ over $\VV$, there is an elementary extension $\mm' \succ \mm$ and a realization
$\uvvec = \uvvec{}_\li \uvvec{}_\ld$ of $\psi(\uxvec)$ with entries in $\VV'$ that are linearly independent over $\VV$.
As with the placeholder variables, the arrow notation means that $\uvvec{}_\li = (v^i_{\lii, k} : i \in \omega, 1 \leq k \leq m)$ and $\uvvec{}_\ld = (v^i_{\ldd, k} : i \in \omega, 1 \leq k \leq n)$.
We also set $\uv^0 := \uv^0_\li\uv^0_\ld$, where $\uv^0_\li = (v^0_{\lii, 1}, \dots, v^0_{\lii, m})$ and $\uv^0_\ld = (v^0_{\ldd, 1}, \dots, v^0_{\ldd, n})$.

Our proof strategy is to construct a $C$-endomorphism $\theta'$ of $\VV'$ that extends $\theta$ in a way that $(\mm', \theta') \models \psi_\theta(\uv^0) \wedge S(\uv^0)$ holds.
We will first define $\theta'$ on a subspace of $\VV'$ and then extend it to all of $\VV'$ using our \standartConstruction{} (Lemma~\ref{lemma_standart_construction}).

\begin{subdefinition}
    Define $V$ as the following subtuple of $\uvvec$:
    $$
    (v_{\lii, k}^i : 1 \leq k \leq m,\; i \in \omega)^\frown (v_{\ldd, k}^i : 1 \leq k \leq n,\; 0 \leq i < \deg(\xi_k)).
    $$
    We extend $\theta$ to an endomorphism $\theta'$ of $\spanA{V \VV}{K}$ by setting $\theta'(v^i_{\lii, k}) := v^{i+1}_{\lii, k}$ and
    $$
    \theta'(v^i_{\ldd, k}) := \begin{cases}
        v^{i+1}_{\ldd, k} & \text{if $i < \deg(\xi_k)-1$} \\
        \sum\nolimits_{l=1}^m P_{k, l}[v_{\lii, l}] + \sum\nolimits_{l=1}^{k-1} Q_{k, l}[v_{\ldd, l}] + u_{k}-\sum\nolimits_{j=0}^{\deg(\xi_k)-1} (\xi_k)_j \cdot v^j_{\ldd, k} & \text{if $i = \deg(\xi_k)-1$} \\
    \end{cases}
    $$
    for all relevant $k$ and $i$, where, for example, $P_{k, l}[v_{\lii, l}] := \sum\nolimits_{j=0}^{\deg(P_{k, l})} (P_{k, l})_j \cdot v^j_{\lii, l}$, and similarly for the other polynomial expressions.
\end{subdefinition}

\noindent Note that everything above is well-defined.
Since $\uvvec$ is linearly independent over $\VV$, so is $V$.
Moreover, since $\deg(Q_{k, l}) < \deg(\xi_l)$ by Claim \ref{lemma_better_assumptions}, the terms $Q_{k, l}[v_{\ldd, l}]$ in the definition above lie in $\spanA{V}{K}$.

\begin{subclaim} \label{lemma_theta_behav_eq}
    The following holds:
    \begin{enumerate}[(i)]
        \item $\theta'^i(v^0_{\lii, k}) = v^i_{\lii, k}$.
        \item $\theta'^i(v^0_{\ldd, k}) = v^i_{\ldd, k}$ for $i \!<\! \deg(\xi_k)$,
        $\xi_k[\theta'](v^0_{\ldd, k}) \!=\! \sum\nolimits_{l=1}^m P_{k, l}[\theta'](v^0_{\lii, l}) + \sum\nolimits_{l=1}^{k-1} Q_{k, l}[\theta'](v^0_{\ldd, l}) + u_{k}.$
        \item Every $v \in \spanA{V \VV}{K}$ can be uniquely written as $v = \sum\nolimits_{k=1}^m \rho_k[\theta'](v^0_{\lii, k}) + \sum\nolimits_{k=1}^n \eta_k[\theta'](v^0_{\ldd, k}) + u$ with $\deg(\eta_k) < \deg(\xi_k)$ for all $k \in \set{1, \dots, n}$ and $u \in \VV$.
        \item Given $f \in \Kp{C<\infty}$ and any polynomial $\eta \in K[X]$, the equality
        $$
        (f^{C+1} \cdot \eta)[\theta'](v^0_{\ldd, k}) = r[\theta'](v^0_{\ldd, k}) + \sum\nolimits_{l=1}^m\! P'_{l}[\theta'](v^0_{\lii, l}) + \sum\nolimits_{l=1}^{k-1}\! Q'_{l}[\theta'](v^0_{\ldd, l}) + u
        $$
        holds, where $r \in K[X]$ is the remainder of $f^{C+1} \cdot \eta$ divided by $\xi_k$, the $P'_l$'s and $Q'_l$'s are polynomials with $\deg(Q'_{l}) < \deg(\xi_l)$ for $l < k$, and $u \in \VV$.
        In the case where $\eta = 0$ or $\xi_k = f^q$ for some $q$ with $0 < q \leq C(f)$, all polynomials on the right-hand side, as well as $u$, are $0$.
        \item For every $k$ with $\xi_k = f^q$ for some $f \in \Kp{0<C<\infty}$ and $q$ with $0 < q \leq C(f)$, we have $v^0_{\ldd, k} \in \Ker(f^C[\theta'])$.
        \item $(\spanA{V\VV}{K}, \theta') \models S(\uv^0)$.
        \item $\psi(\uvvec)$ is obtained from $\psi(\uxvec{})$ by replacing every occurrence of $x^i_{\lii, k}$ with $\theta'^i(v^0_{\lii, k})$ and every occurrence of $x^i_{\ldd, k}$ with $\theta'^i(v^0_{\ldd, k})$.
    \end{enumerate}
\begin{innerproof}
    Points (i), (ii), and (iii) are clear by construction.
    Point (vi) follows from (ii).
    Point (vii) follows from (i), (ii), and the fact that $\psi(\uxvec{})$ is bounded by $S$, i.e., that $x^i_{\ldd, k}$ does not occur in $\psi(\uxvec{})$ for $i \geq \deg(\xi_k)$.
    Point (v) follows from (vi), because $(\spanA{V\VV}{K}, \theta') \models S(\uv^0)$ implies $f^C[\theta'](v^0_{\ldd, k}) = 0$ for all such $k$ by (iii) of Definition \ref{def_param_gen_c_ss}.
    We now show (iv).
    The \remainderRule{} (Lemma~\ref{lemma_remainder_rule}) applied to the equation $\xi_k[\theta'](v^0_{\ldd, k}) = \dots$ from (ii) yields
    \begin{align*}
        (f^{C+1} \cdot \eta)[\theta'](v^0_{\ldd, k}) = r[\theta'](v^0_{\ldd, k}) +\chi[\theta']\Big(\sum\nolimits_{l=1}^m\!\! P_{k, l}[\theta'](v^0_{\lii, l}) +\! \sum\nolimits_{l=1}^{k-1}\!\! Q_{k, l}[\theta'](v^0_{\ldd, l}) + u_{k}\Big).
    \end{align*}
    Now apply Lemma \ref{lemma_bound_term_light} with the conjunction $\bigwedge_{l=1}^{k-1} \xi_l[\theta'](v^0_{\ldd, l}) = \dots$ (again from (ii)) to the term $\chi[\theta'](\dots)$.
    The final assertion in the case $\eta = 0$ or $\xi_k = f^q$ follows from the linear independence of $V$ (and from (v) in the case $\xi_k = f^q$).
\end{innerproof}
\end{subclaim}

\noindent With (vi) and (vii) of Claim \ref{lemma_theta_behav_eq}, we immediately see that $(\mm', \theta') \models \psi_\theta(\uv^0) \wedge S(\uv^0)$ holds once we extend $\theta'$ to an endomorphism defined on all of $\VV'$.

\begin{subclaim}
    The map $\theta'$ is a $C$-endomorphism of $\spanA{V\VV}{K}$.
\begin{innerproof}
    First assume that $C$ is a transcendental kernel configuration.
    In this case, we need to verify that $\Ker(f^C[\theta']) = \Ker(f^{C+1}[\theta'])$ holds for all $f \in \Kp{C<\infty}$.
    Fix $f \in \Kp{C<\infty}$ and $v \in \Ker(f^{C+1}[\theta'])$.
    For ease of notation, set
    $$
    \kk_f := \set{k \in \set{1, \dots, n} : \xi_k = f^q \text{ for some $q$ with $0 < q \leq C(f)$}}.
    $$
    By (iii) of Claim \ref{lemma_theta_behav_eq}, we can write
    $$
    v = \sum\nolimits_{k=1}^m \rho_k[\theta'](v^0_{\lii, k}) + \sum\nolimits_{k=1}^n \eta_k[\theta'](v^0_{\ldd, k}) + u
    $$
    where all the $\rho_k$'s and $\eta_k$'s are polynomials and $u \in \VV$.
    Also, we have $\deg(\eta_k) < \deg(\xi_k)$ for all $k \in \set{1, \dots, n}$.
    By (iv) of Claim \ref{lemma_theta_behav_eq}, for each $k \in \set{1, \dots, n}$ we obtain
    $$
        (f^{C+1} \cdot \eta_k)[\theta'](v^0_{\ldd, k}) = r_{k}[\theta'](v^0_{\ldd, k}) + \sum\nolimits_{l=1}^m\! P'_{k, l}[\theta'](v^0_{\lii, l}) + \sum\nolimits_{l=1}^{k-1}\! Q'_{k, l}[\theta'](v^0_{\ldd, l}) + u'_{k}
    $$
    where $r_k$ is the remainder of $f^{C+1} \cdot \eta_k$ divided by $\xi_k$, each $P'_{k, l}$ is a polynomial, each $Q'_{k, l}$ is a polynomial with $\deg(Q'_{k,l}) < \deg(\xi_l)$, and $u'_k \in \VV$.
    Moreover, if $\eta_k = 0$ or $k \in \kk_f$, then $r_k$, all $P'_{k,l}$'s, all $Q'_{k, l}$'s, and $u'_k$ are $0$.
    We now obtain:
    \begin{align}
        0 = f^{C+1}[\theta'](v) &= \sum\nolimits_{k=1}^m \Big(\rho_k \cdot f^{C+1} + \sum\nolimits_{l=1}^n P'_{l, k}\Big)[\theta'](v^0_{\lii, k})  \label{tag_proofCend_1}\\
        &\quad+ \sum\nolimits_{k=1}^n \Big(r_{k} + \sum\nolimits_{l=k+1}^n Q'_{l, k}\Big)[\theta'](v^0_{\ldd, k})  \label{tag_proofCend_2}\\
        &\quad+ \sum\nolimits_{k=1}^n u'_{k} + f^{C+1}[\theta](u).  \label{tag_proofCend_4}
    \end{align}
    Notice that the $k$-th summand of term (\ref{tag_proofCend_1}) lies in $\spanA{v_{\lii, k}^i : i \in \omega}{K}$, the $k$-th summand of term (\ref{tag_proofCend_2}) lies in $\spanA{v_{\ldd, k}^i : 0 \!\leq\! i \!<\! \deg(\xi_k)}{K}$, and term (\ref{tag_proofCend_4}) lies in $\VV$.
    Since $\spanA{V\VV}{K}$ is the direct sum of all these subspaces, we conclude that all summands of (\ref{tag_proofCend_1}), all summands of (\ref{tag_proofCend_2}), and term (\ref{tag_proofCend_4}) must be $0$.
    This has the following consequences:
    \begin{enumerate}[(i)]
        \item $\eta_k = 0$ for all $k \not\in \kk_f$.
        Suppose not, and choose $k$ maximal with $\eta_k \neq 0$ and $k \not\in \kk_f$.
        Recall that $\deg(\eta_k) < \deg(\xi_k)$, and hence $\xi_k \nmid \eta_k$.
        Also recall that $\xi_k$ either has only irreducible factors in $\Kp{C=\infty}$ or is of the form $g^q$ for some $g \in \Kp{0<C<\infty} \setminus \set{f}$.
        Since $C(f) < \infty$ and $\deg(\xi_k) > 0$, we obtain $\xi_k \nmid f^{C+1} \cdot \eta_k$, and therefore $r_{k} \neq 0$.
        Since $\eta_l = 0$ or $l \in \kk_f$ for every $l > k$, (iv) of Claim \ref{lemma_theta_behav_eq} gives $Q'_{l, k} = 0$ for every $l > k$.
        Hence, the $k$-th summand of term (\ref{tag_proofCend_2}) is just $r_{k}[\theta'](v^0_{\ldd, k})$.
        Using $r_{k} \neq 0$, the inequality $\deg(r_{k}) < \deg(\xi_k)$, (ii) of Claim \ref{lemma_theta_behav_eq} (i.e., $\theta'^i(v^0_{\ldd, k}) = v^i_{\ldd, k}$ for $i < \deg(\xi_k)$), and the linear independence of $V$, we contradict that the $k$-th summand of term (\ref{tag_proofCend_2}) is $0$.
        \item $\rho_k = 0$ for all $k$.
        By (i), we know that $\eta_k = 0$ or $k \in \kk_f$ holds for all $k$.
        Therefore, every $P'_{l, k}$ is also $0$ by (iv) of Claim \ref{lemma_theta_behav_eq}, which implies that the $k$-th summand of term (\ref{tag_proofCend_1}) is equal to
        $
        (\rho_k \cdot f^{C+1})[\theta'](v^0_{\lii, k})
        $.
        As in (i), but with (i) of Claim \ref{lemma_theta_behav_eq}, this implies $\rho_k = 0$.
    \end{enumerate}
    This shows that
    $
    v = \sum\nolimits_{k \in \kk_f} \eta_k[\theta'](v^0_{\ldd, k}) + u.
    $
    By (v) of Claim \ref{lemma_theta_behav_eq}, we obtain
    $$
    0 = f^{C+1}[\theta'](v) = f^{C+1}[\theta](u),
    $$
    and $f^{C}[\theta'](v) =\! f^{C}[\theta](u)$.
    Since $\theta$ is a $C$-endomorphism, which implies $\Ker(f^C) = \Ker(f^{C+1})$, the first equation yields $f^C[\theta](u) = 0$, and therefore $f^C[\theta'](v) = 0$.
    This completes the proof in the transcendental case.

    Assume that $C$ is an algebraic kernel configuration, and let $v \in \spanA{V\VV}{K}$ be given.
    We must show that $\Ker(\mipo(C)[\theta']) = \spanA{V\VV}{K}$.
    By (i) of Definition \ref{def_param_gen_c_ss}, the tuple $\ux_\li$ is empty if $C$ is algebraic.
    Thus (iii) of Claim \ref{lemma_theta_behav_eq} gives
    $$
    v = \sum\nolimits_{k=1}^{n} \eta_k[\theta'](v^0_{\ldd, k}) + u.
    $$
    Also note that $\Kp{C=\infty}$ is empty, since $C$ is algebraic; hence every $\xi_k$ is of the form $f^q$ for some $f \in \Kp{0<C<\infty}$ and some $q$ with $0 < q \leq C(f)$.
    Together with (v) of Claim \ref{lemma_theta_behav_eq}, this yields $v^0_{\ldd, k} \in \Ker(f^C[\theta']) \subseteq \Ker(\mipo(C)[\theta'])$ for some $f \in \Kp{0<C<\infty}$ that depends on $k$.
    Since $\theta$ is a $C$-endomorphism, we also have $\mipo(C)[\theta'](u) = \mipo(C)[\theta](u) = 0$.
    Combining everything, we obtain $\mipo(C)[\theta'](v) = 0$.
\end{innerproof}
\end{subclaim}

\noindent Recall that $\spanA{V\VV}{K}$ is a subspace of the vector space $\VV'$ of an extension $\mm' \succ \mm$.
Without loss of generality, we can assume that $\dim(\VV' / \spanA{V\VV}{K}) \geq \aleph_0$.
Using our \standartConstruction{}, we can now extend $\theta'$ to a $C$-endomorphism of $\VV'$.
Then $(\mm', \theta') \models T^C_\theta$.
Moreover, (vi) and (vii) in Claim \ref{lemma_theta_behav_eq} still hold for this extended $\theta'$.
Therefore, $(\mm', \theta') \models \psi_\theta(\uv^0) \wedge S(\uv^0)$.
Thus $(\mm', \theta') \models T^C_\theta$ is an extension of $(\mm, \theta)$ with $(\mm', \theta') \models \exists \ux \in \VV : \psi_\theta(\ux) \wedge S(\ux)$.
\end{proof}
\end{theorem}

\noindent As mentioned above Theorem \ref{theorem_stronger_left_to_right} this finishes the proof of our characterization of existentially closed models of $T^C_\theta$, i.e., Theorem \ref{theorem_big_characterization}.

\section*{Acknowledgments}
\noindent The author would like to thank Christian d'Elbée for reading earlier drafts of this paper and for providing valuable feedback and suggestions, especially on the introduction.
The author is also grateful to Mike Prest for his comments and suggestions, which led to significant improvements to this paper, especially regarding the characterization of existentially closed models and its proof.
The author would also like to thank Philipp Hieronymi for general support.

\bibliographystyle{alpha} 
\bibliography{refs}

\end{document}